%% file: main.tex
\documentclass[12pt]{article}{}
\pdfoutput=1

\newcommand{\doctitle}{Reduced order models from computed states of physical systems using non-local calculus on finite weighted graphs}
\newcommand{\docauthor}{M. Duschenes \& K. Garikipati}
\newcommand{\docaffil}{Applied Physics, Mathematics, Mechanical Engineering, and Michigan Institute for Computational Discovery \& Engineering\\ University of Michigan}

\newcommand{\docheader}{Reduced order models ... graphs - \docauthor}
\newcommand{\docfooter}{University of Michigan}

\usepackage{main}

\renewcommand{\vec}[1]{\mathbf{#1}}
\newcommand{\vecbar}[1]{\vec{\bar{#1}}}

\newcommand{\nmodm}[2]{#1\ (\mathrm{mod}\ #2)}
\newcommand{\sign}[1][1=]{\text{sgn}(#1)}

\renewcommandx{\square}{{\ooalign{$\sqsubset\mkern2mu$\cr$\mkern1mu\sqsupset$\cr}}}
\newcommandx{\hrectangle}{{\ooalign{$\sqsubset\mkern1mu$\cr$\mkern3mu\sqsupset$\cr}}+}
\newcommandx{\vrectangle}{{\ooalign{$\sqsubset\mkern1mu$\cr$\mkern3mu\sqsupset$\cr}}-}
\renewcommandx{\dotproduct}[2][1=,2=]{#1\cdot#2}
\renewcommandx{\gradient}[2][1={w},2=]{{\nabla}_{#1}#2}
\renewcommandx{\laplacian}[2][1={w},2=]{\gradient[#1]^2#2}
\newcommandx{\nlaplacian}[3][1={w},2=,3=n]{\gradient[#1][]^{#3}#2}
\renewcommandx{\divergence}[2][1={w},2=]{\dotproduct[{\gradient[#1]}][{#2}]}
\newcommand{\varchi}{\raisebox{\depth}{$\chi$}}

\newcounter{Nargs}
\def\numargs#1{%
	\setcounter{Nargs}{0}%
	\foreach \x in #1{%
	\stepcounter{Nargs}%
 }%
	\the\value{Nargs}%
} 

\renewcommand*{\derivative}[3][0]{%
	\def\N{\numargs{{#3}}}%
	\frac{\partial^{\ifnum#1=1 {} \else #1 \fi} #2}{%
	\foreach \x in {#3} {\partial \x}}
}

\newcommand*{\uniderivative}[3][1]{%
	\frac{\partial^{\ifnum#1=1 {} \else #1 \fi} #2}{
	\partial {{}{#3}}^{\ifnum#1=1 {} \else #1 \fi}}
}

\newcommand*{\uninderivative}[3][n]{%
	\frac{\partial^{#1} #2}{
	\partial {{}{#3}}^{#1}}
}

\newcommand*{\nderivative}[4][n]{%
	\frac{\partial^{#1} #2}{\partial#3\cdots\partial#4}
}

\newcommand*{\difference}[3][0]{%
	\def\M{1}%
	\def\N{\numargs{{#3}}}%
	\frac{\delta^{\ifnum#1=1 {} \else #1 \fi} #2}{%
	\foreach \x in {#3} {\delta \x}}
}

\newcommand*{\unidifference}[3][1]{%
	\frac{\delta^{\ifnum#1=1 {} \else #1 \fi} #2}{
	\delta {{}#3}^{\ifnum#1=1 {} \else #1 \fi}}
}

\newcommand*{\unindifference}[3][n]{%
	\frac{\delta^{#1} #2}{
	\delta {{}{#3}}^{#1}}
}

\newcommand*{\ndifference}[4][n]{%
	\def\N{\numargs{{#3}}}%
	\frac{\delta^{#1} #2}{\delta#3\cdots\delta#4}
}

\newcommand*{\Difference}[3][0]{%
	\def\M{1}%
	\def\N{\numargs{{#3}}}%
	\frac{\Delta^{\ifnum#1=1 {} \else #1 \fi} #2}{%
	\foreach \x in {#3} {\Delta \x}}
}

\newcommand*{\uniDifference}[3][1]{%
	\frac{\Delta^{\ifnum#1=1 {} \else #1 \fi} #2}{
	\Delta {{}#3}^{\ifnum#1=1 {} \else #1 \fi}}
}

\newcommand*{\derivativeinline}[3][0]{%
	\def\N{\numargs{{#3}}}%
	{\partial^{\ifnum#1=1 {} \else #1 \fi} #2}/{%
	\foreach \x in {#3} {\partial \x}}
}

\newcommand*{\uniderivativeinline}[3][1]{%
	{\partial^{\ifnum#1=1 {} \else #1 \fi} #2}/{
	\partial {{}#3}^{\ifnum#1=1 {} \else #1 \fi}}
}

\newcommand*{\nderivativeinline}[4][n]{%
	\def\N{\numargs{{#3}}}%
	{\partial^{#1} #2}/{\partial#3\cdots\partial#4}
}

\newcommand*{\differenceinline}[3][0]{%
	\def\M{1}%
	\def\N{\numargs{{#3}}}%
	{\delta^{\ifnum#1=1 {} \else #1 \fi} #2}/{%
	\foreach \x in {#3} {\delta \x}}
}

\newcommand*{\unidifferenceinline}[3][1]{%
	{\delta^{\ifnum#1=1 {} \else #1 \fi} #2}/{
	\delta {{}#3}^{\ifnum#1=1 {} \else #1 \fi}}
}

\newcommand*{\Differenceinline}[3][0]{%
	\def\M{1}%
	\def\N{\numargs{{#3}}}%
	{\Delta^{\ifnum#1=1 {} \else #1 \fi} #2}/{%
	\foreach \x in {#3} {\Delta \x}}
}

\newcommand*{\uniDifferenceinline}[3][1]{%
	{\Delta^{\ifnum#1=1 {} \else #1 \fi} #2}/{
	\Delta {#3}^{\ifnum#1=1 {} \else #1 \fi}}
}

\newcommand{\etal}{\emph{et al.}~}

\begin{document}

\maketitle
\begin{abstract}
\pagestyle{abstract}
\noindent Partial differential equation-based numerical solution frameworks for initial and boundary value problems have attained a high degree of complexity. Applied to a wide range of physics with the ultimate goal of enabling engineering solutions, these approaches encompass a spectrum of spatiotemporal discretization techniques that leverage solver technology and high performance computing. While high-fidelity solutions can be achieved using these approaches, they come at a high computational expense and complexity. Systems with billions of solution unknowns are now routine. The expense and complexity do not lend themselves to typical engineering design and decision-making, which must instead rely on reduced-order models. Here we present an approach to reduced-order modelling that builds off of recent graph theoretic work for representation, exploration, and analysis on computed states of physical systems (Banerjee \etal, \emph{Comp. Meth. App. Mech. Eng.}, \textbf{351}, 501-530, 2019). We extend a non-local calculus on finite weighted graphs to build such models by exploiting first order dynamics, polynomial expansions, and Taylor series. Some aspects of the non-local calculus related to consistency of the models are explored. Details on the numerical implementations and the software library that has been developed for non-local calculus on graphs are described. Finally, we present examples of applications to various quantities of interest in mechano-chemical systems.
\end{abstract}

\newpage
\clearpage
\singlespacing
\thispagestyle{plain}
\tableofcontents
\thispagestyle{plain}

\newpage
\clearpage
\singlespacing
\pagenumbering{arabic}
\pagestyle{document}


\newpage

\section{Introduction} \label{sec:intro}
Natural phenomena that are susceptible to a mathematical physics treatment typically lead to initial and boundary value problems. The solution of partial differential equations that govern these problems has occupied the better part of the last century of research in applied mathematics, leading to numerical methods, and from there to computational science and high performance computing. With Fourier methods, time integration algorithms, discretization frameworks of finite elements, finite difference, Green's functions-based formulations, and hybrids of these techniques among other approaches, it is now possible to attain a degree of fidelity that is sufficient for many applications in mathematical physics. While the need to incorporate additional physics and gain a finer resolution of phenomena will never be met, there is also the recognition that the fidelity that is already achievable comes paired with a complexity that outstrips easy representation and analysis. Reduced-order models must therefore be derived on a principled basis with controlled approximations.

For the purpose of distinguishing the approach adopted here, we take account of areas including projection methods, coarse graining techniques, dimensionality reduction, sub-manifold extraction and numerous others that are arguably encompassed by the field of dynamical systems. Within this field, is the graph theoretic approach for representation, exploration and analysis of computed states of physical systems that we have proposed recently \cite{Banerjee2019}. This formalism is predicated on computed solutions that are typically high-dimensional in the sense that the primary unknown fields are represented by degrees of freedom numbering $\gtrsim O(10^8)$. Functionals on these discrete solutions are extracted as \emph{states}. Some examples are the lift, drag or pressure difference in computational fluid dynamics, measures of load, deformation and total energies in solid mechanics, phase volume fractions and free energies in materials physics, and cross sections in nuclear physics \cite{Kochunas2020}. These states undergo \emph{transitions} as some \emph{parameter} of the system is varied, such as its physical, possibly time-like, parameters, boundary or initial conditions. 

\subsection{Quantities of interest and graph representations} \label{sec:intro_qoi}
The states, of which there is typically a small number $\sim O(10)$ relevant to a system, form a low-dimensional vector of quantities of interest. However, these state vectors can be of high fidelity limited only by the dimensionality of the full-field numerical solution. For completeness of the description, the state vector can be extended to include the parameters referred to above. We have shown that by defining the states to be vertices on a graph, and the transitions to be edges between the vertices, a near one-to-one correspondence is uncovered between the properties of the computed physical system and the elements of graph theory. Weighted edges between vertices can then be assigned, or found through graph theoretic principles, indicating relationships or transitions between states, and the magnitude of such correlations. Reversible linear systems lead to undirected graphs that are fully connected, whereas dissipative dynamical systems are represented by directed trees. Path-dependence can be thought of as being associated with an absence of cycles. Notions of centrality reveal insights to the relations between states and path traversal properties are induced on graphs by equilibrium or dynamical transitions on the physical system.

Given the above graph theoretic framework we seek to extend it to reduced order modelling based on computed states and transitions of physical systems. We take the view of the graph of states (vertices) and transitions (edges) as a discrete manifold describing the physical system. On this manifold we seek to define differential equations and functional representations for the evolution of certain components of the state vector. An illustrative example that we pursue in this communication is of an ordinary differential equation for the phase fractions of a material system, driven by itself and other components of the state vector such as the free energy and strains. Guided by our knowledge of first-order systems in materials physics, it is of interest to compute derivatives on these states, such as of the free energy, with respect to the phase fractions, yielding generalized chemical potentials, or with respect to the strains, yielding stresses. The motivation comes from the rates of change of phase fractions driven by chemical potentials and stresses in materials. An alternative to the differential equation is a direct functional representation, such as for the free energy with phase fractions and strains as arguments. In this case, we are guided by the universality of the Taylor series representation, to which we seek approximations. This too makes necessary the computation of derivatives of state vector components with respect to each other. 

\subsection{Non-local calculus for reduced order modelling on graphs} \label{sec:intro_nonlocalcalculus}
We have adopted a discrete non-local calculus on weighted graphs to develop the reduced order models outlined above. This formalism, which has been applied previously to image analysis \cite{Gilboa2008,Desquesnes2013,Hein2007} begins with definition of a vector space on the graph, followed by the gradient operator, inner products and partial derivatives. Also important is the consistency and rigor reflected in a divergence operator, adjointness relations and a divergence theorem. This non-local theory bears a relation to classical differential calculus, that is determined by the graph weights. Through its self consistently defined inner products, the weights appear as kernels determining the support of integral operators on the graph manifold, and crucially controlling the behavior of partial derivatives relative to their classical definitions. 

This forms a central consideration on theoretical and practical grounds. However, the machinery that translates the non-local calculus on graphs to reduced order models can be formally put in place independently of the correspondence between non-local and local treatments. With functions and non-local partial derivatives in hand, algebraic and non-local differential operators can be introduced and viewed as basis terms in the reduced order representation. Since these can be evaluated at the vertex data points, a stepwise regression approach can be invoked to choose the best combination of operators representing the evolution or dependence sought. Model parsimony can be dialed in by using thresholds on the stepwise regression.

\subsection{Physical systems of interest} \label{sec:intro_physicalsystems}
The system of interest for this work are multi-crystalline solids, where the dynamical quantities of interest are the volume fraction of each phase, which are assumed to undergo first-order dynamics. These dynamics are selected to be forced by a linear functional of a $\ird[3]$ order polynomial of the phase volume fractions themselves, the interfacial lengths between the phases, the number of particles, and the first derivatives of the free energy, with respect to the phase volumes, and the mechanical strains. In order to complete the model of the system to be able to forward solve the first order dynamics, a functional representation for the free energy is found as a Taylor series in the phase volume fraction and mechanical strains, up to $\ith[4]$ order. This example allows the general graph theoretic method to be shown, and analyzed for its effectiveness at developing functional representations, and investigating any numerical issues associated with these derivatives. Further examples of physical systems will be explored in future communications.

Computations on each of these examples relies upon a general graph theory library code that accepts general data with inputs and outputs, and determines appropriate models based on these stepwise regression procedures, as well as visualizes fits with plots of loss curves, fitted curves for various models, and trends in the relative magnitudes of fitted coefficients over the stepwise regression. These examples and library allow conjectures on future applications and examples of this graph theoretic approach to be explored. 

The non-local calculus on finite weighted graphs is laid out in Section \ref{sec:graphtheory}, illustrative physical systems and discussions are presented in Section \ref{sec:physicalsystems}, and conclusions appear in Section \ref{sec:concl}. Appendices outline the implementation details and example calculations for these methods, including details of stepwise linear regression in \cref{app:implementation_stepwiseregression,app:implementation_linearregression}, implementation of the graph theory code library \cref{app:implementation_graphtheorylibrary,app:implementation_gauss}, and error analyses of the method in \cref{app:error,app:error_p}.


\section{Graph theoretic approach} \label{sec:graphtheory}
The major formalism used in this work relies on methods from graph theory \cite{Banerjee2019,West2001,Newman2010}. A graph, denoted by $G = (V,E)$, consists of a set of vertices, $i \in V$, with $\abs{V} = n$. The vertices are connected by a set of edges, $e \in E$, where each edge $e = (i,j)$ is a pair of vertices. A vertex is denoted as being the part of an edge by $i \in e$. The graph can also have global and local attributes on the vertices and edges, as we describe below. Of importance are edge weights $w$, that are generally functions of the local vertex attributes, and admit a notion of distance between states on the graph as will be elaborated on in \cref{sec:graphtheory_weightconstraints}.

As outlined in \cref{sec:intro_qoi} of the Introduction, in our treatment, a state is a low-dimensional vector of quantities of dimension $p$, obtained after solving a high-fidelity computation of dimension $P$, where ideally $p \ll P$. These high fidelity computations are typically direct numerical simulations (DNS), for example the solution of partial differential equations (PDE) using a finite element-like approach. The DNS simulations can yield some functionals, such as volume averaged fields, or resulting scalar outputs. These, collected into a low-dimensional vector will be denoted by $u$. Separately, \emph{parameters} such as physical parameters, or initial and boundary conditions of the DNS can be collected into another low-dimensional vector denoted by $x \in \reals[p]$. Each state is represented as a vertex $i \in V$ on the graph. It follows that the notation $x_{i}$ is unambiguous in labelling a state by its vertex. Transitions between these states, including changes in system parameters, or steps of the numerical solver, allow for connections to be made between vertices on the graph. These transitions therefore can be thought of as inducing edges $e = (i,j)$ between states $i$ and $j$.

Stated formally now, the vertices, $i \in V$ represent states $\{x^{\mu}_{i},u(x_{i})\}$, where, as introduced above, $x_{i}$ are parameters (physical parameters, initial and boundary conditions) and $u(x_{i})$ are quantities of interest, obtained as functionals from the high-dimensional numerical solution of a system of PDEs. The graphs naturally inherit the mathematical structure of the PDEs \cite{West2001,Newman2010}, and even confer new machinery to it \cite{Banerjee2019}.

\begin{figure}[ht]
	\centering
	\includegraphics[width=0.6\textwidth,trim={0 0 0 6cm},clip]{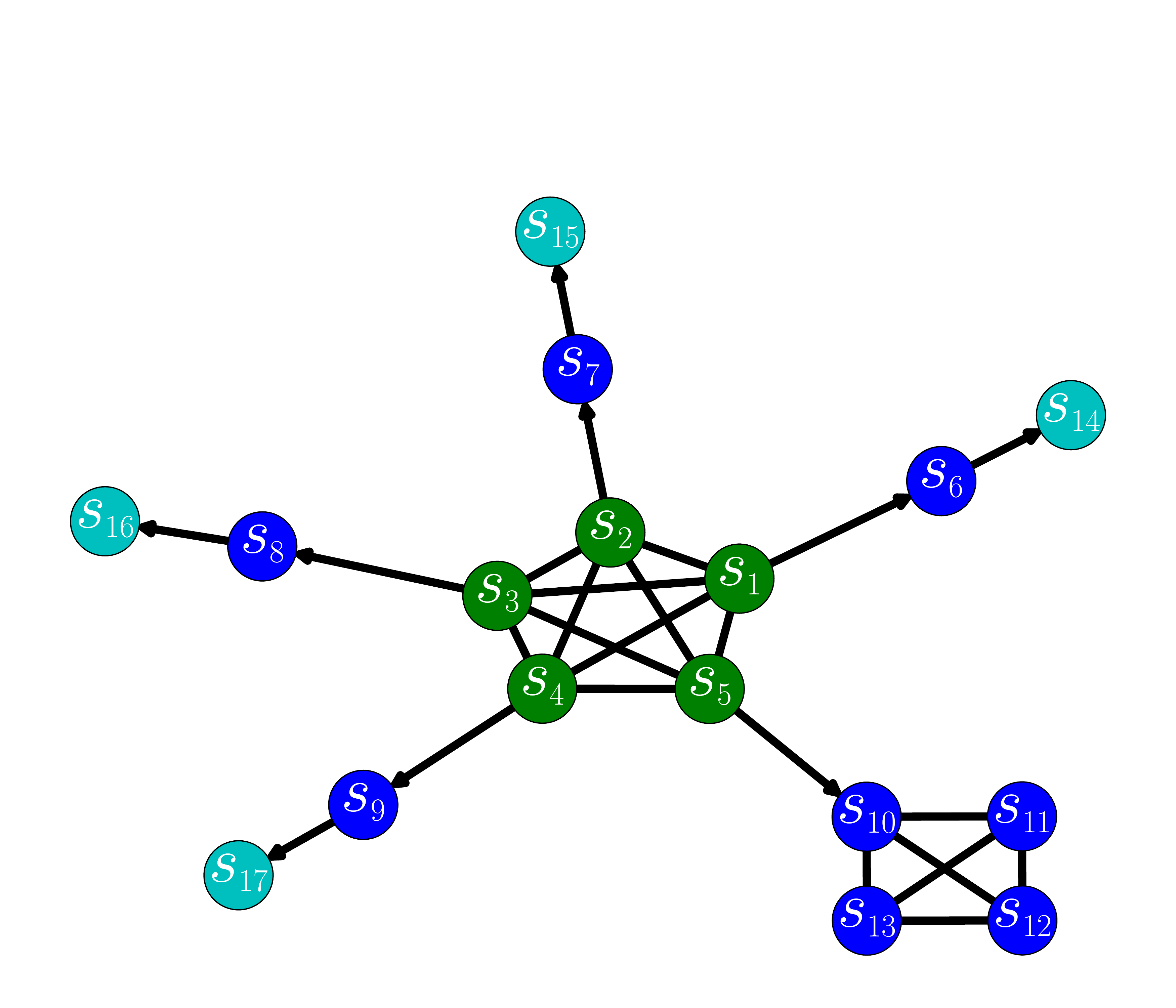}
	\caption{Example graph, with a central clique, and branching clique subgraphs. Vertex coloring may represent various local attributes.}
	\label{fig:cliquegraph}
\end{figure}

\begin{figure}[ht]
	\centering
	\includegraphics[width=0.6\textwidth,trim={0 0 0 6cm},clip]{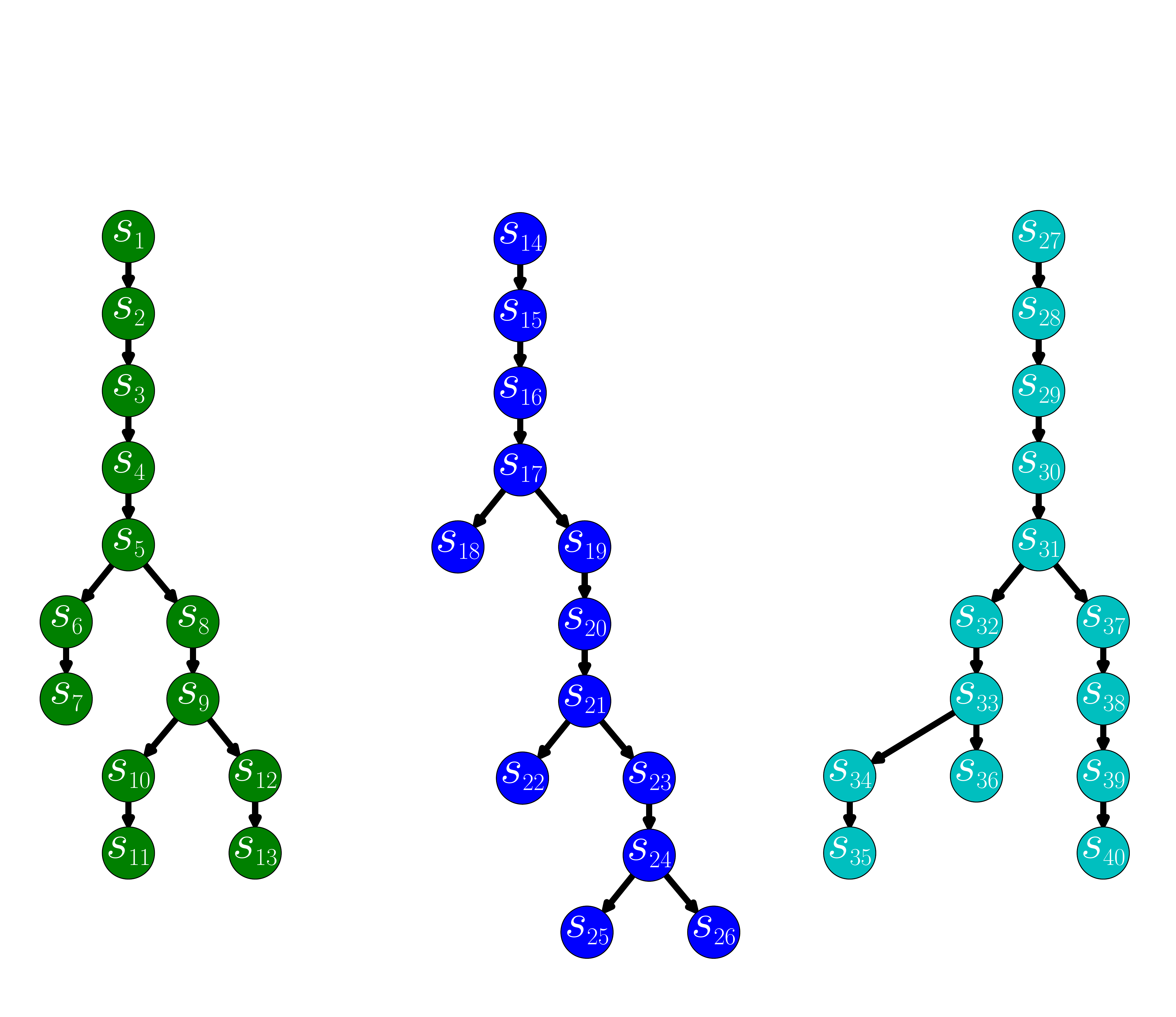}
	\caption{Example directed tree graphs.}
	\label{fig:treegraph}
\end{figure}

Several reoccurring graph structures arise in formal graph theory, including fully connected clique graphs, such as the subgraphs shown in \cref{fig:cliquegraph}, or minimally connected tree graphs, such as shown in \cref{fig:treegraph}. There are intricate correspondences between the physical processes being studied and the resulting graph structure. As identified by Banerjee \etal\cite{Banerjee2019}: directed graphs and less densely connected tree graphs correspond to processes with causality or irreversibility induced by energy dissipation, whereas undirected graphs correspond to reversible processes. Fully connected clique graphs generally represent linear processes in equilibrium. These correspondences, centrality measures and graph traversal will be revisited in future communications.

\subsection{Non-local calculus on finite, weighted graphs} \label{sec:graphtheory_nonlocalcalculus}
As outlined in \cref{sec:intro_nonlocalcalculus} of the Introduction, we seek to develop reduced order models in the form of differential equations for the evolution of some component of the state vector on the graph, or as functional representations. In the latter case, we explore forms that are motivated by the Taylor Series of local, differential calculus. In both cases, we need to first define a calculus on the discrete manifold that is the graph. Specifically, we adopt a discrete, non-local calculus on finite, weighted graphs \cite{Hein2007,Gilboa2008,Elmoataz2008,Desquesnes2013,Lozes2014}.

Gilboa \etal define a discrete calculus, consisting of non-local operators, based on differences between states $i,j \in V$, and an edge weight $w(x_{i},x_{j})$\cite{Gilboa2008}. The edge weights $w$ can be defined in a number of ways, and $w(x_{i},x_{j}) = w(x_{j} - x_{i})$ is natural for our purposes. In the case of symmetric, undirected graphs, a radially symmetric weight, $w(x_{i},x_{j}) = w(r(x_{i},x_{j}))$, where $r(x_{i},x_{j}) = \abs{x_{j} - x_{i}}$ is also appropriate. As will be discussed, there are several constraints placed on the allowed forms of the weight functions, to ensure that the non-local calculus converges to true differential operators in certain limits. 

Scalars, $u(x_{i})$ correspond to a function at a single vertex $i$. Vectors $v(x_{i},x_{j})$ are functions on vertex pairs $i,j$ or the edge $e = (i,j)$. Accordingly, we have a vector space $\mathnormal{L}$ such that there is a mapping $V\times V \mapsto \mathnormal{L}$ or $E\mapsto \mathnormal{L}$. Integration on the discrete manifold $G$ is a sum over vertices $i \in V$, and is relevant to operations such as inner products, which we define to be normalized by the size of the graph, $\abs{V} = n$. There should be no cause for confusion between the state vectors $x_{i} \in \reals[p]$ represented as vertices $i \in V$, and the edge vectors $v(x_{i},x_{j}) \in \mathnormal{L}$. These are distinct manifestations of vector spaces that are not copies of each other. We also note that for this treatment, there are no self edges and so vectors $v(x_{i},x_{j})$ at vertex $i$ are defined as being between vertices distinct from $i$, and $j \neq i$.

\subsubsection{Non-local calculus definitions} \label{sec:graphtheory_nonlocalcalculus_definitions}
The non-local gradient operator $\gradient[][{u}] = \gradient[{w}][[{u}]](x_{i},x_{j}) $, with respect to the weight $w$, can now be defined as the vector functional of scalars $u$ at state $x_{i}$, and represents the vector of weighted differences with \textit{all other vertices} $j \in V$:
\begin{equation}
	 \gradient[{w}][[{u}]](x_{i},x_{j}) \equiv [u(x_{j})-u(x_{i})]\sqrt{w(x_{i},x_{j})},\quad \forall ~ j \neq i \in V.
\end{equation}
The inner product between scalars $u_{\alpha}$ and $u_{\beta}$ is defined as the sum over the vertices:
\begin{align}
	\langle u_{\alpha},u_{\beta}\rangle \equiv&~ \frac{1}{n}\sum_{i\in V} u_{\alpha}(x_{i})u_{\beta}(x_{i}),
\end{align}
and the contraction between vectors $v_{\alpha}$ and $v_{\beta}$ at vertex $i$ is defined as the sum over disparate vertices:
\begin{align}
	[\dotproduct[{v_{\alpha}}][{v_{\beta}}]](x_{i}) \equiv&~ \frac{1}{n-1}\sum_{j \in V\setminus \{i\}} v_{\alpha}(x_{i},x_{j})v_{\beta}(x_{i},x_{j}).
\end{align}
Inner products, and norms between pairs of vectors are likewise defined. 
From these definitions, in the case of symmetric weights, there are parallels to standard differential calculus. There is an adjoint relation between the gradient and the divergence operators, as well as a divergence theorem, leading to a definition of the Laplacian being found to be self-adjoint.

\noindent For this work, the most important operators are the partial derivatives of $u(x_{i})$ with respect to $x^{\mu}$, the $\ith[\mu]$ component of the state $x$. Here, we denote differential calculus partial derivatives with $\partial$, and non-local partial derivatives with $\delta$. The argument of the function in $u(x_{i})$ also makes clear that the function's partial derivative is evaluated at vertex $i$. These non-local operations can be used to rigorously define the partial derivatives as the contraction between the function gradient and the specifically defined unit vector,
\begin{align}
	\difference[1]{u(x_{i})}{x^{\mu}} ~=&~ [\dotproduct[{\gradient[{w}][{[{u}](x_{i})}]}][{\hat{x}^{\mu}_{i}}]] \label{eq:diff{1}} \\
	 ~=&~~ \frac{1}{n-1}\sum_{j \in V\setminus \{i\}} [u(x_{j}) - u(x_{i})](x^{\mu}_{j} - x^{\mu}_{i})w(x_{i},x_{j}), 
\end{align}
where the unit vector at state $i$ in the $\mu$ direction is defined as a vector gradient between states $i$ and $j$:
\begin{equation}
	\hat{x}^{\mu}(x_{i},x_{j}) \equiv \gradient[{w}][[x^{\mu}]](x_{i},x_{j}) = [x^{\mu}_{j} - x^{\mu}_{i}]\sqrt{w(x_{i},x_{j})}.
\end{equation}
We will proceed with the weights $w(x_{i},x_{j})$ being chosen such that the unit vectors are normalized:
\begin{equation}
 	[\dotproduct[{\hat{x}^{\mu}}][{\hat{x}^{\nu}}]] = \delta^{\mu\nu}, \label{eq:unitnorm}
\end{equation}
where $\delta^{\mu\nu}$ is the Kronecker delta. As will be discussed, this constraint together with the inner products also being normalized by the volume of the space, leads to other important behaviors of the convergence and consistency of these non-local derivatives to their counterparts in differential calculus.

Higher order derivatives also can be computed as approximations to $\nderivative[q]{u(x)}{x^{\mu_{0}}}{x^{\mu_{q-1}}}$ using an extension of \cref{eq:diff{1}} and can be defined recursively:
\begin{align}
\ndifference[q]{u(x_{i})}{x^{\mu_{0}}}{x^{\mu_{q-1}}} ~=&~ 
[\dotproduct[{\gradient[{w}][{[ 
	{\dotproduct[{\gradient[{w}][{[
	{\cdots[\dotproduct[{\gradient[{w}][{[{\dotproduct[{\gradient[{w}][{[u(x_{i})]}]}][{\hat{x}^{\mu_{0}}}] }]}]}][{\hat{x}^{\mu_{1}}}]]\cdots}]}]}][{\hat{x}^{\mu_{q-2}}}]}]}]}][{\hat{x}^{\mu_{q-1}}}]] \label{eq:diff_p}
\end{align}

It will now be shown, that under specific constraints on the weight functions, the above partial derivatives are, to leading order, the corresponding partial derivatives of differential calculus. However certain properties of differential calculus, such as the commutativity of mixed partial derivatives, do not automatically hold to all orders of approximation with all possible choices of edge weights. Fortunately, with weights that decay suitably, issues of commutativity and convergence can be guaranteed.

\subsection{Constraints on edge weights} \label{sec:graphtheory_weightconstraints}
When using the given non-local calculus definitions, there is an arbitrary normalization of the quantities that can be chosen. For numerical purposes, and to ensure the operator definitions are consistent and independent of the size of the graph, as discussed, all integrations over the graph are normalized by the size of the discrete space we are in $\abs{V} = n$. In what follows repeated Greek indices $\mu \in \{ 0,\dots p-1\}$ will follow Einstein summation convention.

For the non-local calculus, we will observe order-by-order how the partial derivatives converge, given derived constraints on the weight function $w(x_{i},x_{j})$ for the $p$ dimensional representation of the state vector $x_{i}$. We begin by examining the normalization condition in \cref{eq:unitnorm} for this graph, with the normalization constant denoted as $W[p,n]$
\begin{align}
	\sum_{j \in V\setminus \{i\}} (x^{\mu}_{j}-x^{\mu}_{i})(x^{\nu}_{j}-x^{\nu}_{i})w(x_{i},x_{j}) = W^{\mu\nu}_{}[p,n] = (n-1)\delta^{\mu\nu} ~~~ \forall~ i. \label{eq:weightnorm}
\end{align}
Due to the quadratic form of the constraint, posed as a second order moment of this weight function, it is natural to assume that the weights are decaying, and take a form that decays at least as rapidly the inverse square of the distance between graph vertices:
\begin{align}
	w(x_{i},x_{j}) = O\left(\frac{1}{\abs{x_{j}-x_{i}}^2}\right).
\end{align}
For $p=1$, the form of the weight may be chosen as precisely the relation above, and for $x = x^{1}$,
\begin{align}
	w(x_{i},x_{j}) = \frac{1}{\abs{x_{j}-x_{i}}^2}.
\end{align}
For $p>1$ dimensional graphs, the form of the weight is less obvious if the Kronecker delta representation is to be extended or generalized. The weight function must decay rapidly enough to minimize consistency error and ensure that off-diagonal terms on the left-hand side of \cref{eq:weightnorm} are zero; however, too rapid a decay may fail to satisfy the finite diagonal constraints.

In the case of radially symmetric metrics $w({r_{ij}})$, where ${r_{ij}}^2 = \sum_{\mu} \abs{x^{\mu}_{j}-x^{\mu}_{i}}^2$, the normalization constraint can be written as 
\begin{align}
	\sum_{j \in V\setminus \{i\}} {r_{ij}}^2 w({r_{ij}}) = W[p,n] = (n-1)p ~~~ \forall~ i,
	\label{eq:discreteweightnormalization}
\end{align}
which follows by computing the trace in \cref{eq:weightnorm}. 

In generalizing \cref{eq:weightnorm} to the case of continuous manifolds, we consider the case of graphs whose states $x_{i}$ have support over a subset of $\reals[p]$, where the discrete sums approach integrals. Here we desire a form for the weights to guarantee that all off-diagonal terms, $\mu \neq \nu$ are identically zero. In this setting, \cref{eq:discreteweightnormalization} leads to the form
\begin{align}
	\int dr~ r^{p+1} w(r) = W[p,R] = p\frac{V_p}{S_p}R^p = R^p, \label{eq:weightrad}
\end{align}
where $V_p$ and $S_p$ are the volume and surface area of the unit ball in $p$ dimensions, $S_p = pV_p$, and $R$ is the radius of the sphere that encloses the graph.

In this continuous case of a graph whose states $x_{i}$ are dense in $\mathbb{R}^p$, the off-diagonal terms in \cref{eq:weightnorm} can be shown to cancel due to symmetry considerations, since $w(r)$ is an even symmetric function. Letting $x^{\mu} = rf^{\mu}(\theta)$, where $f$ are angular functions, the spherical symmetry ensures that only diagonal terms are non-zero in the tensor components
\begin{align}
	W^{\mu\nu}[p,R] &\equiv \int d^{p}y~ (y^{\mu} - x^{\mu})(y^{\nu} - x^{\nu})w(x,y).
\end{align}
For example in $p=2$ dimensions, where $x^{1} = r\cos{\theta}, x^{2} = r\sin{\theta}$, the angular integral trivially evaluates to zero. 

However, in the case of graphs that are not dense in $\reals[p]$, with radius
\begin{align}
	R = \max_{ij} \frac{{r_{ij}}}{2}
\end{align}
in order for these off-diagonal terms to cancel in these finite integrals for an arbitrary graph structure, in addition to being symmetric functions, the weight functions must also decay rapidly so that the odd, antisymmetric terms are negligible.

\noindent It will now be assumed that for $p>1$, the weights decay faster than the above assumed inverse square relationship:
\begin{align}
	w(r) = o\left(\frac{1}{r^2}\right).
\end{align}
There still remains some choice in the exact form of the weights, and for this work, we will consider radially symmetric weights that decay like a Gaussian function $\sim e^{-r^2}$, or decay polynomially $\sim 1/r^2$. Based on these forms of the weights, there are some free parameters that can be calculated or chosen. We choose the weights to have the form $w(r) = C ~\tilde{w}(r/\sigma)$, where $C(p,R)$ is the weight scale, and $\sigma(p,R)$ is the weight decay parameter, or an appropriate length scale in the graph manifold. Based on \cref{eq:discreteweightnormalization}, the weights must scale as the inverse of length-squared, and so we will discuss the weight properties in terms of the scaled weight constraints 
\begin{align}
\widetilde{W}[p,z] = \frac{W[p,R]}{C\sigma^{p+2}}, \label{eq:scaled_weight_constraint}
\end{align}
where $z = {R}/{\sigma}$ is a scaled radius. The weight constraints therefore take the form of $p$-order moments of the weight function:
\begin{align}
	\int_0^z dy~ {y}^{p+1} \tilde{w}(y) = \widetilde{W}[p,z]. \label{eq:constraint_scaled}
\end{align}

\subsubsection{Gaussian weights}\label{sec:graphtheory_weightconstraints_gaussian}
We first may consider Gaussian decaying weights
\begin{align}
	w(r) = Ce^{-\frac{1}{2}\frac{r^2}{\sigma^2}}, \label{eq:weight_gauss}
\end{align}
where in the case of graphs with a finite radius $R$, the weight is a truncated Gaussian for $0\leq r\leq R$. The scaled weight constraints are the $p+1$ order moments of the Gaussian distribution. These moments have a known recursion relationship \cite{Orjebin2014} for $p>0$:
\begin{align}
	\widetilde{W}[p,z] ~=&~ p\widetilde{W}[p-2,z] - z^pe^{-\frac{1}{2}z^2}, \label{eq:recurrence_gauss}
\end{align}
which also can be expressed in the analytic form
\begin{align}
	\widetilde{W}[p,z] ~=&~ p!!\left[\widetilde{W}[\nmodm{p}{2},z] - z^{p}e^{-\frac{1}{2}z^2}\sum_{i=0}^{\floor{p/2}-1} \frac{z^{-2i}}{(p-2i)!!}\right].
\end{align}

\subsubsection{Polynomial weights}\label{sec:graphtheory_weightconstraints_polynomial}
We also may consider polynomially decaying weights
\begin{align}
	w(r) = C\left(\frac{\sigma}{r}\right)^{2+\epsilon}, \label{eq:weight_poly}
\end{align}
where $0\leq\epsilon<p$ is the super-quadratic decay scaling. As previously discussed, the weights must decay faster than quadratically, however cannot decay so fast as for the integrals to vanish over discrete spherical volumes. Here the scaled weight constraints are tractable integrals for $p>0$:
\begin{align}
	\widetilde{W}[p,z] = \frac{1}{p-\epsilon} z^{p-\epsilon}.	\label{eq:constraint_poly}
\end{align}
Values for $\widetilde{W}[p,z]$ for radially symmetric Gaussian and polynomially decaying weights, with a finite radius, continuous graph are shown in \cref{tab:weights}. Here, $\Phi(x)$ is the standard normal cumulative distribution function.
\begin{table}[htb!]
	\centering
	\begin{tabular}{|c|c|c|} \hline\tableheight{20pt}
		$p$ & Gaussian & Polynomial \\ \hline\tableheight{20pt}
		$-2$ & $0$ & Undefined \\ \hline\tableheight{20pt}
		$-1$ & $\sqrt{2\pi}\left[\Phi(z) - \Phi(0)\right]$ & Undefined\\ \hline\tableheight{20pt}
		$0$ & $1 - e^{-\frac{1}{2}z^2}$ & Undefined \\ \hline\tableheight{20pt}
		$1$ & $\sqrt{2\pi}[\Phi(z) - \Phi(0)] - ze^{-\frac{1}{2}z^2}$ & $\frac{1}{1-\epsilon} z^{1-\epsilon}$ \\ \hline\tableheight{20pt}
		$2$ & $2 - (2+z^2)e^{-\frac{1}{2}z^2}$ & $\frac{1}{2-\epsilon} z^{2-\epsilon}$ \\ \hline
	\end{tabular}
	\caption{Weight constraint values for a radially symmetric weight $\widetilde{W}[p,z]$, given a finite radius $z = {R}/{\sigma}$, continuous graph in $p$ dimensions.}
	\label{tab:weights}
\end{table}

\subsubsection{Free parameters in weight functions}\label{sec:graphtheory_weightconstraints_freeparameters}

The free parameters $C$ and $\sigma$, and other hyperparameters such as the scaling $\epsilon$, may now be constrained based on for example \cref{eq:constraint_scaled} and \cref{eq:constraint_poly}, and there is some freedom, depending on the graph of interest. Generally, the decay parameter of the weight $\sigma$ can be chosen, after fixing the scale
\begin{align}
 C ~=&~ \frac{z^p}{\sigma^2}\frac{1}{\widetilde{W}[p,z]}, \label{eq:weightscale}
\end{align}
where the scaling of the weight by $1/\sigma^2$ is evident.

\begin{figure}[hpt]
\centering
\begin{subfigure}[t]{0.49\textwidth}
	\centering
	\includegraphics[width=0.6\textwidth]{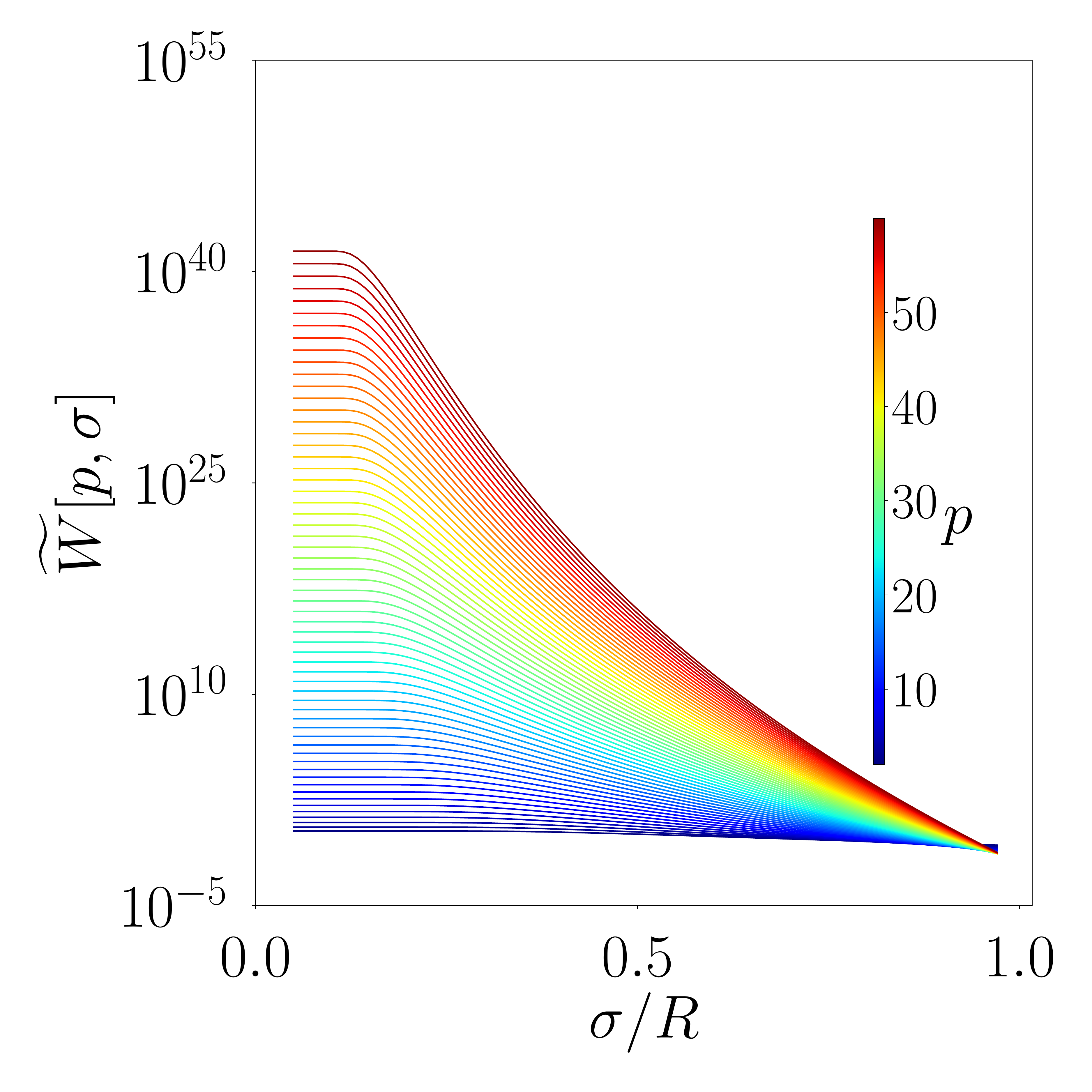}
	\subcaption{Gaussian}
	\label{fig:Wtilde_p_sigma_gauss}
\end{subfigure}
\hfill
\begin{subfigure}[t]{0.49\textwidth}
	\centering
	\includegraphics[width=0.6\textwidth]{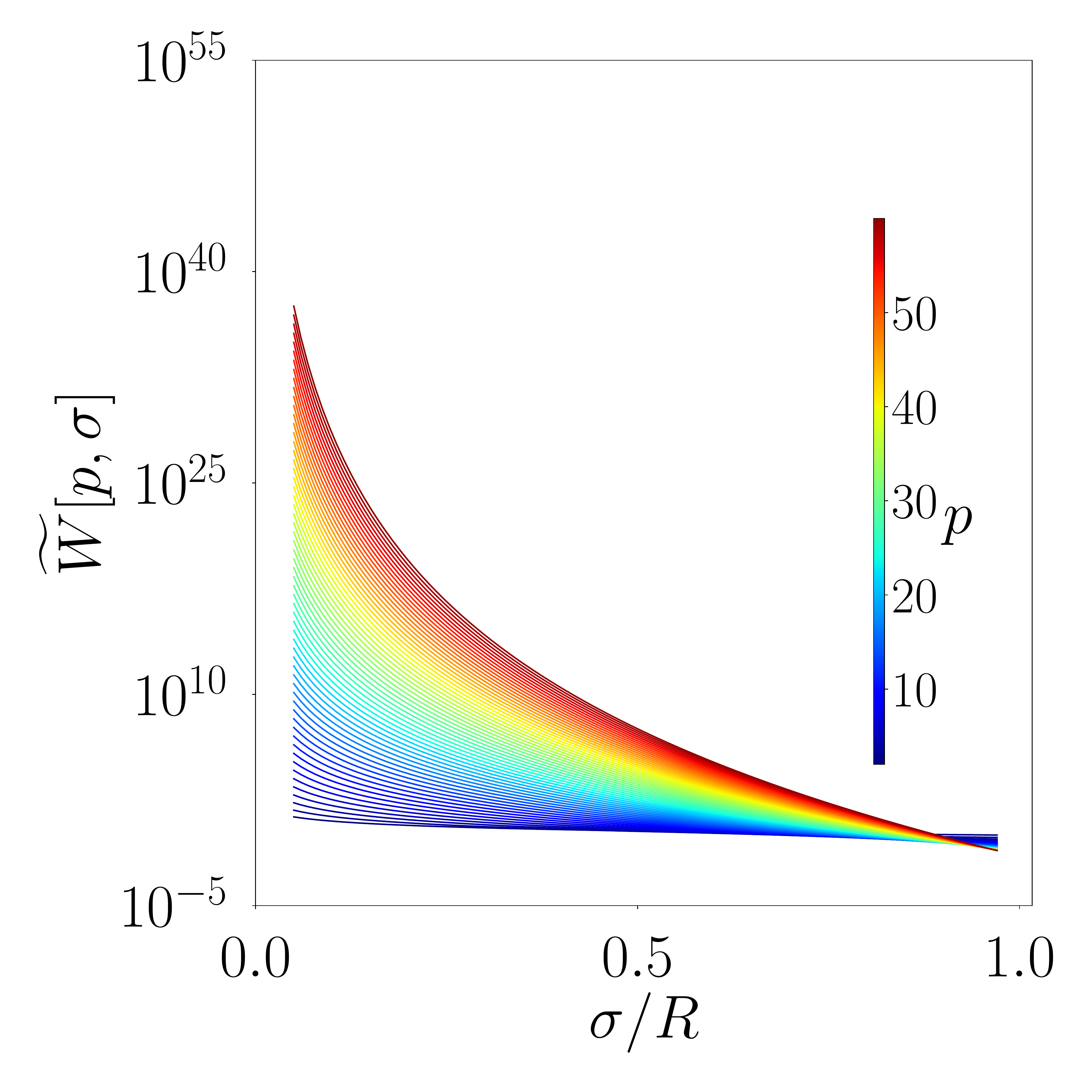}
	\subcaption{Polynomial}
	\label{fig:Wtilde_p_sigma_poly}
\end{subfigure}
\caption{Scaled weight constraints for various $p$ values, as a function of weight decay parameter $\sigma$; $\epsilon = {p}/{2}$ for the polynomial weights.}
\label{fig:Wtilde_p_sigma}
\end{figure}

\begin{figure}[hpt]
\centering
\begin{subfigure}[t]{0.49\textwidth}
	\centering
	\includegraphics[width=0.6\textwidth]{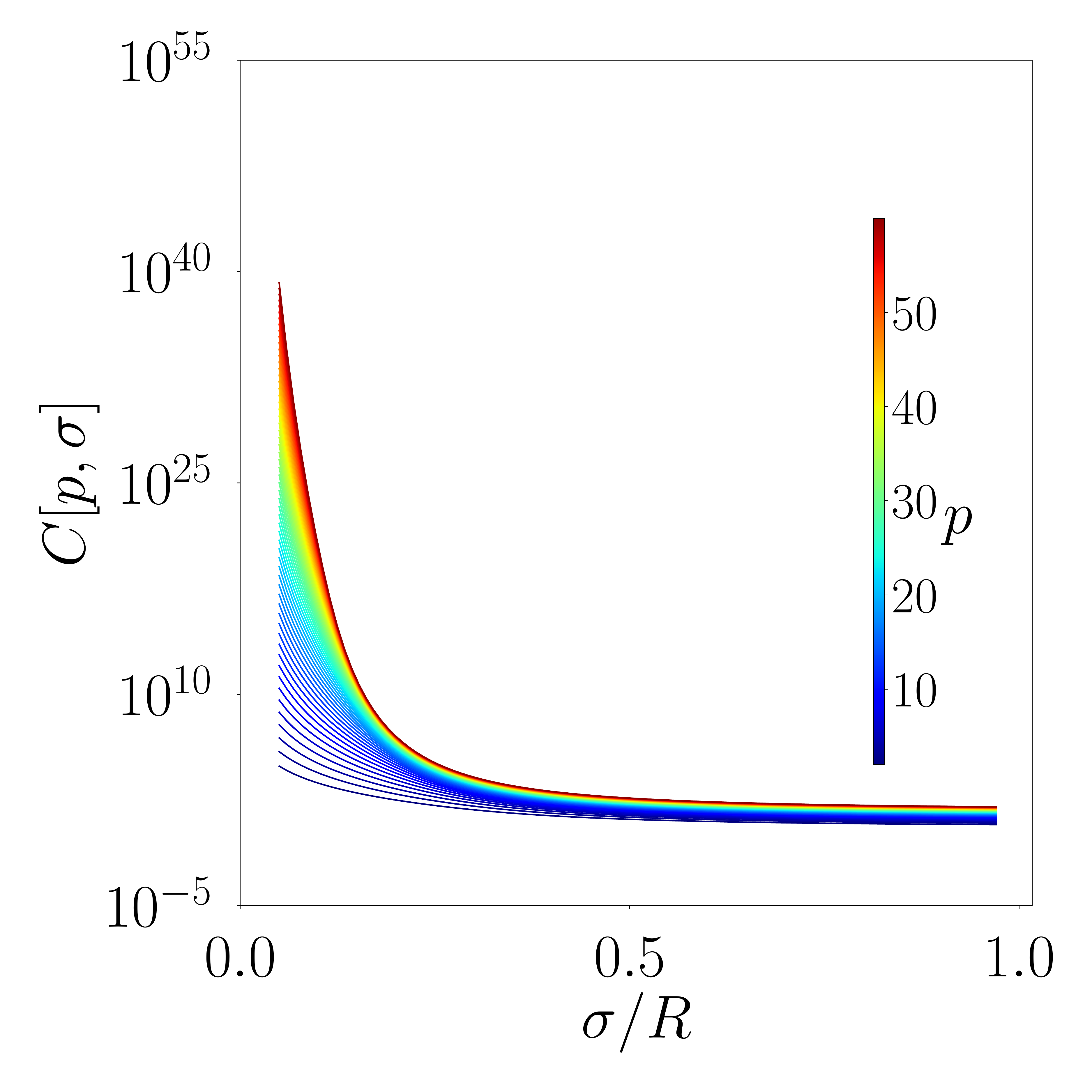}
	\subcaption{Gaussian}
	\label{fig:C_p_sigma_gauss}
\end{subfigure}
\hfill
\begin{subfigure}[t]{0.49\textwidth}
	\centering
	\includegraphics[width=0.6\textwidth]{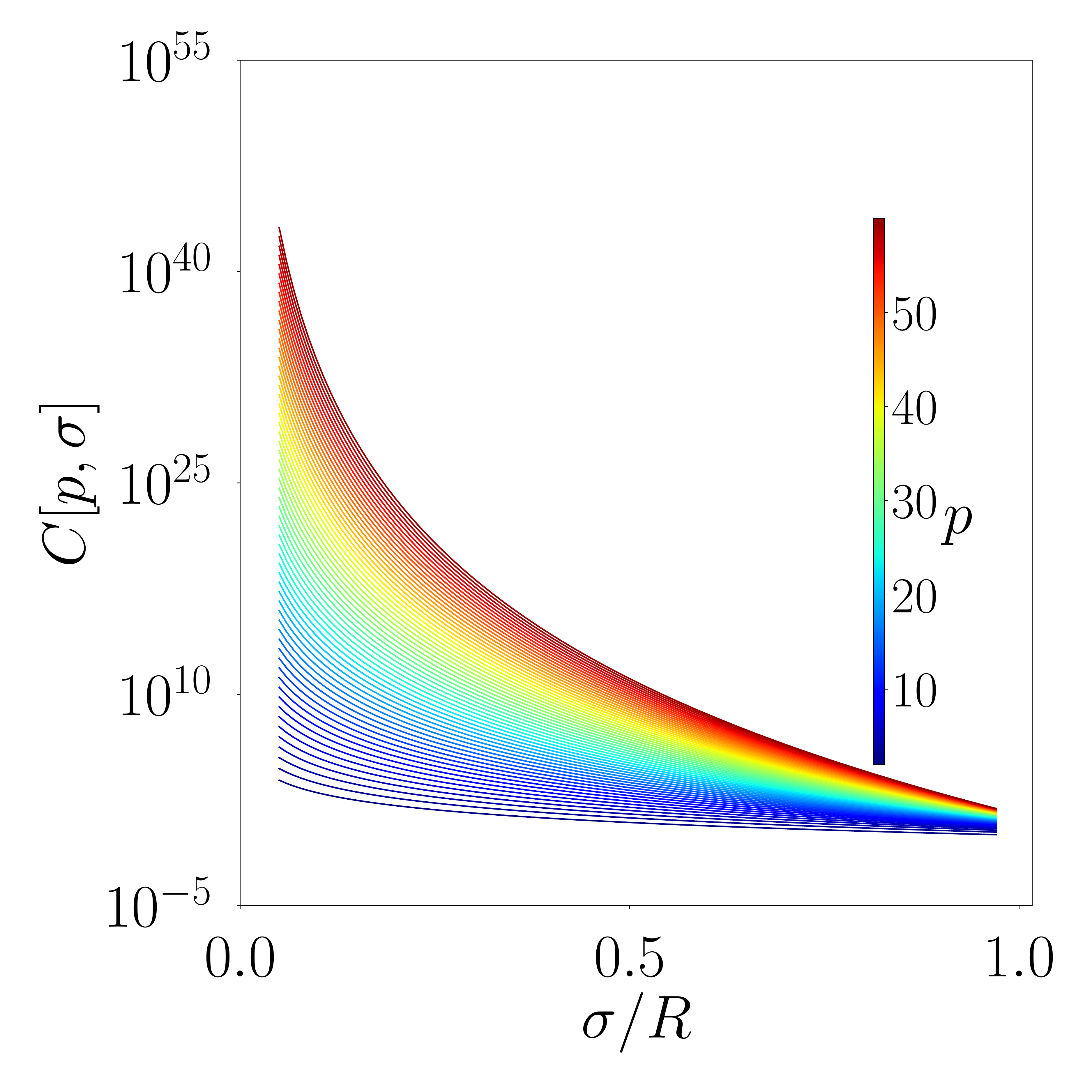}
	\subcaption{Polynomial}
	\label{fig:C_p_sigma_poly}
\end{subfigure}
\caption{Weight scales for various $p$ values, as a function of weight decay parameter $\sigma$; $\epsilon = {p}/{2}$ for the polynomial weights.}
\label{fig:C_p_sigma}
\end{figure}

\begin{figure}[hpt]
\centering
\begin{subfigure}[t]{0.49\textwidth}
	\centering
	\includegraphics[trim={0 2 0 0},clip,width=0.8\textwidth]{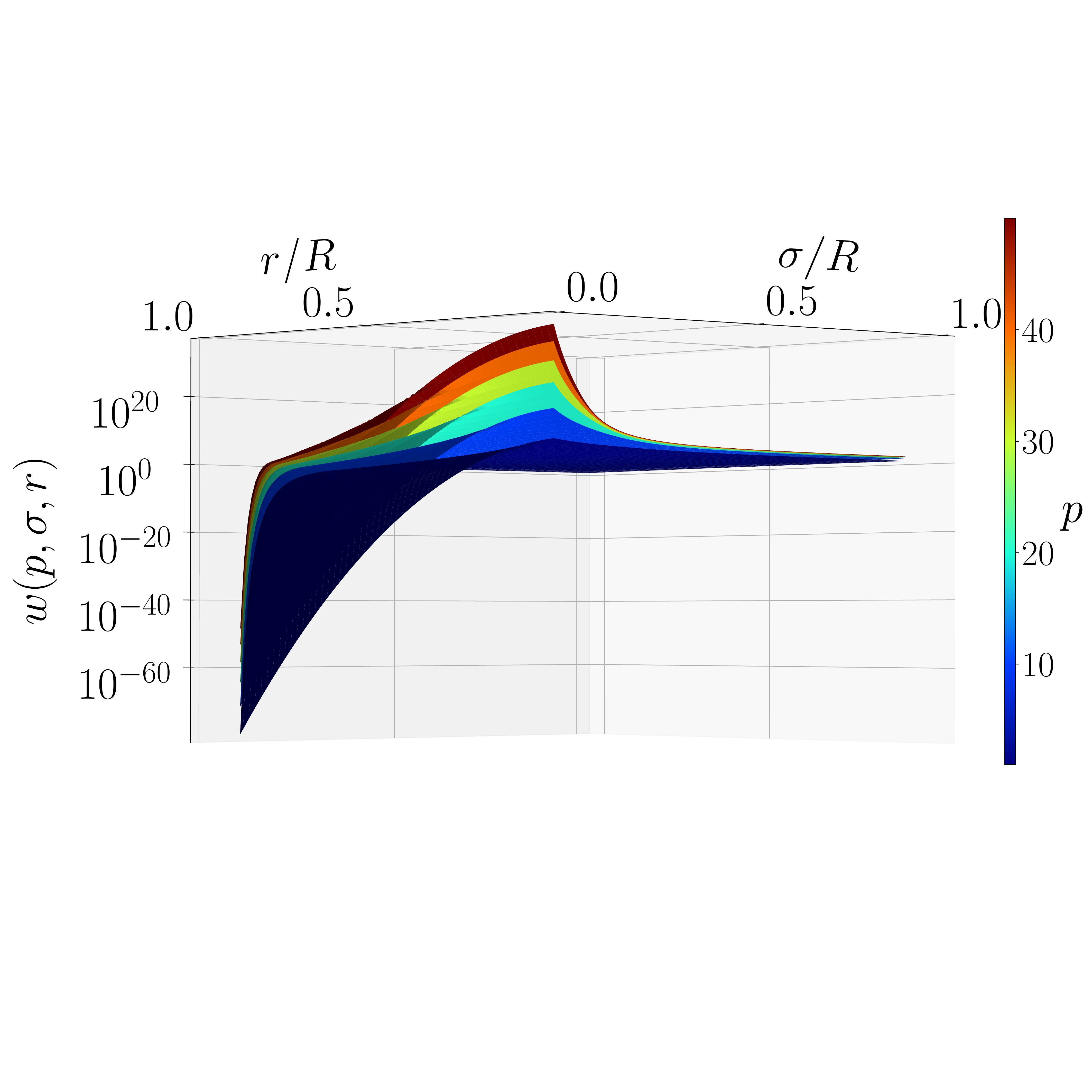}
	\vspace{-1.5cm}
	\subcaption{Gaussian}
	\label{fig:w_p_sigma_gauss}
\end{subfigure}
\hfill
\begin{subfigure}[t]{0.49\textwidth}
	\centering
	\includegraphics[trim={0 -10pt 0 0},clip,width=0.8\textwidth]{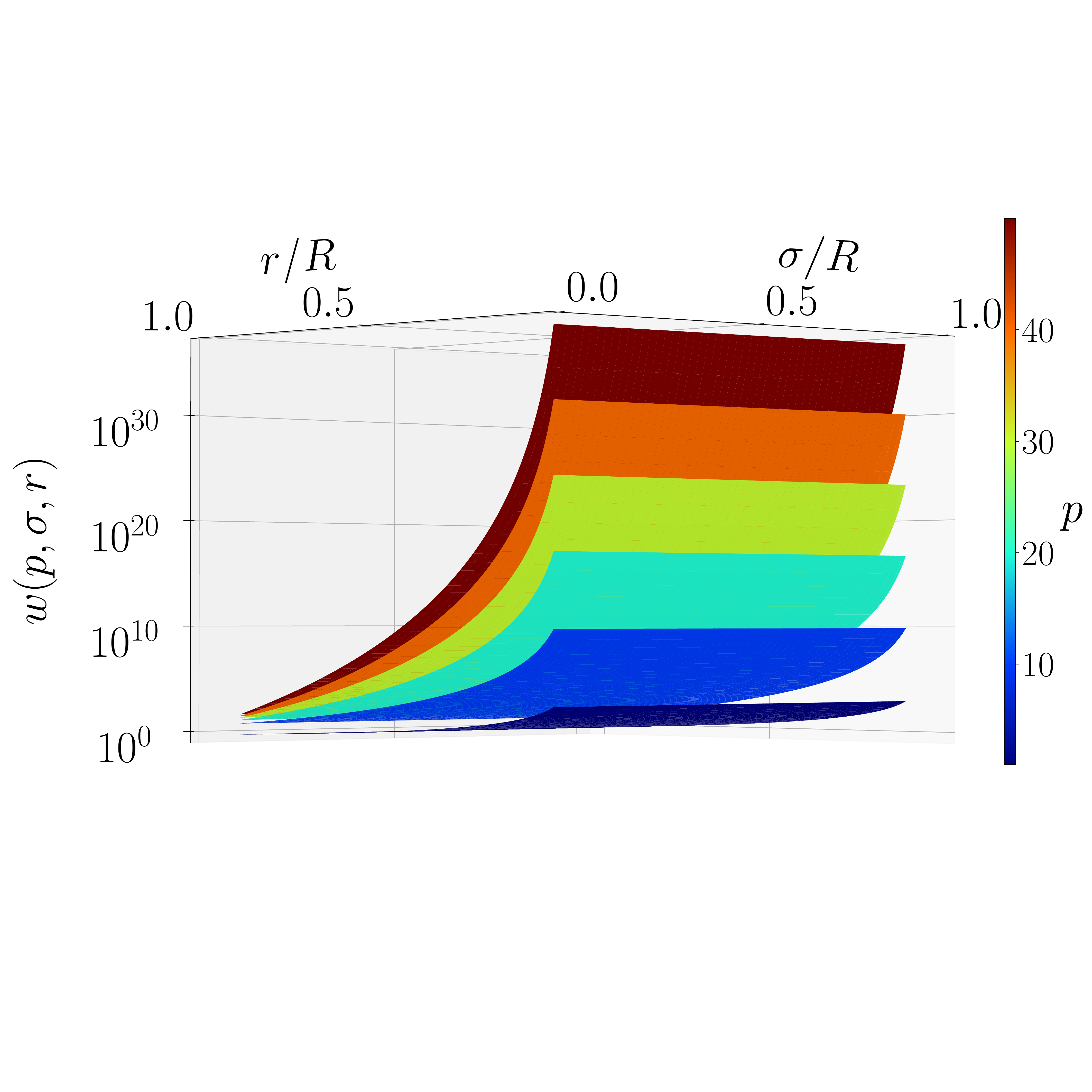}
	\vspace{-1.5cm}
	\subcaption{Polynomial}
	\label{fig:w_p_sigma_poly}
\end{subfigure}
\vspace{-0.5cm}
\caption{Weight functions for various $p$ values, as a function of weight decay parameter $\sigma$ and radius $r$; $\epsilon = {p}/{2}$ for the polynomial weights.}
\label{fig:w_p_sigma}
\end{figure}
Using Gaussian and polynomially decaying weights, the calculated weight constraints $\widetilde{W}[p,\sigma]$, and the weight scalings $C[p,\sigma]$ are shown in \cref{fig:Wtilde_p_sigma} and \cref{fig:C_p_sigma}. For the polynomial weights, $\epsilon = {p}/{2}$. Here it can be seen that the Gaussian weight constraints are constant for small $\sigma > 0$, and then decrease for $\sigma > {R}/{z_{\textrm{inflection}}}$, where $z_{\textrm{inflection}}^2 = p+{3}/{2} + \sqrt{2p+{9}/{2}}$ is the point at which the Gaussian case of the integrand in \cref{eq:constraint_scaled} changes curvature and tends to zero. The Gaussian weight scalings also decrease polynomially in $\sigma$, as per \cref{eq:weightscale}. The polynomial weight constraints decay as ${1}/{\sigma^{p-\epsilon}}$, and the polynomial weight scalings decay as ${1}/{\sigma^{2+\epsilon}}$. 

Importantly, for a fixed $\epsilon$, the polynomial weight scalings scale much slower with $\sigma$ and $p$ compared to the Gaussian weight scalings due to this polynomial form. For $\epsilon=0$ and a given $p$ value, how the polynomial scalings $C$ vary with $\sigma$ exactly cancel how the weight itself varies with $\sigma$. Therefore the polynomial weight function is constant with respect to $\sigma$:
\begin{align}
	w(r) = (p-\epsilon)\frac{R^{\epsilon}}{r^{2+\epsilon}}, \label{eq:poly_weight_full}
\end{align}
and suggests this weight form is more easily applied to graphs with different weight scales. In contrast, the Gaussian weights scale as ${1}/{\sigma^{2}} e^{{R^2}/{2\sigma^2}}$ with respect to $\sigma$ and increases exponentially for $\sigma \ll R$. Both of these calculations show the weight values also increasing with $p$, indicating that functions that depend on a large number of variables may still have very large weights, regardless of the $\sigma$ dependence. 

We note that the Gaussian weights bear connections to the diffusion kernel, which furnishes the probability distribution between states on a graph in the context of discovering the data embedding in a continuous manifold. Currently the form of the weights also has a single length scale $\sigma$, and the manifold is assumed to be homogeneous along each state variable dimension. Future communications will investigate the effects of using each independent equation along each dimension in \cref{eq:discreteweightnormalization} to determine constraints on length scales in each direction $\sigma^{\mu}$. For the case of discrete manifolds, certain approximations on these continuous weight constraints must be made. Please refer to \cref{app:implementation_gauss} of the Appendix for implementation details of these approximations.

\subsection{Consistency of non-local derivatives} \label{sec:graphtheory_consistency}
With this insight into the form of the weights, and their constraints, we now examine the consistency of the partial derivatives in this non-local calculus. The analysis will follow a similar approach to any numerical scheme, substituting Taylor series of differential calculus terms for various non-local terms, to determine the leading order behavior. Due to the weight constraint, we will see that depending on the decay rate of the weights, the leading order error is potentially first order in the spacing of the vertices.

\noindent For the first order derivatives the weight constraint yields
\begin{align}
	\difference[1]{u(x_{i})}{x^{\mu}} 
	~=&~ \derivative[1]{u(x_{i})}{x^{\mu}} + O\left(\derivative[2]{u(x_{i})}{x^{\nu},x^{\theta}}(x^{\nu}_{j}-x^{\nu}_{i})(x^{\theta}_{j}-x^{\theta}_{i})(x^{\mu}_{j}-x^{\mu}_{i}) w(x_{i},x_{j})\right). \label{eq:non-local1}
\end{align}
For second order derivatives, each non-local first derivative can be expanded in the form found above in \cref{eq:non-local1}, yielding
\begin{align}
	\difference[2]{u(x_{i})}{x^{\mu},x^{\nu}} 
	~=&~ \derivative[2]{u(x_{i})}{x^{\mu},x^{\nu}} + O\left(\derivative[2]{u(x_{i})}{x^{\theta},x^{\epsilon}} (x^{\theta}_{j}-x^{\theta}_{k})(x^{\epsilon}_{j}-x^{\epsilon}_{k})(x^{\mu}_{j}-x^{\mu}_{k})(x^{\nu}_{k}-x^{\nu}_{i}) w(x_{k},x_{j})w(x_{i},x_{k})\right).\label{eq:non-local2}
\end{align}
Therefore, given this weight constraint, the first order non-local derivatives to leading order correspond with the corresponding local derivatives. Although the second order non-local derivatives do not contain first order terms, the second order terms contain the other second order derivatives along the other dimensions. It can be concluded that in order for the non-local calculus definitions for the partial derivatives to equal the local calculus definitions to leading order, a constraint on the higher order moments of the weight function for each order of derivative is necessary. These constraints are assumed to be of the form of products of all differences in state variables at all combinations of dimension and pairs of vertices, weighted by the edge weight between those vertices:
\begin{align}
\ndifference[q]{u(x_{i})}{x^{\mu_{0}}}{x^{\mu_{q-1}}}  = \nderivative[q]{u(x_{i})}{x^{\mu_{0}}}{x^{\mu_{q-1}}} + O(x_{j}-x_{i}) \longleftrightarrow & \\ 
\sum_{j_{0}\dots j_{q-1} \in V\setminus \{i\}} (x^{\mu_{q-1}}_{j_{q-1}}-x^{\mu_{q-1}}_{j_{q-2}}) \cdots (x^{\mu_{0}}_{j_{0}}-x^{\mu_{0}}_{i}) w(x_{j_{q-2}},x_{j_{q-1}}) \cdots w(x_{i},x_{j_{0}}) ~=&~  \nonumber \\
(n-1)^q \delta^{\mu_{q-1}\mu_{q-2}}\cdots \delta^{\mu_{0}\mu_{i}}. \nonumber
\end{align}

As alluded to in the previous example of decaying Gaussian and polynomial weights, an alternative to choosing weights based on these difficult to construct higher order moment constraints, is to choose weights that decay faster than any finite polynomial, and still satisfy the primary constraint in \cref{eq:weightnorm}. Given a fast enough decaying weight, this will allow several issues to be resolved. First, all off-diagonal terms in the primary constraint will decay rapidly. Second, all higher order non-local derivatives will equal their differential derivative counterparts to leading order due to the weights decaying faster than the polynomials in the Taylor series. 

\subsubsection{Compromises with a non-local calculus approach}\label{sec:graphtheory_consistency_compromise}

However, due to this calculus involving sums over the entire graph, the non-local nature means the terms beyond the local calculus derivatives in the consistency relationships above in \cref{eq:non-local1,eq:non-local2} may be unreasonably large. This is particularly true if the weight scales are large as we have seen for the Gaussian weights. Therefore there may be no real similarity in the values and trends of the non-local derivatives compared to their differential counterparts. However, from studies of the total scale of the weights, due to the polynomially decaying weights having a closed form that is independent of the decay parameter, these weights are reasonable scales of order $O(1)$ for small $p$. 

Due to these numerical issues, we choose to use the polynomially decaying weights for the current studies, with $\epsilon = {p}/{2}$. Although not necessarily similar in form, the non-local derivative values will be similar in magnitude to what would be expected for differential derivatives, and decay reasonably quickly with $\epsilon>0$. Future communications will study in depth the effect of using different forms and parameter choices for the weights, and how these large non-local error terms can be understood and controlled. 

This extends the compromise regarding the weight functions to be made, between obeying the consistency and constraint equations, and having scales of derivatives that correspond with those from differential calculus in the model. As to be discussed when these methods are applied to physical systems, the model fitting procedure may be able to overcome the large scales of these derivatives, such as associating small or large linear coefficients to terms in the model that contain these large values. Please refer to \cref{app:implementation_linearregression} of the Appendix for details on the regression implementation.

This non-local calculus has the objective of being a general, accurate framework for any modelling procedure, and the exact form of the non-local error term will be developed. An important aspect of numerical method development is an analysis of the effect of data sparsity on model error. As can be seen in \cref{eq:non-local1}, the non-local derivatives have a persistent error term for all weight functions with non-local extent, and whether this is detrimental to subsequent reduced order modelling is an important question in determining the validity of this approach. 

For example, a modified Taylor series, as discussed in \cref{app:implementation_modifiedtaylorseries} of the Appendix, may be constructed using a training dataset as a functional representation, where a different Taylor series model is developed at each possible base point. Please refer to \cref{app:error_mesh} for details concerning the mesh of data used. The point-wise error for evaluation points within a distance $h$ of each base point model can then be computed, and refining the spacing of the data, we are able to get a sense of how the method error scales with $h$. A complete error analysis for this modified Taylor series approach can be referred to in \cref{app:error} for $p=1$, and in \cref{app:error_p} for higher dimensions. As discussed in \cref{app:error_pointwisederivative,app:error_pointwisgamma,app:error_totalerror} of the Appendix, in addition to the expected effects of a $\ith[k]$ order Taylor series having error that scales as $h^{k+1}$, the non-local calculus derivatives, and any fit model parameters are also dependent on the spatial distribution of the training and testing data.

With this understanding of the behaviors of the non-local calculus and partial derivatives, reduced order models can now be developed for systems of interest. The modelling will consist of finding a basis of operators for a given model, where the operators are defined in terms of this non-local calculus. In the following section, an example of a physical systems is presented, and several reduced order modelling approaches are examined.


\section{Physical systems of interest} \label{sec:physicalsystems}
\subsection{Microstructures in a gradient-regularized model of non-convex elasticity} \label{sec:physicalsystems_microstructures}
We apply the graph theoretic approach to studying reduced order models for the mechanochemical response of solids that develop distinct regions of local phases over their evolution. This physical system is driven by a free energy density function that is parameterized by composition (chemical) and strain (elastic) variables. The underlying functional form of the free energy is non-convex in both these quantities and is regularized by composition and strain gradients. Microstructures develop as phases and symmetry-breaking structural variants arise corresponding to negative eigenvalues of the Hessian of the free energy density in the strain-composition space. These models are ill-posed, however; a condition arising from the absence of penalization on inter-phase and inter-variant interfaces and manifesting in PDE solutions that are pathologically discretization-dependent. Well-posedness and penalization of interfaces are restored by free energy functionals $\psi_{\textrm{tot}} = \psi_{\textrm{hom}} + \psi_{\textrm{grad}}$, dependent on the composition and strain gradients. When decoupled, the chemical component of this problem causes spinodal decomposition and Ostwald ripening, described by the Cahn-Hilliard equation \cite{Cahn1958}. The non-linear elasticity problem is described by a non-convex strain energy density stabilized by incorporating strain gradient effects as described by many authors. 

Our treatment follows Toupin \cite{Toupin1962} and has appeared as decoupled, gradient-regularization of non-convex, non-linear elasticity \cite{Rudraraju2014,Teichertetal2017}, as well as mechanochemical spinodal decomposition \cite{Rudraraju2016,Sagiyamaetal2016,SagiyamaGarikipati2017b}. The free energy density functional is $\psi= \psi_{\textrm{tot}}$, with a dependence on the fields $\{c,\vec{E},\gradient[]{c},\gradient[]{\vec{E}}\}$ where $c$ is the composition and $\vec{E}$ is the Green-Lagrange strain tensor. This dependence, however, is expressed in terms of symmetry-coding strain parameters that are introduced below.

At low temperatures, these systems consist of stable phases related to the local crystalline symmetry, specifically cubic, or tetragonal symmetry. The system considered is restricted to $d=2$ spatial dimensions, and so the material will consist of phases with either square or rectangular structure. This pattern of regions of phases of a field are generally referred to as microstructures. This system is further modeled as a binary mixture, with scalar composition $0\leq c \leq 1$ of one of the two constituents.

Additionally, the introduction of mechanical strain couplings to the composition allows for a degeneracy in the rectangular phase. Linear combinations of the strain components $\vec{E}$ yield relevant parameters $\vec{e}$ relating to the original high symmetry of the square phase, namely $e_{1,2} \equiv ({E_{11} \pm E_{22}})/{\sqrt{2}}, e_{6} \equiv \sqrt{2} E_{12}$, that allow this degeneracy to be determined. Specifically, the sign of the $e_{2}$ field can determine the orientation of the rectangular phase, either in the horizontal "plus", or the vertical "minus" variant. The order parameters will thus be defined to be $\eta = (c,e_{2})$, and the microstructure can be described by its phase volume fractions $\varphi_{\alpha}$, where $\alpha = \{\square,\hrectangle,\vrectangle\}$ corresponds to the phase of the material.

\subsubsection{Free energy description} \label{sec:physicalsystems_microstructures_freeenergy}
\begin{figure}[ht]
	\centering
	\includegraphics[width=0.6\textwidth]{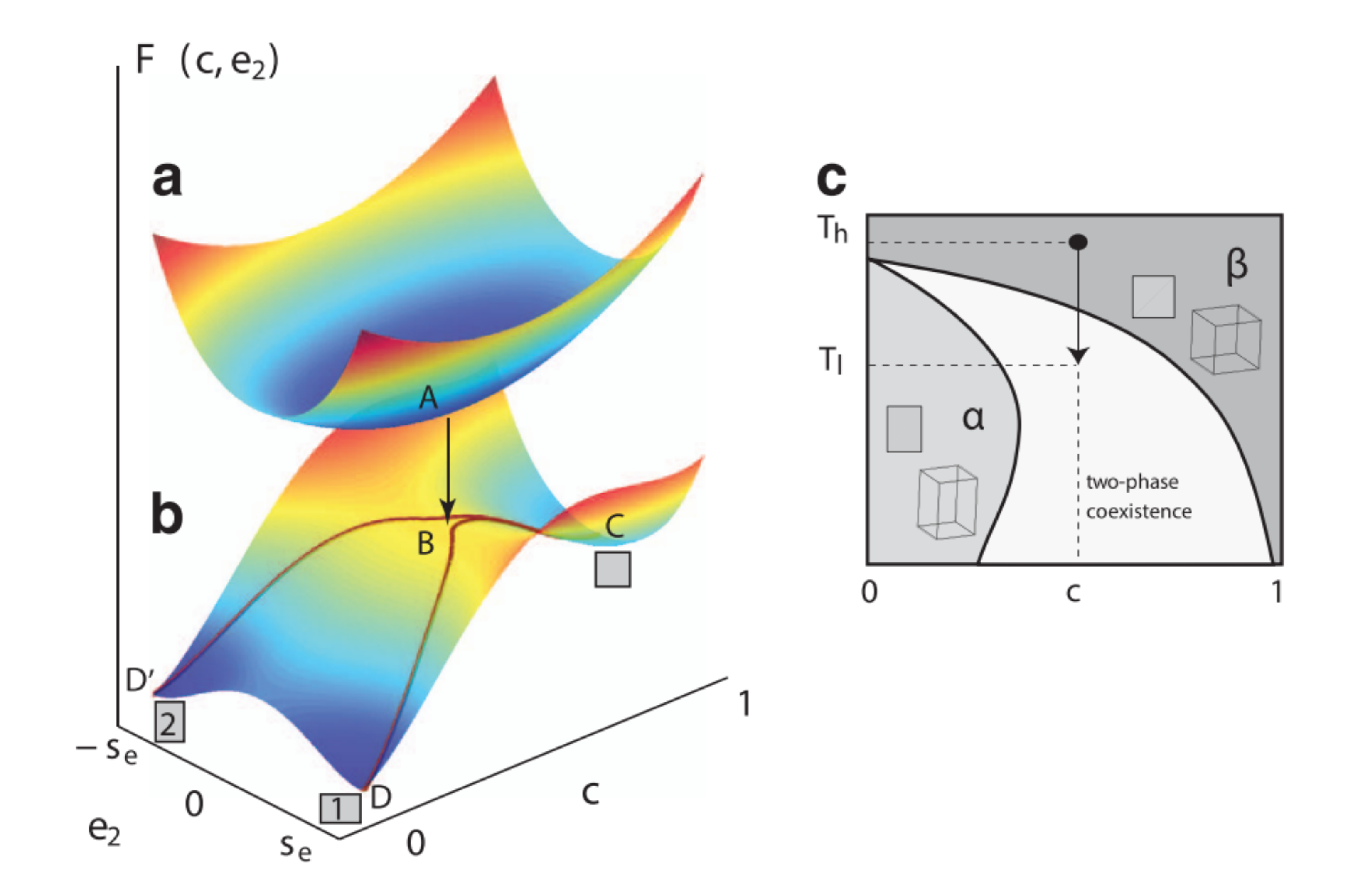}
	\caption{Free energy density landscape over the composition $c$ and strain component $e_{2}$ order parameters, showing changes in curvature and resulting phases of square and rectangular lattice symmetries, from Rudraraju \etal \cite{Rudraraju2016}. On the right, the effect of temperature and composition on the phase diagram of stable regions of purely square or rectangular phases, or a mixture of phases is shown.}
	\label{fig:energyphase}
\end{figure}

To describe these phase transitions, we take the description used by Zhang \& Garikipati\cite{Zhang2020}. The total free energy density description of the system $\psi$ can be comprised of a homogeneous component, $\psi_{\textrm{hom}}(c,\vec{e})$, that is not entirely concave up over the entire $c-e_{2}$ domain, yielding potential wells with stable phases. An example two dimensional energy landscape from Rudraraju \etal \cite{Rudraraju2016} is shown in \cref{fig:energyphase}.

During the mechanical-chemical coupled dynamical processes, within the domain where the free energy is concave, such as around the saddle point $B$ in \cref{fig:energyphase}, phase separation will tend to occur, where a mixture of the stable phases, will yield a lower energetic state. This can be seen from a typical lever law analysis for a binary mixture, to obtain the phase volume fractions $\varphi_{\alpha}$. This simple phase separation within the non-convex region of the free energy density is mechano-chemical spinodal decomposition, where phases tend to separate. Furthermore, the gradient regularization that represents the energy of interfaces also manifests in a number $N_{\alpha}$ of each variant, and total interfacial lengths $l_{\alpha}$.

With the composition $c$ and symmetry-coding strain parameter $\vec{e}$ and their gradients, a non-uniform free energy density component $\psi_{\textrm{grad}}(c,\vec{e},\gradient[]{c}, \gradient[]{\vec{e}})$ allows for gradient regularization, typically through (coupled) quadratic terms in the order parameter gradients. These terms yield phase regions that are limited by interfaces being less energetically favorable, and so equilibrium states will tend to have less numbers of domains of phases, and furthermore, larger phase regions grow at the expense of smaller ones. This effect is known as Oswald ripening \cite{Rudraraju2016}.

This gradient contribution, in addition to the homogeneous contribution yields the total free energy functional of the order parameters over the material volume $\Omega_t$ at time $t$,
\begin{equation}
	\Psi[c,\vec{e}] = \int_{\Omega_t} dV~ (\psi_{\textrm{hom}}+\psi_{\textrm{grad}}) + \int_{\partial \Omega_t} dS~ \dotproduct[{\vec{T}}][{\vec{u}}], \label{eq:Psi_tot}
\end{equation}
where $\vec{T}$ are mechanical traction terms, depending on which boundary conditions are specified, and $\vec{u}$ is the spatial displacement field of the material.

A variational formalism $\difference[1]{\Psi}{{c,\vec{e}}}$ can then be applied to yield the relevant generalized chemical potentials. As discussed, for this particular work, the homogeneous free energy density is a smooth polynomial function of the order parameters,
\begin{equation} 
	\psi_{\textrm{hom}}(c,\vec{e}) \equiv 16 \alpha_{c} c^4 - 32\alpha_{c} c^3 + \alpha_{c} c^2 + 2\frac{\alpha_{e}}{{\beta_{e}}^2}({e_{1}}^2 + {e_{6}}^2) + \frac{\alpha_{e}}{{\beta_{e}}^2}{e_{2}}^4 + 2\frac{\alpha_{e}}{{\beta_{e}}^2}(2c-1){e_{2}}^2,
	\label{eq:psi_hom}
\end{equation}
where $\{\alpha_{c},\alpha_{e},\beta_{e}\}$ are constant coefficients.

The gradient contributions will take the form of strictly quadratic forms in the order parameter gradients,
\begin{equation}
	\psi_{\textrm{grad}}(\gradient[]{c},\gradient[]{\vec{e}}) = \frac{1}{2} \gradient[]{c} \cdot\vec{\kappa} \gradient[]{c} + \frac{1}{2} \gradient[]{\vec{e}} \cdot \boldsymbol{\gamma} \gradient[]{\vec{e}} + \gradient[]{c}\cdot \boldsymbol{\theta} \gradient[]{\vec{e}}, \label{eq:psi_grad}
\end{equation}
and for this work, it is assumed that there are no gradient couplings between the composition and strain fields $\boldsymbol{\theta} = 0$, and only constant isotropic tensors $\vec{\kappa} = \kappa \textrm{I}$ and $\boldsymbol{\gamma} = \gamma \vec{I}$ for the gradient terms. Furthermore, the additive terms mean the total free energy density can be separated into the sum of independent mechanical and chemical terms, plus a coupling term:
\begin{equation}
	\psi= \psi_{\textrm{mech}} + \psi_{\textrm{chem}} + \psi_{\textrm{couple}}.
\end{equation}

When stationarity of the free energy is imposed in a variational framework, the Euler-Lagrange equations yield the weak form of the balance of linear momentum, and further variational arguments lead to the strong form. Using this variational approach, based on this form of the free energy, and given its dependence on the strain and composition fields, the governing PDE's are chosen to be coupled Cahn-Hilliard dynamics for the composition \cite{Cahn1958}, and non-linear gradient elasticity equations in equilibrium \cite{Toupin1962,Rudraraju2016}.

\subsubsection{Cahn-Hilliard dynamics} \label{sec:physicalsystems_microstructures_CahnHilliard}
Cahn-Hilliard dynamics are first order for the composition and take the form of a transport equation
\begin{align}
	\derivative[1]{c}{t} + \divergence[][{\vec{J}}] ~=&~ 0, \label{eq:pde_cahnhilliard} 
\end{align}
where
\begin{align}
	\vec{J} ~=&~ -\vec{L}\gradient[][\mu], \label{eq:cflux}
\end{align}
and $\vec{L} = \vec{L}(c,\vec{e})$ is a mobility transport tensor, and is assumed to be isotropic $\vec{L} = \textrm{L}\vec{I}$. The variational treatment yields the chemical potential $\mu = \difference[1]{\Psi}{c}$, whose gradient yields the flux in a form that is guided by the thermodynamic dissipation inequality \cite{degrootmazur1984}, and leads to there being fourth-order mechano-chemically coupled equations. 

\subsubsection{Toupin model of elasticity} \label{sec:physicalsystems_microstructures_Toupin}
The strain gradient elasticity can be modeled in equilibrium as equations for the first Piola-Kirchhoff $\vec{P}$ and higher order $\vec{B}$ stress tensors,
\begin{align}
	\vec{P} ~=&~~ \derivative[1]{\psi}{\vec{F}} \label{eq:piola1}\\
	\vec{B} ~=&~~ \derivative[1]{\psi}{\gradient[]{\vec{F}}} \label{eq:piola3}
\end{align}
which are derivatives of the strain energy and are conjugate to the deformation gradient $\vec{F}$. Here the Green-Lagrange strains are $\vec{E} = {1}/{2}\left(\vec{F}^{T}\vec{F} - \vec{I}\right)$ and gradients $\gradient[]$ are derivatives with respect to the material indices. The elasticity governing equations, shown here in the strong form, are
\begin{align}
	\divergence[]{\vec{P}} - \divergence[]{(\divergence[]{\vec{B}})} = 0, \label{eq:pde_elasticity}
\end{align}
plus additional Dirichlet and higher order traction boundary conditions.

\subsubsection{Direct numerical simulations} \label{sec:physicalsystems_microstructures_DNS}
\begin{figure}[hpt]
\centering
\includegraphics[width=0.6\textwidth]{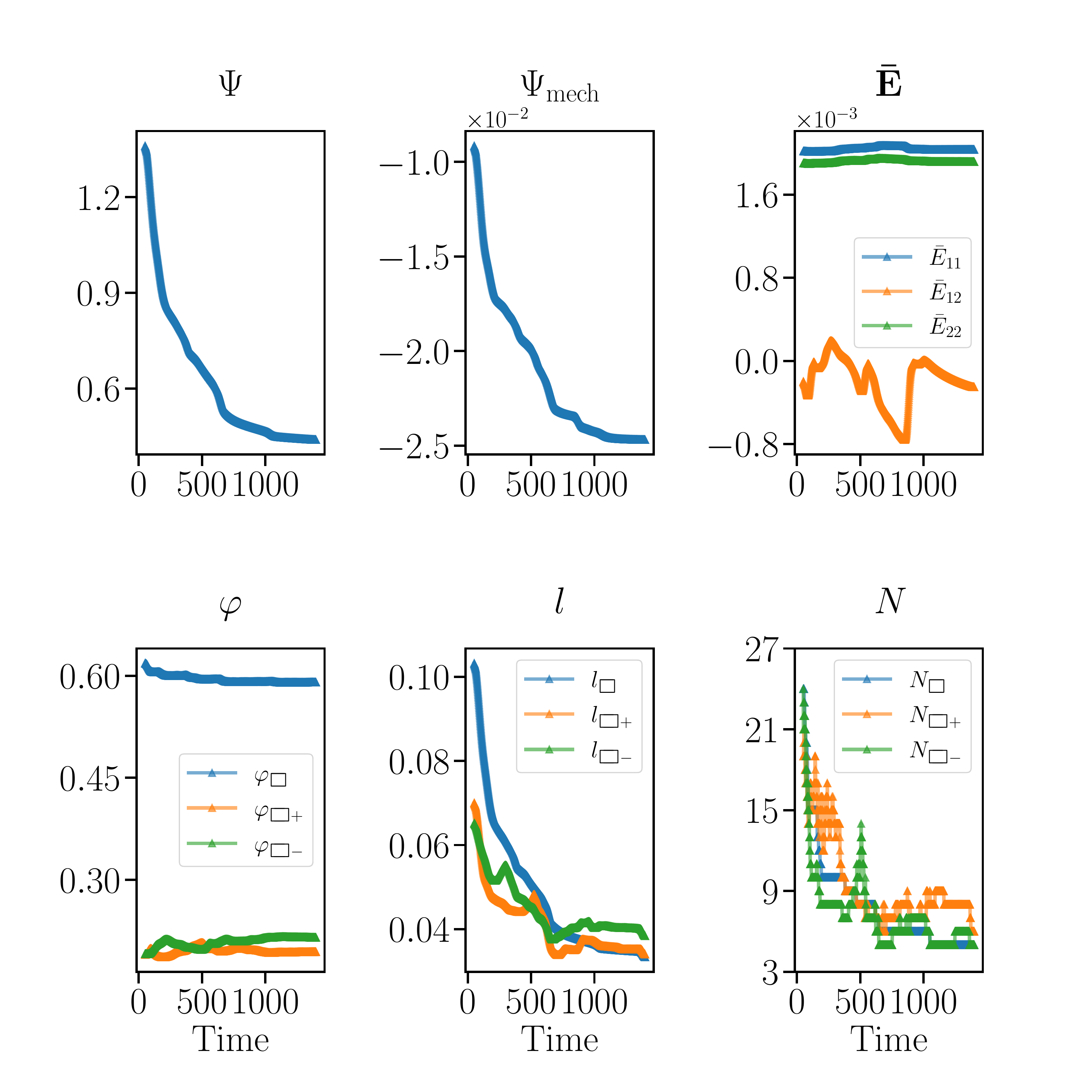}
\caption{DNS calculated state vector components.}
\label{fig:DNSobservables}
\end{figure}

These governing PDEs in \cref{eq:pde_cahnhilliard,eq:pde_elasticity} may now be solved, yielding time series data, as per Zhang \& Garikipati \cite{Zhang2020}. The weak form is chosen here for implementation, motivated by numerical methods such as finite elements and the isogeometric approaches. \cite{Rudraraju2014,Rudraraju2016} A given system is described by the static mechanical Dirichlet boundary conditions for the displacement fields $\vec{u}$, uniform initial $e_{2} = 0$ fields, corresponding to the material being initially in the high temperature square phase, random initial compositions $c = 0.46 \pm 0.05$, and zero chemical flux boundaries. The DNS proceeds by time stepping with steps of $\Delta t = 10^{-5}~\textrm{s}$, and iteratively solving the PDE's at each time step, with tolerances placed on the residual of the solver. The data generated for this work is at times $50-1500 ~\mu\textrm{s}$, with example DNS observables shown in \cref{fig:DNSobservables}. In this work, we consider a single boundary and initial condition and one evolution of the microstructure and aim to model its quantities of interest. In future communications we will consider the effectiveness of the graph theoretic approach at producing general models for a family of microstructures with different boundary and initial conditions.

To develop the graph theoretic model for the given mechano-chemical system, the state vector is assembled from quantities of interest that are functionals of the high-dimensional solution to \cref{eq:pde_cahnhilliard,eq:pde_elasticity}. These quantities, denoted as barred variables unless otherwise specified, are effective volume averaged quantities $\bar{f}$ over the material domain. They therefore contain no spatial dependence of their high dimensional counterparts $f(\eta(x))$. Chemical quantities such as the local composition and chemical potential are represented by the phase volume fractions, interfacial lengths, and numbers of phase domains
\begin{align}
\{c(x),~\mu(x)\} ~\to&~ \{\varphi_{\alpha},~l_{\alpha},~N_{\alpha} \}_{\alpha =\{\square,\hrectangle,\vrectangle,\dots~\}}, \nonumber \\
\intertext{mechanical quantities such as the deformation gradients, Green-Langrange strains, and Piola-Kirchoff stresses are represented by their volume averaged quantities}
\{\vec{F}(x),~\vec{E}(x),~\vec{P}(x),\dots\} ~\to&~ \{\vecbar{F},~\vecbar{E},~\vecbar{P},\dots\}, \nonumber \\
\intertext{and free energy quantities are represented by their total values}
\{\psi(x),~\psi_{\textrm{mech}}(x),~\psi_{\textrm{chem}}(x),\dots\} ~\to&~ \{\Psi,~\Psi_{\textrm{mech}},~\Psi_{\textrm{chem}},\dots \}.
\end{align}
The phase volume fractions $\varphi$ are calculated by summing over volume elements that are within the associated phase, as per the phase diagram in \cref{fig:energyphase}; the interfacial lengths $l_{\alpha}$ are calculated by computing the total lengths of the interfaces between the square-rectangular-(plus,minus) fields; the number of domains $N_{\alpha}$ are counted through a gradient edge detection algorithm to identify the number of distinct domains of each field. The volume like terms of $l_{\alpha}$ and $N_{\alpha}$ are essential for distinguishing changes in topology of the microstructure, such as when field domains merge, separate, appear, or disappear. The specific forms of the free energy $\Psi, \Psi_{\textrm{mech}}$ are the integrated free energy densities over the material volume. For the strain quantities, they are chosen to be defined in terms of the volume averaged Green-Lagrange strains
\begin{align}
	\vecbar{E} = \frac{1}{2}\left(\vecbar{F}^{T}\vecbar{F} - \vec{I}\right),
\end{align} 
where it is the deformation gradient $\vec{F}$ that has been explicitly volume averaged
\begin{align}
\vecbar{F} = \int_{\Omega_t} dV~ \vec{F}.
\end{align}

Due to the high dimensional data being generated by a finite element DNS method, each of the spatially dependent quantities can be thought of as fluctuating about their volume averaged value, and could be represented by being expanded as a Taylor series
\begin{align}
f(\eta) = \bar{f} + \epsilon + \dotproduct[{\derivative[1]{f(\tilde{\eta})}{\eta}}][{(\eta-\tilde{\eta})}] + O((\eta-\tilde{\eta})^2) \label{eq:local_avg_quantities}
\end{align}
in the order parameters, where $\tilde{\eta} = \argmin_{\eta} \abs{f(\eta)-\bar{f}}$, such that $f(\tilde{\eta}) = \bar{f} + \epsilon$. This polynomial series means any dependencies of the spatially dependent quantities on the order parameters will become higher order polynomial dependencies in the volumed averaged quantities; $\bar{f}$ will be a higher order polynomial in $\bar{\eta}$. For example, the free energy density contains up to fourth order polynomials of the order parameters, however the total free energy will be a function of much higher polynomials of the effective order parameter quantities. These effects will be discussed when determining forms of the reduced order models while being informed of the original dependencies that yield the training data.

The relevant $p$-dimensional state vector $x = \{x^{\mu}\}$ associated with each vertex $i$ of the graph is chosen to be
\begin{equation}
	x = \{\Psi,\Psi_{\textrm{mech}} ,\vecbar{E},\varphi_{\alpha},l_{\alpha},N_{\alpha}\}, \label{eq:observables}
\end{equation}
and vertices in this graph are associated with the time index. From physical arguments, the microstructure states at different times are all related, and therefore the graph is considered to be fully connected when considering edges and edge weights.

\subsection{Reduced order modelling} \label{sec:physicalsystems_reducedordermodelling}
To determine a model for the system, two different aspects of the model are considered, the time dependence of the phase volume fractions, when zero additional mechanical stress is applied to the system, and the free energy dependence on the strain and phase volume fractions. The forms of these reduced order models will follow the general principals of how the various fields are coupled from this variational framework. This choice of PDE's involves polynomial couplings between the composition and strain fields in the free energy density in \cref{eq:psi_hom}, however the constitutive relations for the chemical potential do not involve explicit derivatives of the free energy density with respect to the strains or strain gradients, and the constitutive relations for the elasticity do not involve derivatives of the free energy density with respect to the compositions or composition gradients. These polynomial and derivative couplings indicate important quantities of interest are those dominated by either algebraic or differential forms of the chemical and mechanical effects.

\subsubsection{Selection of physically relevant basis} \label{sec:physicalsystems_reducedordermodelling_Basis}
Given the microstructure data is generated from explicit forms for the free energy density in \cref{eq:psi_hom,eq:psi_grad}, these exact forms give an indication for which algebraic and differential terms should be proposed for the model. The free energy contains no explicit spatial dependence, and contains monomials up to $c^4$ and ${{e_{\alpha}}}^4$ order, plus coupled polynomial terms that are up to $c {e_{2}}^2$ order, plus quadratic gradient terms. The chemical potential and associated current are of the form
\begin{align}
	\mu ~=&~ \derivative[1]{\psi}{c} - \kappa \laplacian[][{c}] \\
	\divergence[][{\vec{J}}] ~=&~ -L\left[\uniderivative[3]{\psi}{c}\dotproduct[{\gradient[]{c}}][{\gradient[]{c}}] + 2\derivative[3]{\psi}{c,c,{e_{\alpha}}}\dotproduct[{\gradient[]{c}}][{\gradient[]{{e_{\alpha}}}}] + \derivative[3]{\psi}{c,{e_{\alpha}},{e_{\beta}}}\dotproduct[{\gradient[]{{e_{\alpha}}}}][{\gradient[]{{e_{\beta}}}}] \right.\\
	&~~~~~~~~~~~ + \left.\uniderivative[2]{\psi}{c}\laplacian[][{c}] +\derivative[2]{\psi}{c,{e_{\alpha}}}\laplacian[][{e_{\alpha}}] - \kappa \laplacian[][]\laplacian[][{c}]\right] \nonumber.
\end{align}
\noindent The equations are solved in the weak form over the material domain, and therefore any spatial dependencies, and particularly any spatial derivatives of the order parameters, are removed, yielding terms that are proportional to the parameter being differentiated, similar to \cref{eq:local_avg_quantities}. For example $\int_{\Omega_t} dV~ \nlaplacian[][{\eta}][n] \sim O(\eta)$. Therefore the Cahn-Hilliard first order dynamics for the composition should be of the order
\begin{align}
	\derivative[1]{c}{t} &\sim O(\uniderivative[3]{\psi}{c} c^2) + O(\derivative[3]{\psi}{c,c,{e_{\alpha}}} c {e_{\alpha}})
	+ O(\derivative[3]{\psi}{c,{{e_{\alpha}}},{{e_{\beta}}}}~ {{e_{\alpha}}}{e_{\beta}}) \label{eq:cahnhilliard_approx}	\\
	&~ + O(\uniderivative[2]{\psi}{c} c) + O(\derivative[2]{\psi}{c,{e_{\alpha}}} {e_{\alpha}}) + O(c) \nonumber.
\end{align}
For the elasticity equations, given the free energy and its derivatives, upon being integrated over the material volume are polynomials in the order parameters, it will be assumed that a functional representation of the free energy can be expressed as a finite polynomial series in the phase volume fractions and volume averaged strain parameters.

When performing this reduced order modelling of the volume averaged quantities, we will be informed by the various equations and relations used to produce this data in determining which operators should be present in the models. Both the specific terms present, the form of the higher dimensional representations, such as polynomial expansions, differential relationships, or general dependencies, as well as the number of terms required to form the model will inform the hypothesis of the model. Specific implementation details are discussed in \cref{app:implementation_stepwiseregression} of the Appendix, as well as in Wang \etal \cite{Wang2019,Wang2020} and Kochunas \etal \cite{Kochunas2020}. 

Here, all partial derivatives are approximated using the non-local calculus definitions in \cref{eq:diff_p}, with the exception of the time derivatives, which are calculated using the standard backwards Euler finite difference method. It is also important to note when computing these operators and forming a basis, finite datasets mean numerically some terms may be not linearly independent. This will affect the span of this basis, and if there are more linearly dependent terms in the model, some models with different subsets of operators may be degenerate. This is an effect that will be discussed for the specific models below.

\subsubsection{Weighted residual loss functions} \label{sec:physicalsystems_reducedordermodelling_lossfuncs}
An important aspect of the fitting procedure is the choice of loss function to characterize the error between the known data and predictions. For general models, we choose the root mean square error with the $L_{2}$ norm. The $L_{2}$ norm loss generally minimizes fluctuations in the data and leads to smooth curves, whereas the $L_{1}$ norm loss generally minimizes error at specific points, and the $L_{\infty}$ norm bounds the maximum error of a single point. For highly oscillatory data, such as mechano-chemical systems which can change their local phase rapidly during their transient response, a weighted loss function with several norms is more appropriate. Here, we define the total loss function as a weighted sum of loss functions
\begin{align}
	\ell = \sum_{\chi \in \{1,2,\infty,\cdots\}} w_{\chi}\cdot \ell_{\chi}, \label{eq:lossfunc}
\end{align}
where $\chi$ denotes the type of norm $L_{\chi}$ used in the loss function, with weighting $w_{\chi}$. Each loss function $\ell_{\chi}$ may have its own sub-domain to evaluate the loss. In general, we choose the $\ell_{2}$ domain to be over the entire dataset. For oscillatory problems we choose the $\ell_{1}$ loss to be evaluated at the local maxima and minima of the data, through a specified number of passes of a peak finding algorithm.

\subsubsection{Regression approaches} \label{sec:physicalsystems_reducedordermodelling_regression}
The general modelling approach will consist of proposing a basis of algebraic and differential operators. System inference methods, specifically linear stepwise regression approaches, allow the model to be made as parsimonious as possible, while remaining an accurate representation of the underlying physics. 

Here, we choose to perform regression using Ordinary Least Squares (OLS), as well as Ridge regression (Ridge)\cite{Lv2009} and Lasso regression (Lasso) \cite{Tibshirani1996}, schemes which constrain the norm of the coefficients. Regularization is added to the regression with the objective of improving the robustness of the method, particularly for low rank basis of operators that induce very high condition numbers in the linear problems to be solved. In the case of Ridge regression, it can be directly shown that regularization increases the rank and decreases the collinearity of the data. Please refer to \cref{app:implementation_linearregression_ridge} and \cref{app:implementation_linearregression_lasso} of the Appendix for details on these regression approaches and the residual equations being minimized. 

The Ridge and Lasso regression approaches generally use a cross validation scheme to determine the optimal ratio $\lambda$ of the primary constraint and the additional norm constraints in the total residual loss function. At each iteration of the stepwise regression, we perform an independent cross validation to determine the optimal $\lambda$ hyper-parameter for each given reduced basis of operations with a proposed operator removed. This approach increases the computational complexity at least linearly with the size of the data, and linearly with the number of possible hyperparameters, however potentially ensures more optimal coefficients that are away from a local basin of the loss function. We observe that the Lasso regression induces sparsity in the regression coefficients due to using a $1$-norm constraint on the coefficients. This conflicts with the stepwise regression approach, and appears less suitable for these types of problems. We therefore will be comparing the OLS and Ridge approaches.

This effect of the independent cross validations yielding optimal coefficients for a given proposed model at a stepwise iteration, but not accounting for the overall stepwise procedure that is selecting which coefficients to remove from the model, will be seen to affect the overall trends of the model loss. There is potentially a competing interest between the penalization of the loss function with the coefficient norm constraints, and the desire to remove insignificant terms and produce the most effective model. 

\textbf{Remark}: It is also important to acknowledge that the cross validation may be performed using a different optimality metric compared to the loss function for the statistical tests in the stepwise regression. The global minima of the two metrics do not necessarily coincide, and as will be discussed in \cref{sec:physicalsystems_vol}, this greatly affects the model, and the order of operators present after the overall stepwise regression is complete. Ultimately, we desire a robust, accurate, and efficient method, and the following example studies indicate the compromises to be made between these objectives.

\subsection{Phase volume fraction dynamics}\label{sec:physicalsystems_vol}
To model the phase volume fraction dynamics, we are guided by the Cahn-Hilliard equation in \cref{eq:pde_cahnhilliard}, and the order of terms present in \cref{eq:cahnhilliard_approx}. Here the composition evolves under first-order dynamics driven by a conservative flux that itself has a non-linear, and high-order gradient dependence on the composition and strains. Importantly, the flux is dependent on the chemical potential of $c$ and stresses $\vecbar{E}$ that arise as variational derivatives of the free energy with respect to the composition and strains, respectively in \cref{eq:cflux}. We therefore seek evolution equations for the volume phase volume fractions $\varphi$ that are also first-order in time, and dependent on polynomials of products of the free energy derivatives up to $\difference[3]{\Psi}{\varphi_{\alpha},\varphi_{\beta},\bar{E}_{\xi \kappa}}$, with phase dependent terms in $\{\varphi_{\alpha},l_{\alpha}, N_{\alpha},\vecbar{E}\}$. We further restrict the possible dependencies based on the form of the free energy in $d$ spatial dimensions, a function of a single chemical parameter and ${d(d+1)}/{2}$ mechanical parameters, which we choose to be $\bar{\varphi} = \varphi_{\square}$, and $\vecbar{E}$. 

We also assume that any dependencies on the microstructure parameter terms $l_{\alpha},N_{\alpha}$ show up only in the terms that contain polynomials of the state variables, the last term in \cref{eq:cahnhilliard_approx}, and not in terms that contain derivatives. These polynomial terms without derivatives arise in the Cahn-Hilliard equation from the interfacial gradient terms in the free energy. The interfacial lengths and particle numbers, here defined to be $\bar{l}_{1} = l_{\square},~\bar{N}_{1} = N_{\square},~\bar{l}_{2} = {(l_{\vrectangle}+
l_{\vrectangle}-l_{\square})}/{2},~\bar{N}_{2} = N_{\hrectangle}+N_{\vrectangle}$ will suffice to model the $\kappa$ and $\gamma$ interfacial terms.

\subsubsection{Proposed model for phase volume fraction dynamics}\label{sec:physicalsystems_vol_proposedmodel}
The model for the phase volume fraction first order dynamics is assumed to be of the following form
\begin{equation}
	\derivative[1]{\bar{\varphi}}{t} = f(\nderivative[q]{\Psi}{\bar{\varphi}}{\bar{E}_{\alpha \beta}},\bar{\varphi},\bar{E}_{\alpha \beta},\bar{l}_{\alpha}, \bar{N}_{\alpha}), \label{eq:model_vol}
\end{equation}
where $\nderivative[q]{\Psi}{\bar{\varphi}}{\bar{E}_{\alpha \beta}}  = \{\derivative[2]{\Psi}{\bar{\varphi},\bar{E}_{\alpha \beta}},\derivative[2]{\Psi}{\bar{\varphi},\bar{\varphi}},\derivative[3]{\Psi}{\bar{\varphi},\bar{E}_{\alpha \beta},\bar{E}_{\xi \kappa}},\derivative[3]{\Psi}{\bar{\varphi},\bar{\varphi},\bar{E}_{\alpha \beta}},\uniderivative[3]{\Psi}{{\bar{\varphi}}}\}$ are the relevant derivatives for this model.

The function $f$ for this work is an algebraic, linear function of polynomials of the state vector components and derivatives of the free energy
\begin{equation}
	f = \mspace{-20mu}\sum_{{\substack{a = \left\{a_{\nderivative[q]{\Psi}{\bar{\varphi}}{\bar{E}_{\alpha \beta}}}\cdots~ a_{\bar{N}_{\beta}}\right\}\\\abs{a}\leq A\\ 0 \leq a_{\nderivative[q]{\Psi}{\bar{\varphi}}{\bar{E}_{\alpha \beta}}} \neq a_{\bar{l}_{\alpha}} \neq a_{\bar{N}_{\alpha}} \leq 1}}} \mspace{-40mu} \gamma^{a}~ ({\nderivative[q]{\Psi}{\bar{\varphi}}{\bar{E}_{\alpha \beta}})^{a_{\nderivative[q]{\Psi}{\bar{\varphi}}{\bar{E}_{\alpha \beta}}}}({\bar{\varphi}})^{a_{\bar{\varphi}}} ({\bar{E}_{\alpha \beta}})^{a_{\bar{E}_{\alpha \beta}}}({\bar{E}_{\xi \kappa}})^{a_{\bar{E}_{\xi \kappa}}}({\bar{l}_{\alpha}})^{a_{\bar{l}_{\alpha}}} ({\bar{N}_{\beta}})^{a_{\bar{N}_{\beta}}}} \label{eq:model_volfunc}
\end{equation}
where the powers of the quantities are represented by components of the multicomponent index $a$, with the indices being denoted by the quantities themselves. Here, summation over repeated indices in a given term should not be assumed. Instead, the function is explicitly a sum over all combinations of products of powers of the quantities, including all combinations of components of the tensor quantities. The powers are constrained such that the derivative $a_{\nderivativeinline[q]{\Psi}{\bar{\varphi}}{\bar{E}_{\alpha \beta}}}$, interfacial length $a_{\bar{l}_{\alpha}}$, and number $a_{\bar{N}_{\alpha}}$ terms are all separate from each other and never appear more than once in a term, as per \cref{eq:cahnhilliard_approx}. The polynomial terms are up to order $A=3$ powers, and there are a maximum of $308$ terms in this model. 

\subsubsection{Numerical implementation of phase volume fraction dynamics}\label{sec:physicalsystems_vol_numerical}
Given this basis for the terms that drive the phase volume fraction dynamics, stepwise regression is performed using a mixture of the $\ell_{2}$ and $\ell_{1}$ loss functions, as per \cref{eq:lossfunc}. Various loss weights are chosen and a single pass of a peak finding algorithm is used to find the points at which the $\ell_{1}$ loss is computed. The resulting loss curves with a comparison of OLS and Ridge regression fitting are shown in \cref{fig:vol_loss_OLS_Ridge17}. From initial leave-one-out cross-validation studies, it is found that any ridge parameter $\lambda > 0$ improves the rank of the matrix, and for this problem, the fit losses increase with increased $\lambda$. Therefore, $\lambda = 10^{-17}$ is fixed to perform Ridge regression. The optimal direct solver for Ridge regression for this problem is found to be an SVD based solver, whereas OLS regression uses a standard least squares solver. Please refer to a description of the graph theory library in \cref{app:implementation_graphtheorylibrary} of the Appendix for details on the specific direct linear solvers.

Initial studies are also performed using an $R^2$ correlation coefficient for the cross validation metric, while using the weighted loss function for the stepwise regression statistical criteria. This choice of different metrics results in significantly higher errors, of the order $O(10^{-1})$ and the stepwise regression procedure yields a significantly more oscillatory loss curve, as opposed to the desired monotonically increasing curve. The optimal coefficients given the correlation coefficient or the weighted loss metric do not coincide and the cross validation and stepwise regression are competing to the detriment of the overall method. Therefore all cross validation studies will be performed using an identical loss function for all metrics.

Each choice of loss weighting yields similar trends in the loss curves increasing with increased model parsimony, offset by how much the $\ell_{1}$ contributes to the total loss. The high frequency oscillations, as per \cref{fig:DNSobservables} in the phase volume fractions suggests why the fits with solely $\ell_{1}$ have the lowest loss. When looking at the exact ordering of operators that the stepwise regression presents as most important to the model, the different loss weights do not result in identical orderings, even at low numbers of operators where the loss curves are more similar. However, the general trends of which operators are present at various stages of the stepwise regression, denoted by the various plateaus in the loss curves, are very similar across choices of loss weights. We therefore choose to focus on analyzing the fits where solely the $\ell_{1}$ loss is used.

\subsubsection{Reduced order model for phase volume fraction dynamics}\label{sec:physicalsystems_vol_finalmodel}
Looking at the plateaus and sharp increases in the weighted loss curves, particularly for the OLS model in \cref{fig:vol_loss_OLS} over the stepwise iterations, the first plateau extends from $308$ down to approximately past $187$ operators, which corresponds to the rank of the matrix of this basis of operators. In general, as the model becomes more parsimonious, the loss should increase monotonically, however in these cases of lower rank basis, there may be some decreases when certain operators are rejected from the model. The second plateau ends at approximately $47$ operators, which is also the point in the stepwise iterations where the model first has a term that contains a derivative. Practically all iterations between $47$ and $308$ have derivative terms, and this suggests that the purely algebraic polynomial terms in $\{\varphi_{\alpha},l_{\alpha}, N_{\alpha},\vecbar{E}\}$ are over-fitting the model. The derivatives are not seen to play as important a role in the first order dynamics. However, the suggested importance of the algebraic polynomial terms may also be explained by higher order polynomials of the volume averaged quantities being required to represent dynamics controlled by spatially varying fields, as per \cref{eq:local_avg_quantities}. 

With loss weights $w_2 = 0,~ w_1 = 1$, the first 10 terms from the stepwise regression for the first order dynamics of $\bar{\varphi}$ are found using OLS and Ridge regression to be
\input{./figures/model_vol_OLS_Ridge17.tex}
It can be seen that both models consist of solely algebraic terms. The OLS and Ridge approaches also generate models with similar terms, with the Ridge regression model having 9 terms compared to 8 terms for OLS regression that contain predominantly strain dependencies. 

The first order dynamics, and fits for models with $20$, $50$, and $200$ terms is shown in \cref{fig:vol_fit_OLS_Ridge17}. For the practically full model with $200$ terms, the fits correspond to the data very well at early times, but have difficulty matching all of the small oscillations at late times. As the model becomes more parsimonious to $50$ terms where the sharp increase in the loss occurs, the general trend of the dynamics is captured, however this model misses capturing the full amplitude of the major peaks compared to the more complex models. Finally, with $20$ terms, which corresponds to the number of terms present in the original Cahn-Hilliard equation in $d=2$ dimensions, the large peaks, and small oscillations are not entirely captured, however the general trends are still shown. As discussed, the algebraic polynomial terms are potentially over-fitting, however manage to retain the general behavior of the phase volume fractions changing over time.

\subsubsection{Regression approaches for phase volume fraction dynamics}\label{sec:physicalsystems_vol_regression}
\begin{figure}[hpt]
\centering
\begin{subfigure}[t]{0.49\textwidth}
	\centering
	\includegraphics[width=\textwidth]{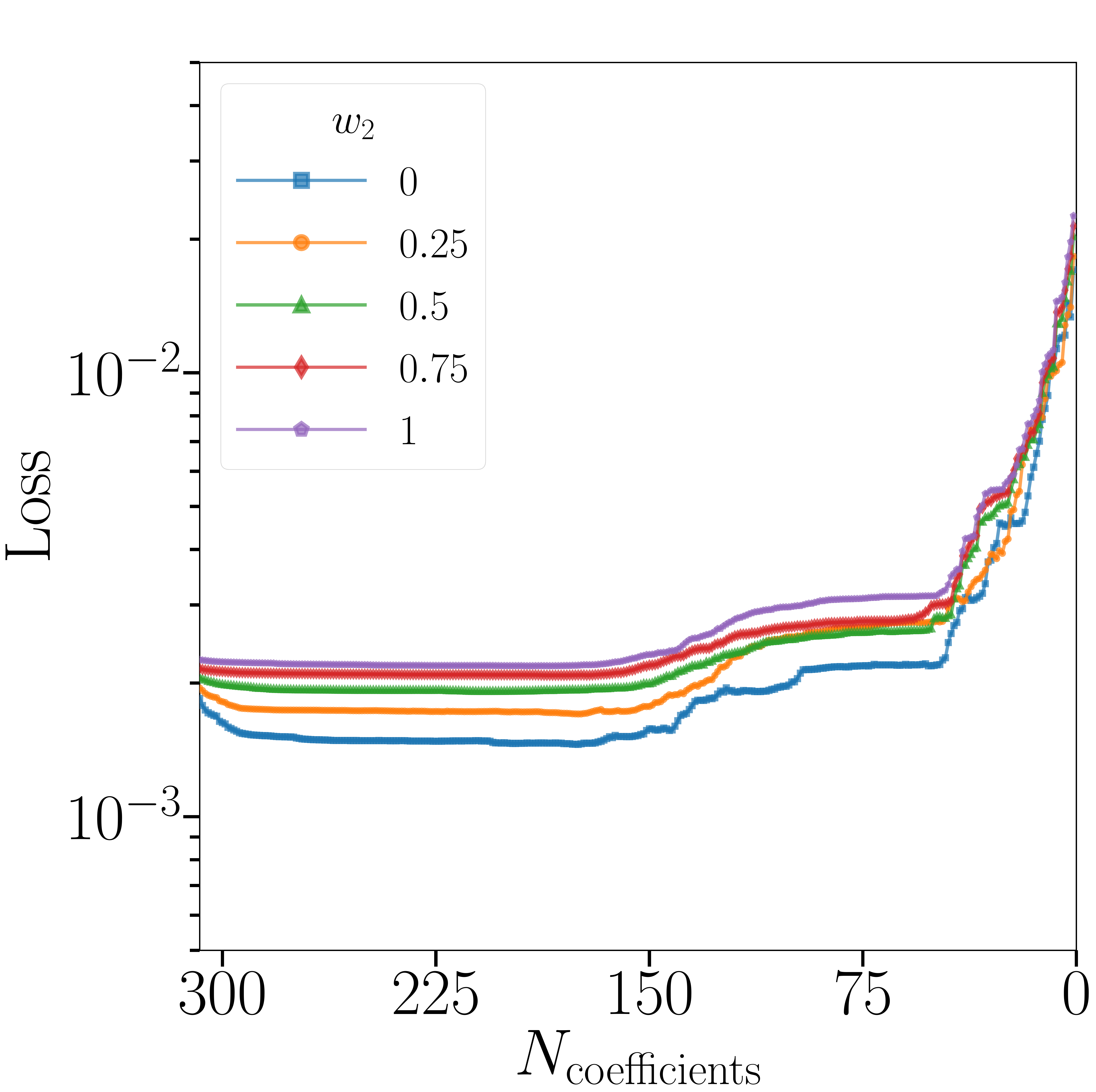}
	\subcaption{OLS regression.}
	\label{fig:vol_loss_OLS}
\end{subfigure}
\hfill
\begin{subfigure}[t]{0.49\textwidth}
	\centering
	\includegraphics[width=\textwidth]{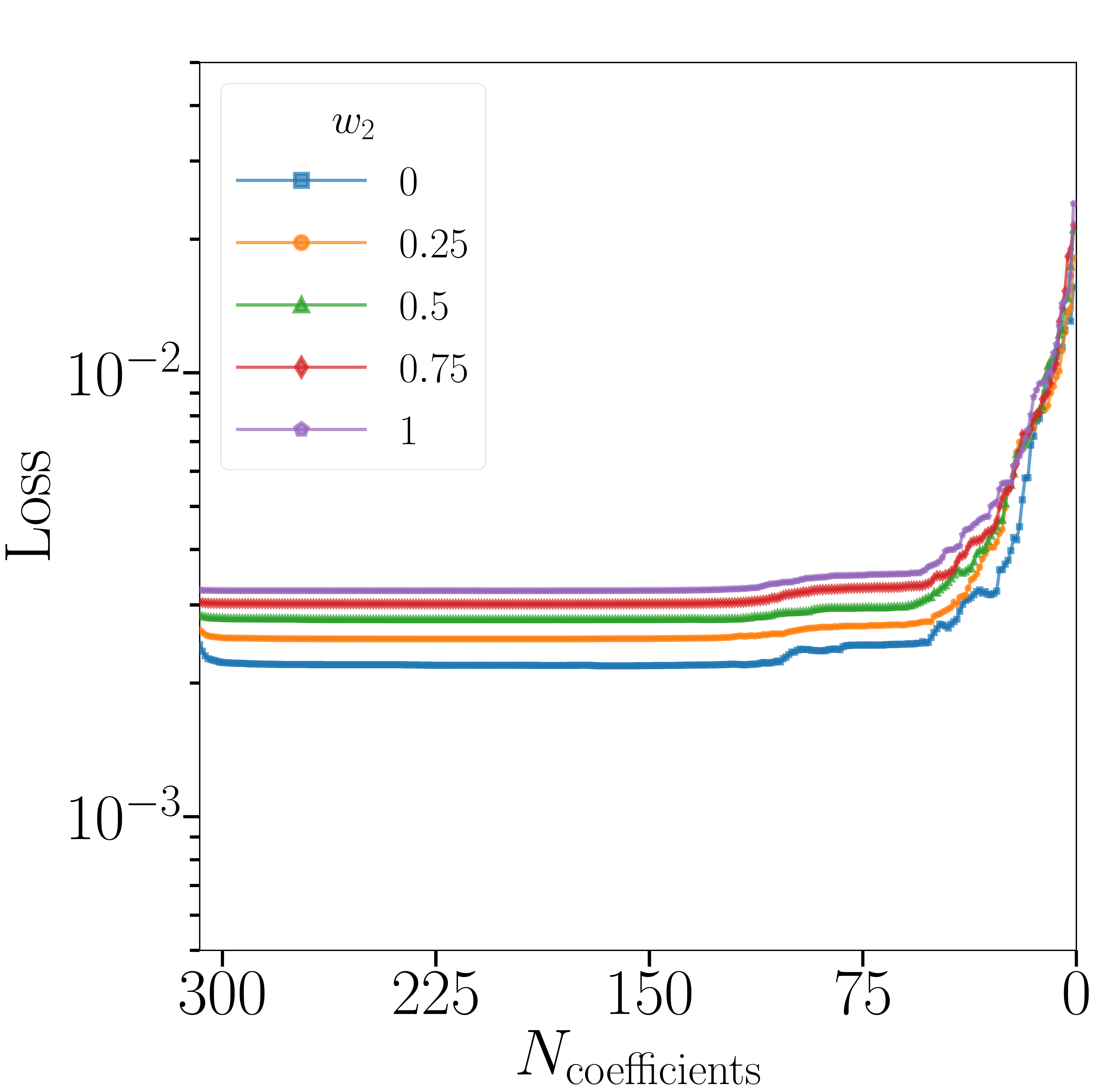}
	\subcaption{Ridge regression with $\lambda = 10^{-17}$.}
	\label{fig:vol_loss_Ridge17}
\end{subfigure}
\caption{Stepwise regression loss curves for the phase volume fractions first order dynamics, labeled by their weighting of the $\ell_{2}$ loss in the weighted residual loss function.}
\label{fig:vol_loss_OLS_Ridge17}
\end{figure}

\begin{figure}[hpt]
\centering
\begin{subfigure}[t]{0.49\textwidth}
	\centering
	\includegraphics[width=\textwidth]{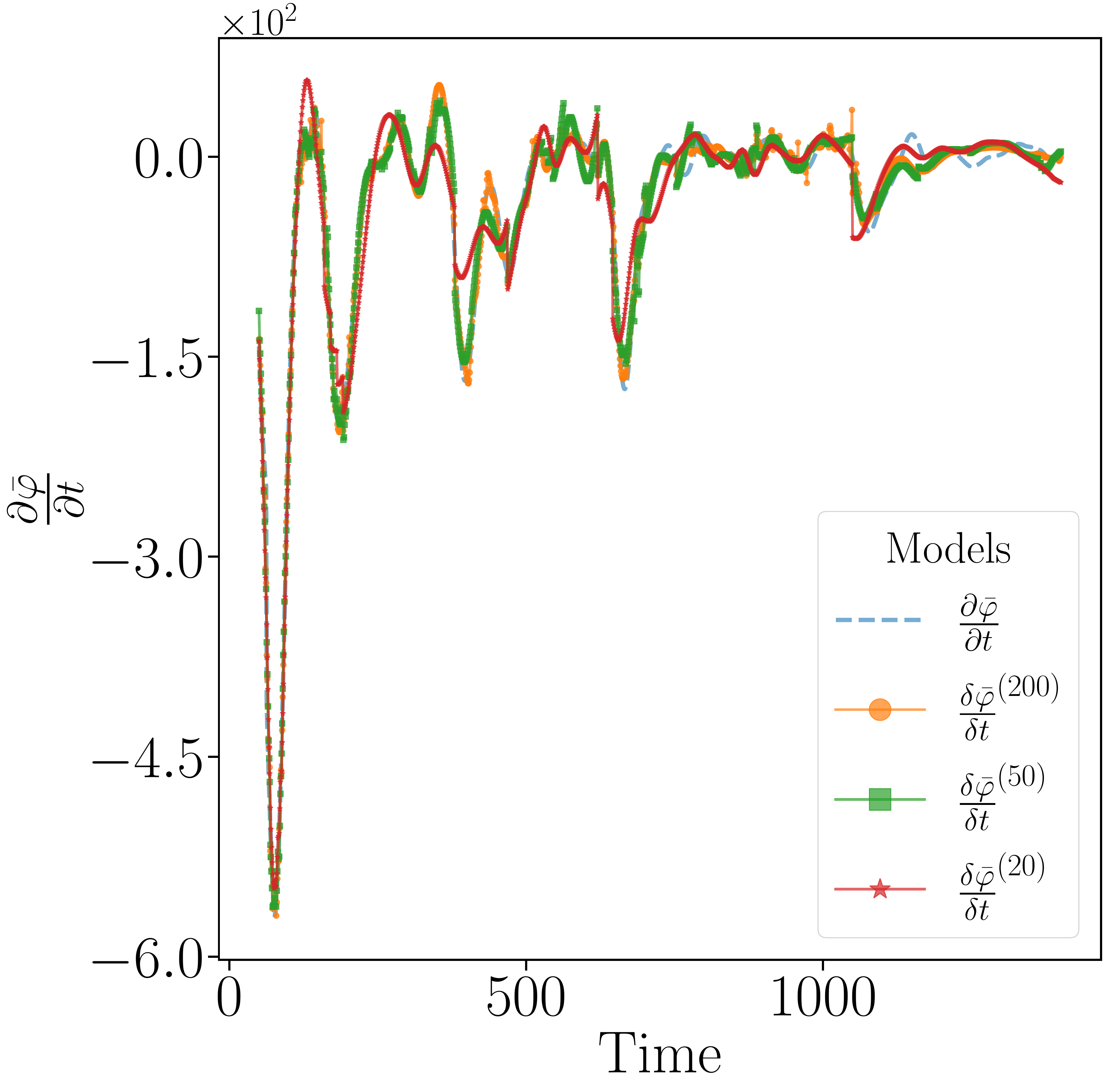}
	\subcaption{OLS regression.}
	\label{fig:vol_fit_OLS}
\end{subfigure}
\hfill
\begin{subfigure}[t]{0.49\textwidth}
	\centering
	\includegraphics[width=\textwidth]{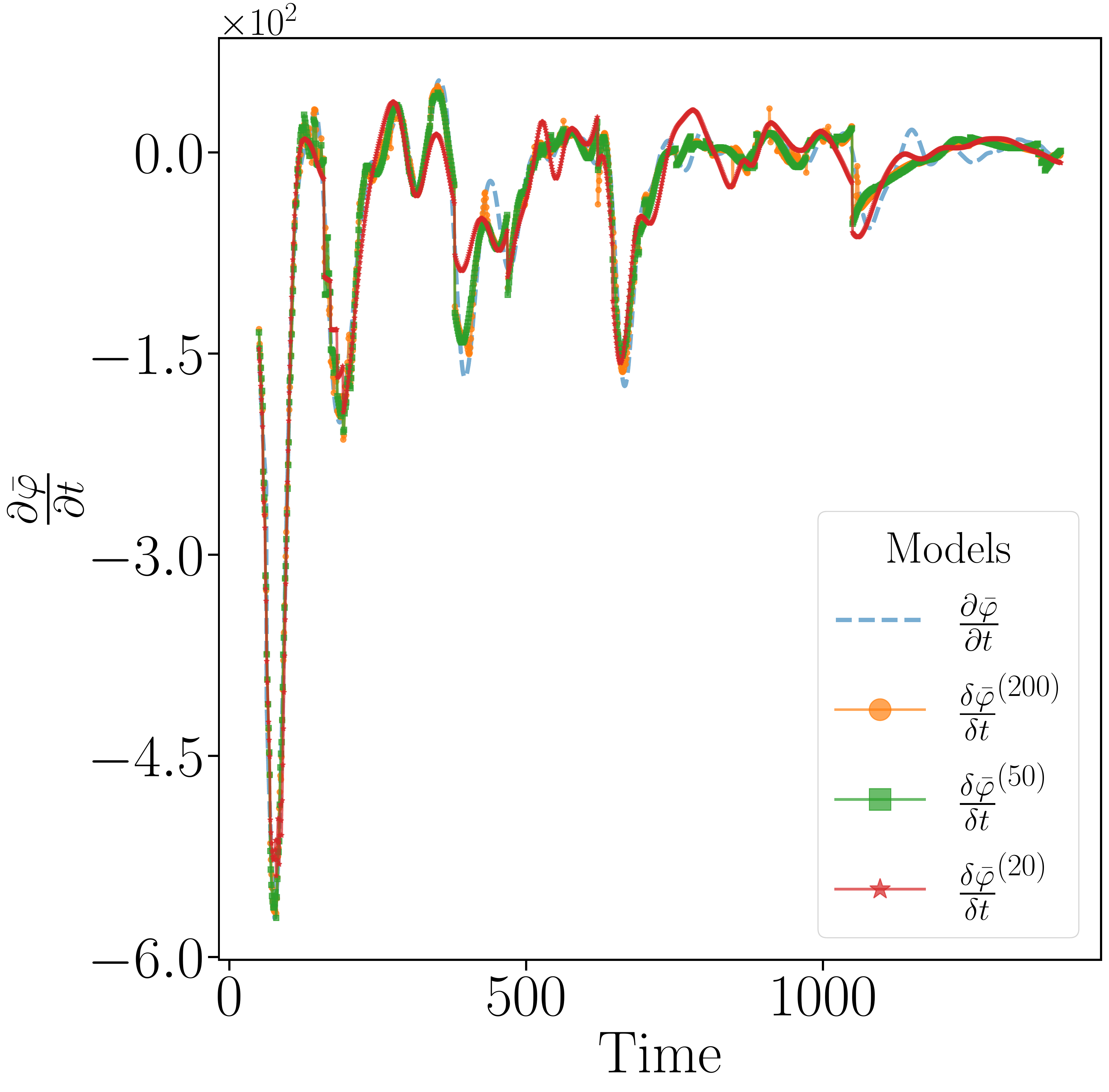}
	\subcaption{Ridge regression with $\lambda = 10^{-17}$.}
	\label{fig:vol_fit_Ridge17}
\end{subfigure}
\caption{First order dynamics of phase volume fractions for $20$, $50$, and $200$ terms in the forcing function, and use of the $\ell_{1}$ loss in stepwise regression. Backward Euler first time derivative data is shown with the dashed blue curve.}
\label{fig:vol_fit_OLS_Ridge17}
\end{figure}

When comparing the performance of the OLS and Ridge regression approaches, both approaches are shown to yield similar trends in the loss curves in \cref{fig:vol_loss_OLS_Ridge17}. However the Ridge regression loss is shown to be greater, by less than an order of magnitude for more complex models than the OLS regression. The resulting fits for the Ridge regression in \cref{fig:vol_fit_OLS_Ridge17} are also shown to capture the general trends of fluctuations quite well, however do not fit the sharper oscillations at later time quite as well as the OLS approach. This increased Ridge loss is reasonable considering the additional constraints, and is balanced by improving the smoothness of the loss curves over the stepwise iterations. The smoother Ridge curves subsequently have less plateaus, and the losses for different loss weights are more distinct. This approach appears to be less affected by changes in the rank of the matrix. Both approaches also behave similarly at the most parsimonious models with less than approximately $47$ terms. These comparable results indicate that Ridge regression is the most appropriate approach for future studies, increasing the robustness of the modelling to more collinear data, and still capturing the essential physical trends. Future cross validation studies, as well as the use of additional constraints on the coefficients in conjunction with the stepwise regression will be useful to determine how to best select the most parsimonious model.

The importance of choosing a physics informed basis is shown clearly by this example. Other choices of basis, such as including other derivatives, or a different set of basis variables like the other phase volume fractions, were shown to have much worse regression losses, and were much less intuitive. Future applications of this graph theoretic method should ensure the forms of the original models in the higher dimensional data are considered when determining how to capture effective physical processes.

\subsection{Free energy functional representation} \label{sec:physicalsystems_psi}
To model the free energy of the system, we are guided by the polynomial form of the free energy density in \cref{eq:psi_hom} and write 
\begin{align}
	\Psi = \Psi(\bar{\varphi},{\vecbar{E}}) \label{eq:model_psi}
\end{align}
as a function of the single chemical and ${d(d+1)}/{2}$ mechanical variables in $d$ spatial dimensions. 

\subsubsection{Proposed model for free energy functional representation}\label{sec:physicalsystems_psi_proposedmodel}

The form of this model suggests our functional representation should be a modified Taylor series at time $j$, about a base value at the time $i$:	
\begin{align}
	\Psi(\bar{\varphi}_{j},\vecbar{E}_{j}) = \Psi(\bar{\varphi}_{i},\vecbar{E}_{i})
	~&+~ \gamma^{\bar{\varphi}} \difference[1]{\Psi(\bar{\varphi}_{i},\vecbar{E}_{i})}{\bar{\varphi}} \Delta {{}{\bar{\varphi}}}_{ij} ~+~ \gamma^{\bar{E}_{\alpha \beta}} \difference[1]{\Psi(\bar{\varphi}_{i},\vecbar{E}_{i})}{\bar{E}_{\alpha \beta}} \Delta {{}\bar{E}_{\alpha \beta}}_{ij} \label{eq:model_psitaylor} \\
	~&+~ \gamma^{{\bar{\varphi} \bar{\varphi}}} \unidifference[2]{\Psi(\bar{\varphi}_{i},\vecbar{E}_{i})}{{\bar{\varphi}}} {\Delta {{}{\bar{\varphi}}}_{ij}}^2 \nonumber\\
	 ~&+~ \gamma^{{\bar{E}_{\alpha \beta}}{\bar{E}_{\xi \kappa}}} \difference[2]{\Psi(\bar{\varphi}_{i},\vecbar{E}_{i})}{{\bar{E}_{\alpha \beta}},{\bar{E}_{\xi \kappa}}} \Delta {{}{\bar{E}_{\alpha \beta}}}_{ij} \Delta {{}{\bar{E}_{\xi \kappa}}}_{ij} \nonumber\\
	 ~&+~ \gamma^{{\bar{\varphi}}{\bar{E}_{\alpha \beta}}} \difference[2]{\Psi(\bar{\varphi}_{i},\vecbar{E}_{i})}{{\bar{\varphi}},{\bar{E}_{\alpha \beta}}} \Delta {{}{\bar{\varphi}}}_{ij} \Delta {{}{\bar{E}_{\alpha \beta}}}_{ij} \nonumber \\
	 ~&+~ O(\Delta {{}{\bar{\varphi}}}_{ij}^3) + O(\Delta {{}{\bar{E}_{\alpha \beta}}}_{ij} ^3), \nonumber
\end{align}
where ${\Delta x_{ij}} = x_{j} - x_{i}$ represents the change in the $\varphi_{\alpha}$ or $\vecbar{E}$ state variable between times $j$ and $i$. Please refer \cref{app:implementation_modifiedtaylorseries} of the Appendix for details on how the Taylor series with fit linear coefficients $\gamma$ are formed. Similar to the model for the phase volume fraction evolution, we choose to restrict the dependencies to be $q=4$ variables of a single chemical and three mechanical order parameters in $d=2$ spatial dimensions. We further restrict the Taylor series to be up to $k = \ith[4]$ order, and there are $p = {(q^{k+1}-1)}/{(q-1)} = 341$ possible terms in this model. Here, the expansion is about the initial time $i=0$, and all partial derivatives in the Taylor series are approximated using the non-local calculus definitions in \cref{eq:diff_p}.

\subsubsection{Numerical implementation for free energy functional representation}\label{sec:physicalsystems_psi_numerical}
From initial leave-one-out cross-validation studies, it is found that any ridge parameter $\lambda > 0$ improves the rank of the matrix, and for this problem, decreases the fit losses compared to the OLS method. The fit losses also increase with increased $\lambda$. Therefore, $\lambda = 10^{-17}$ is fixed to perform Ridge regression. Similarly to for the phase volume fractions, the optimal direct solver for Ridge regression for this problem is found to be an SVD based solver, and the OLS regression is performed with a standard least squares solver. 

The Ridge regression loss curve in \cref{fig:psi_loss_OLS_Ridge17} is slightly smoother than the OLS regression, and unlike in the phase volume fractions problem, does not increase the loss with $\lambda > 0$. The regularization is shown to be able to increase the rank of the operators, and the robustness of the method, while still producing an accurate model. Here, the free energy is a monotonically decreasing function and is very smooth, so the loss function is the standard $\ell_{2}$ error.

\subsubsection{Reduced order model for free energy functional representation}\label{sec:physicalsystems_psi_finalmodel}
\begin{figure}[hpt]
\centering
\includegraphics[width=0.5\textwidth]{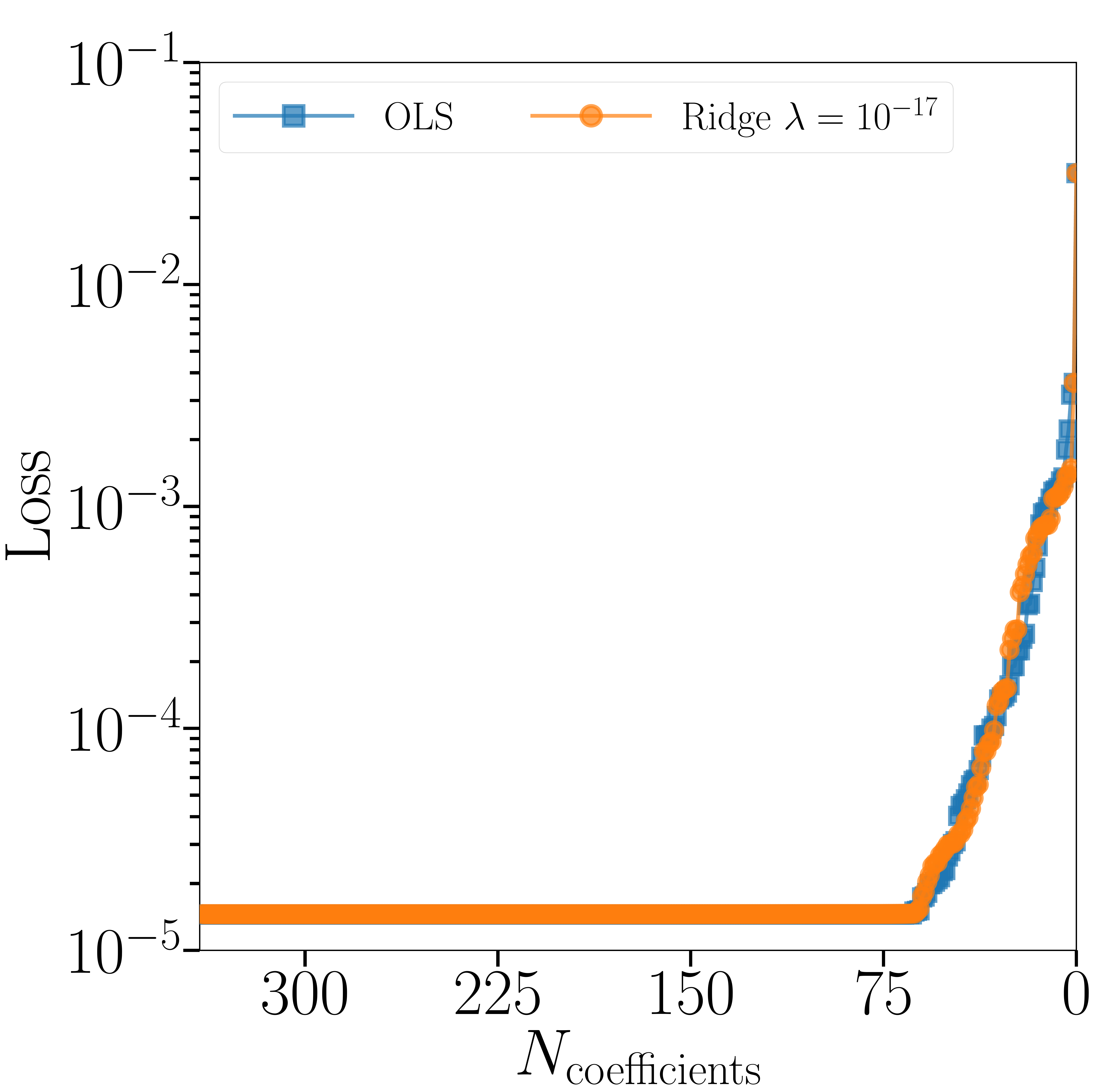}
\caption{Stepwise regression loss curve for the $\ith[4]$ order free energy Taylor series functional representation, using OLS and Ridge regression.}
\label{fig:psi_loss_OLS_Ridge17}
\end{figure}

From the loss curves in \cref{fig:psi_loss_OLS_Ridge17}, the sharp increase in the loss occurs at around $70$ terms, which corresponds to the rank of this basis of Taylor series terms. The loss increasing sharply for models with less than this number of terms indicates the importances of these terms in this truncated Taylor series. The loss for the full model is of order $O(10^{-5})$, comparable to other methods \cite{Zhang2020}. 

The first 10 terms from the stepwise regression for the $\ith[4]$ order Taylor series functional representation for $\Psi$ are found for OLS and Ridge regression to be:
\input{./figures/model_psi_OLS_Ridge17.tex}

\noindent Here, we have imposed for model consistency, that the coefficient associated with the base Taylor series term is $\gamma_0 = 1$. The most important higher order terms for both models are primarily second and third order derivatives, and both the OLS and Ridge models depend most strongly on the second derivative with respect to phase volume fraction. The OLS model has only one fourth order term, whereas the Ridge model has two fourth order terms, and a first order term. This possibly corresponds to the original free energy density $\psi$ in \cref{eq:psi_hom} being second, third and fourth order powers of the order parameters. However care must be taken on the interpretation of the results of the stepwise regression as being directly induced by the forms of the models for the higher dimensional data. In this example, the stepwise regression indicates a direct correspondence between the form of the reduced order model for $\Psi$, and the form of the spatially dependent $\psi$, however this correspondence should not always be assumed. Physics informed initial choices of basis will ensure the most appropriate reduced order model is found, and later interpretation is specific to each problem.

During the modelling process, it is observed that this Taylor series basis without any regularization is very singular, with rank $70 \ll p$, and condition number on the order of $O(10^{18})$. This can be attributed to the sparsity of the data, as well as potentially the form of the free energy, being a monotonically decreasing function, and derivatives with respect to various quantities being possible numerically linearly dependent. This singular nature of this basis means standard linear regression, computed through an SVD algorithm to solve the normal equations, as per \cref{app:implementation_linearregression} of the Appendix, sets certain singular values to zero, as per a given tolerance. The low rank of this matrix means many singular values are close to zero, and floating point issues related to the machine precision may occur.

\begin{figure}[hpt]
\centering
\begin{subfigure}[t]{0.49\textwidth}
	\centering
	\includegraphics[width=\textwidth]{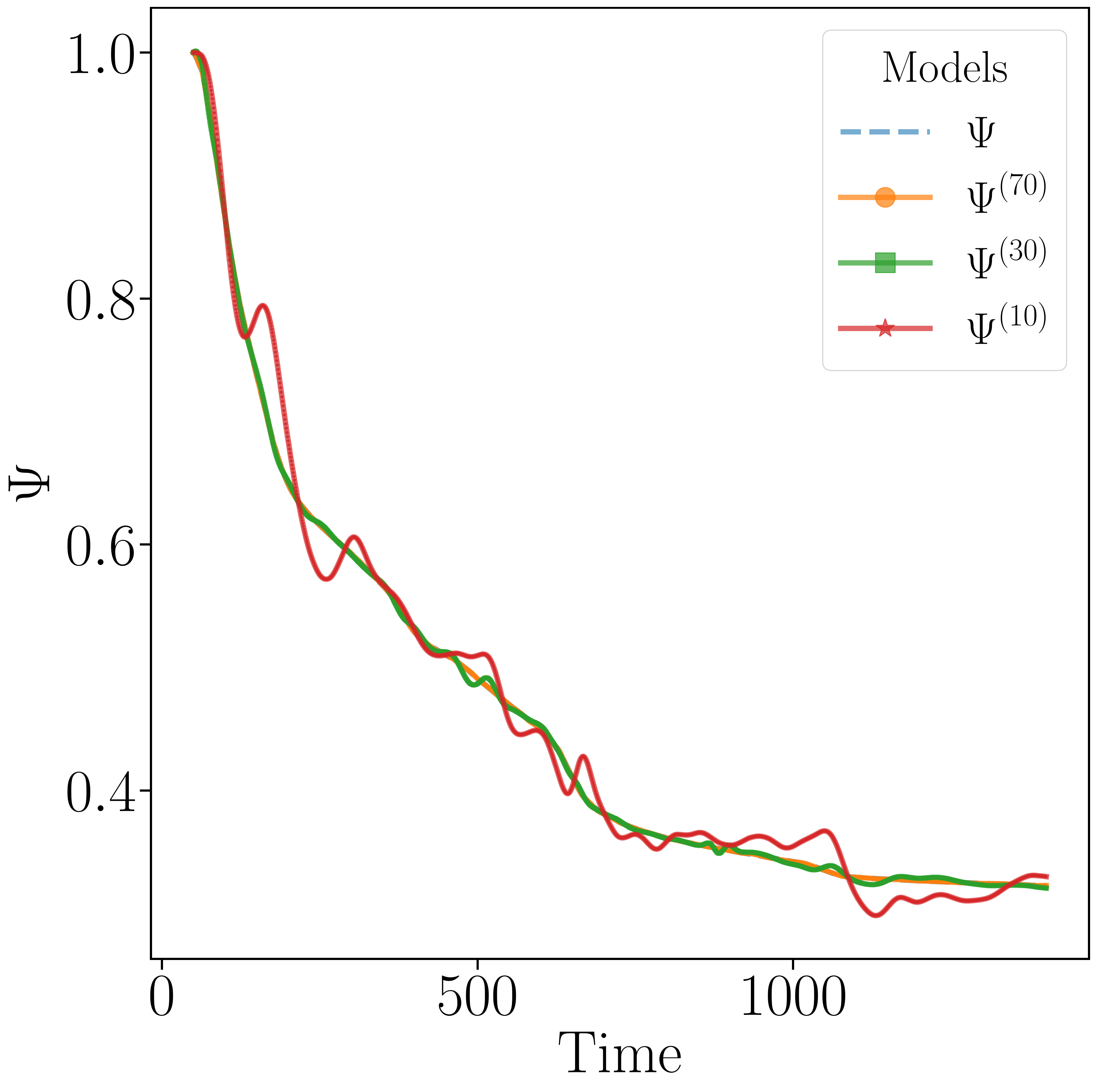}
	\subcaption{OLS regression.}
	\label{fig:psi_fit_OLS}
\end{subfigure}
\hfill
\begin{subfigure}[t]{0.49\textwidth}
	\centering
	\includegraphics[width=\textwidth]{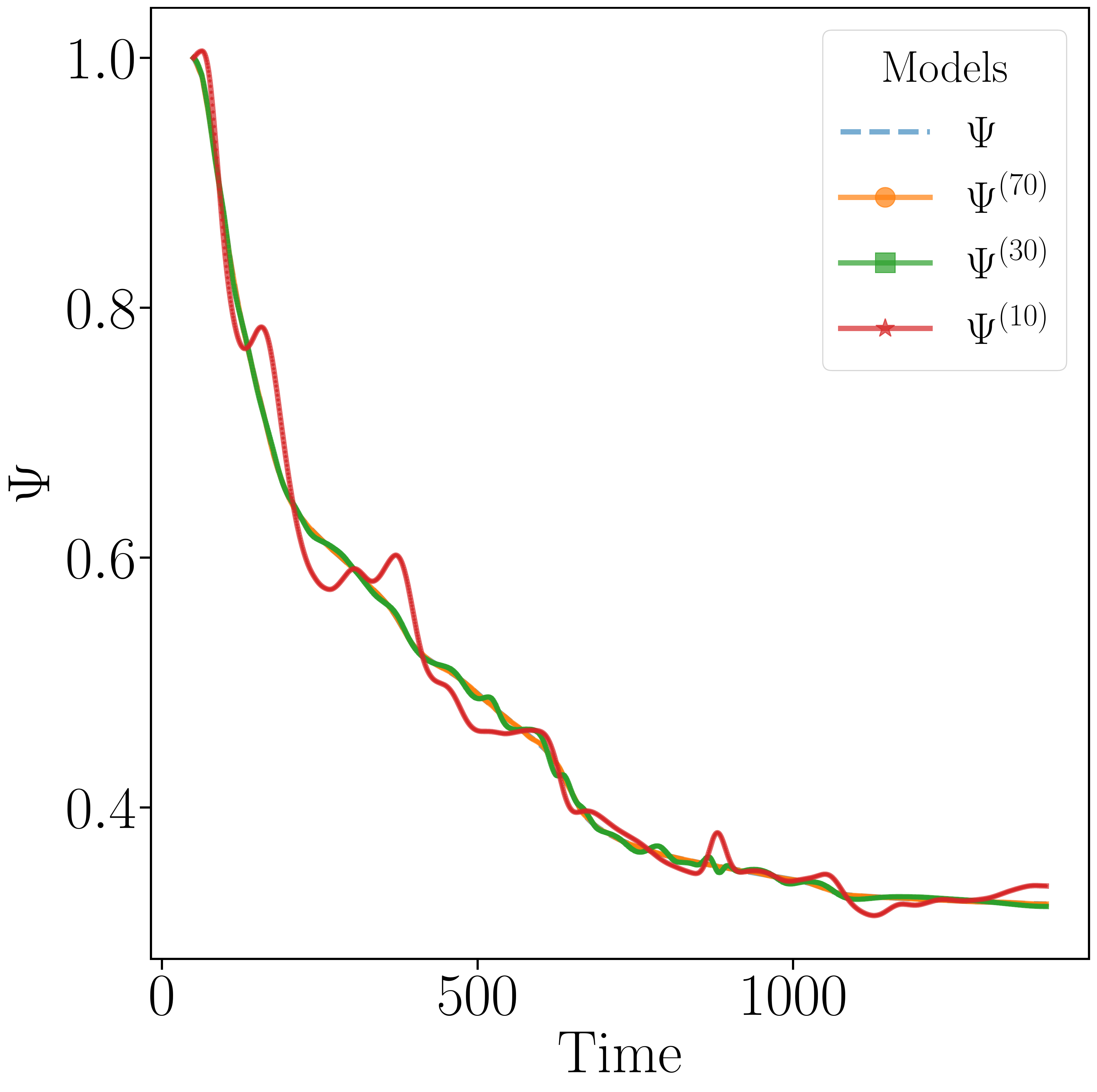}
	\subcaption{Ridge regression with $\lambda = 10^{-17}$.}
	\label{fig:psi_fit_Ridge17}
\end{subfigure}
\caption{Fitted curves for the $\ith[4]$ order free energy Taylor series functional representation, with $10$, $30$, and $70$ terms. DNS free energy data is shown with the dashed blue curve, and is essentially coincident with the higher order models.}
\label{fig:psi_fit_OLS_Ridge17}
\end{figure}

Regarding the fits for the free energy, Taylor series with $10$, $30$, and $70$ terms are shown in \cref{fig:psi_fit_OLS_Ridge17}. Models beyond the rank of this basis fit the free energy very well and model a very smooth function. More parsimonious models capture the overall decrease of the free energy, however oscillations in the phase volume fractions and strains are not able to be as smoothed out. Free energies at earlier times also appear to be fit better by all models. When comparing the OLS and Ridge fits, the fits are almost identical for more complex models, however the OLS method has smaller oscillations for the most parsimonious models. This Taylor series approach is shown to be a logical and effective basis of terms to model this smooth free energy function. 

\section{Conclusions}
\label{sec:concl}
In this work, we develop a graph theoretic approach for reduced order modelling of physical systems. After extending a non-local calculus as the foundation of this approach, we observe there are challenges and compromises when imposing consistency with differential calculus definitions of partial derivatives. Constraints placed on the graph's edge weights lead to numerical issues, but also interpretation issues if the non-local derivatives have error constants that are consistently large. A more in depth comparison of the behavior of non-local and differential derivatives is necessary in future works. 

Applying this calculus to studying volume averaged quantities from high dimensional direct numerical simulations by representing these quantities on a graph shows the intuitive, general, and effective nature of this approach. With an example physical system of the phase evolution of a mechano-chemical multi-crystalline solid, we see that using a physics informed basis of operators can model both differential equations, and functional representations. Highly oscillatory or smooth monotonically evolving functions can be fit using linear regression and an adequate number of terms, with errors on the order of $O(10^{-4})$. As a comparison, recent studies of the same microstructure system using neural networks \cite{Zhang2020} achieve errors on the order of $O(10^{-6})$.

We also note that near the beginning of the stepwise regression when the models contain the full set of operators, there is sometimes, for example in \cref{fig:vol_loss_OLS_Ridge17}, decrease in the loss function as the model becomes initially more parsimonious. This suggests that although we are choosing an intuitive, physics informed basis, there may be some operators that worsen the model, and must be rejected empirically. A more rigorous study to more precisely determine which operators from the appropriate basis should be included in the initial model will be conducted in future communications.

\textbf{Remark}: There remains an open question of how accurate these graph theoretic methods can be. This question may be related to the non-local nature of the derivatives, and whether they are appropriate to model strictly locally behaving systems. Whether or not these methods can achieve the error tolerances of other direct numerical simulations, or methods such as neural networks, will show the true effectiveness of this approach. An error analysis in \cref{app:error,app:error_p} also indicates that the error in a one dimensional linear Taylor series model only scales linearly with the sparsity of the data, compared to quadratically for standard interpolation methods, and further only scales sub-linearly for higher dimensional systems. This sub-optimal error scaling is directly attributed to the non-local calculus definitions and the form of the edge weights, and suggests modifications to these definitions may improve the scaling. 

There is also the compromise of the robustness and consistency of this modelling approach versus its accuracy. As discussed in \cref{app:implementation_linearregression} of the Appendix, the improvement of the multicollinearity of the basis of operators, and ensuring a consistent order of operators yielded by the stepwise regression, from regularization approaches may come at the cost of decreased loss. The physical systems studied indicate a physically appropriate basis of operators may still be numerically collinear, and regularization aids in the robustness and consistency of the modelling, without significant loss in accuracy.

We note, however, that in any reduced-order model, the goal is a parsimony of representation, which is achieved here by $O(10)$ components in the state vector, at the cost of some loss of fidelity. This approach appears to require $O(10^2)$ terms in order to get accurate fits of the data, however these higher order models involve higher order derivatives that may suffer numerical errors due to their large values. This compromise of model accuracy versus model complexity is another open question to be studied. To justify the requirements of these higher order models, we note that since the original PDE itself in the Cahn-Hilliard equation has nearly 20 terms when fully expanded, this number could be regarded as a lower bound for the complexity of the reduced-order models. Then, with the significant saving in dimensionality of the state vector over that of a DNS, even $O(10^2)$ operators in the reduced-order model is not unreasonable. This is especially reasonable when this representation is required to reproduce the dynamics of the vastly higher dimensional DNS.

Given these findings, this graph theoretic approach is shown to be very general, and many applications to a wide range of physical system, and system scales are possible. There are also several interesting studies that can be conducted in future works regarding the methods themselves. First it will be important to develop edge weight constraints for finite, discrete graphs, as opposed to current use of constraints for finite continuous graphs. Rigorous approximations for these discrete sums over the volume of the graph will possibly constrain the numerical issues brought on by the consistency demands, and will allow the local and non-local calculus to be better compared. The forms of the decaying edge weights can also be generalized to include several decay parameters that reflect independent behaviors of the state vector components. Finally, aspects of the stepwise regression, including the effects of different loss functions, and how the independence of operators in the basis plays a role in how and which operators are rejected from the model, will aid in the fitting process. \\

\section*{Acknowledgments}
\noindent MD and KG gratefully acknowledge the support of the National Science Foundation, United States grant \textbf{\#1729166}. "DMREF/GOALI: Integrated Framework for Design of Alloy-Oxide Structures".

\pagestyle{bibliography}
\bibliography{main}

 
\newpage
\begin{appendices} \label{sec:appendix}
\renewcommand\thesection{\Alph{section}}
\renewcommand\theequation{\getcurrentref{section}.\arabic{equation}}
\renewcommand\thefigure{\getcurrentref{section}.\arabic{figure}}



\section{Implementation and methodology}
\label{app:implementation}
\setcounter{equation}{0}
\setcounter{figure}{0}
Given some basis of possible algebraic and differential operators, system inference methods can be employed to identify the reduced order model with the desired combination of expressiveness and parsimony, similar to the procedures developed by Wang \etal \cite{Wang2019,Wang2020} and by Kochunas \etal \cite{Kochunas2020}. A graph theory library has been developed that implements this procedure, and performs regression using backwards stepwise regression to develop various representations for systems of interest.

\subsection{Stepwise regression} \label{app:implementation_stepwiseregression}
Backwards stepwise regression consists of iteratively performing regression on an \textit{a priori} given basis of $p$ operators, with model parameters $\gamma$, and a statistical criterion, the \textit{F-test}, to determine which operators are most, or least relevant to the model.
Operators can then be iteratively removed from the model, and further regression on this reduced basis is performed to confirm the (lack of) relevance of the operator. 
This procedure can be repeated until a minimal set of operators yields the most accurate model, and it can be posited that this minimal basis adequately represents the system of interest.
Physical insight also can aid in determining which operators will be relevant, and to help fine tune the statistical criteria to remove or keep operators.

In this work, we will be constrained to fitting models using linear regression, and it is assumed that there is a linear relationship between the operators $X \in \reals[n \times p]$ and outputs $y \in \reals[n \times d]$.
Here, $\gamma \in \reals[{p \times d}]$ is an array of $p$ coefficients and the problem can be posed as finding $\gamma$ such that: 
\begin{equation}
 y = X\gamma,
\end{equation}
where operators in $X$ are removed by fixing the associated coefficient for that operator to zero. All data is scaled by appropriate length scales to be dimensionless, and then further normalized over the $n$ data as
\begin{align}
	y = X\gamma ~&\to~ \tilde{y} = \tilde{X}\tilde{\gamma} \\
	\tilde{X} = X\alpha,\quad\tilde{y} ~=&~ y\beta,\quad\tilde{\gamma} = \alpha^{-1}{\gamma}\beta, \label{eq:normalizationconditions}
\end{align}
where $\alpha \in \reals[{p \times p}]$ and $\beta \in \reals[d \times d]$ are diagonal normalization matrices.

In addition to being scaled and normalized, particularly for this mechano-chemical example, where there is rapid oscillation between phases of the material over short time scales, the data must be filtered to smooth out rapid fluctuations, but retain all essential physics. A Gaussian filter with finite window and standard deviation is passed multiple times over the data.

The stepwise regression has several possible approaches that can be followed, depending how much of your total dataset you want to include in the training process. One approach is individual regression for each given dataset, using only data and operators for an individual, independent model. Another approach is composite regression using the concatenation of all datasets, where an optimal generalized model for all datasets is found by choosing operators based on their importance across all data. The composite method chooses the order to reject operators that optimize the loss across all groups at each iteration. The individual method may reject operators in such an order that the predicted loss at a given iteration is lower than the composite method loss at that same iteration, and the choice of approach is problem dependent.

\subsection{Linear regression} \label{app:implementation_linearregression}
The reduced order modelling followed in this work is centered around linear regression. Here, we initially solve the resulting normal equations of the normalized data, yielding an Ordinary Least Squares (OLS) loss function
\begin{align}
	\ell_{\textrm{OLS}} =& \left(y-X\gamma\right)^{T}\left(y-X\gamma\right)_{} \label{eq:olsloss} \\
	\intertext{to be minimized}
	\gamma_{\textrm{OLS}} =& \argmin_{\tilde{\gamma}} ~ \ell_{\textrm{OLS}}(\tilde{\gamma}). \label{eq:olslossmin}
\end{align}
This OLS regression can be solved for example via the SVD decomposition, yielding
\begin{align}
	\gamma_{\textrm{OLS}} ~=&~ \left(X^{T}X\right)^{-1}X^{T}y.  \label{eq:lqscoef}
\end{align}
There are several implementation issues that must be addressed concerning the operators in the data $X$, the rank of this matrix, and how this affects the fitting. For the given system of interest, where there are $n \sim O(10^{3})$ data points, and $p \sim O(10^{2})$ operators, it is found that this data matrix is generally very singular, with ranks on the order of $r < p$, and condition numbers on the order of $O(10^{18})$. The derivative terms appear to be very collinear for this dataset. This system is reaching equilibrium through monotonic minimization of the free energy, suggesting the state variable components evolve in a correlated fashion. This singular behavior means alternatives to ordinary least squares are necessary. We choose to pursue a general Tikhonov regularization scheme\cite{Lv2009,Tibshirani1996}, where the least squares loss function is supplemented with other quadratic constraints on the linear coefficients
\begin{align}
	d - C\gamma = 0,
\end{align}
where $C$ is a full-rank regularizing matrix of constraints, and $d$ is a biasing term. Non-linear constraints
\begin{align}
	c(\gamma) = 0,
\end{align}
may also be potentially added. The total loss function now takes the form
\begin{align}
	\ell ~=~&~~~~ \left(y-X\gamma\right)^T\left(y-X\gamma\right) \label{eq:constraintloss} \\ 
	~+&~ \lambda~\left(d - C\gamma\right)^{T}\left(d - C\gamma\right) \nonumber \\
	~+&~ \mu~c(\gamma). \nonumber
\end{align}

This resulting loss function is parameterized by the Lagrange multipliers $\lambda,\mu$, allowing tuning, or determining through cross validation or information criteria \cite{Wang2007}, the optimal ratio between the data and regularizing constraints that improves the matrix rank, and the resulting loss function. Cross validation studies are typically done through partitioning the dataset into training and testing data, and we focus on performing leave-out-out, or $k=n$ cross validation. Here, one of the $n$ data points is not used to fit the model, and instead is used to test the model trained using the $n-1$ points. The optimal tunable hyperparameters are determined from finding the minimum averaged loss over all possible partitionings of the data. There is freedom in the choice of the linear and non-linear constraints, allowing certain linear coefficients to be constrained more than others. Non-linear regularizing approaches such as Lasso \cite{Tibshirani1996} perform a form of feature selection by restricting the allowed magnitude of coefficients, and this can be incorporated into the general model selection offered by stepwise regression.

\subsubsection{Ridge regression} \label{app:implementation_linearregression_ridge}
In the case of generalized Ridge regression \cite{Lv2009}, the regularizing function is the $L_{2}$ norm of the coefficients and the total loss function has the form
\begin{align}
	\ell_{\textrm{Ridge}} ~=&~ \left(y-X\gamma\right)^T\left(y-X\gamma\right) + \lambda~\gamma^{T}C^{T}C\gamma\label{eq:ridgeloss}
\end{align}
where there is zero bias. For this work we perform standard Ridge regression where $C=I$. This loss function is strictly quadratic in $\gamma$, yielding a closed form expression for the optimal coefficients
\begin{align}
	\gamma_{\textrm{Ridge}} ~=&~ \left(X^{T}X + \lambda C^{T}C\right)^{-1}X^{T}y. \label{eq:ridgecoef}
\end{align}
This solution is found from minimizing \cref{eq:ridgeloss}, and ensuring this is a minimum by finding the second derivative Hessian of the loss function
\begin{align}
\derivative[2]{\ell_{\textrm{Ridge}}}{\gamma^{\mu},\gamma^{\nu}} = 2(X^{T}X)_{\mu\nu} + 2\lambda(C^{T}C)_{\mu\nu} >0. \label{eq:ridgehessian}
\end{align}
The addition of $\lambda > 0$ results in the matrix $X^{T}X$ being more diagonally dominant, and thus higher rank, yielding a better pseudoinverse approximation for the coefficients. However, there is a bias-variance trade-off \cite{Geurts2009} of satisfying the additional constraints and improving the coefficient estimates, but incurring greater losses.

\subsubsection{Lasso regression} \label{app:implementation_linearregression_lasso}
In the case of generalized Lasso regression \cite{Tibshirani1996}, the regularizing function is the constrained $L_{1}$ norm of the linearly transformed coefficients and the total loss function has the form
\begin{align}
	\ell_{\textrm{Lasso}} ~=&~ \left(y-X\gamma\right)^T\left(y-X\gamma\right) + \lambda~\sum_{\mu} \abs{C_{\mu\nu}\gamma^{\nu}}\label{eq:lassoloss}
\end{align}
where there is zero bias. For this work we perform standard Lasso regression where $C=I$. This loss function is non-linear and has no closed form expression for the optimal coefficients, and a coordinate descent algorithm is used. The use of the $L_{1}$ means convex optimization can still be performed. We can approach $L_{1}$ by a sequence of norms that come $\epsilon$ close to $L_{1}$, but the limit lies outside this sequence. We can see that this approach to $L_{1}$ gives us greater regularization around $\gamma = 0$.

The second derivative Hessian of the loss function is not defined at $C_{\mu\nu}\gamma^{\nu} = 0$, but can be written as a generalized function using sub-gradients \cite{Clarke1983}.

The sub-gradient of $u(x)$ is defined as the set of all real vectors $v \in \reals[p] $ such that at a non-differentiable point $x_{0}$
\begin{align}
	u(x)-u(x_{0}) \geq&~ \dotproduct[v][{(x-x_{0})}] \nonumber \\ 
	\intertext{where in the case of $p=1$ dimension, the scalar $v$ lies within the interval}
	v \in&~ [a,b], \nonumber
	\intertext{such that $a,b$ are the left and right side derivative limits:}
	a,b =& \lim_{x \to x^{\mp}_{0}} \frac{u(x)-u(x_{0})}{x-x_{0}}.
\end{align}
In the case of the absolute value function used in Lasso regression, $v_{\abs{x}} \in [-1,1]$, and we choose the sub-gradient to be $v_{\abs{x}} = \sign[{x}]$, the vector of the signs of the components of $x$.

The first derivative of the Lasso loss function follows by representing the sign function $\sign[x] = 2H(x)-1 $ in terms of the Heaviside distribution $H(x)$, and invoking the theory of distributions
\begin{align}
\derivative[1]{{\ell_\textrm{Lasso}}}{\gamma^{\mu}} &= 
	2 (X^\text{T}X)_{\mu\nu}\gamma^{\nu} - 2 X_{\nu\mu}y^{\nu} + \lambda \sum\limits_\nu \left(2H(C_{\nu\theta}\gamma^{\theta}) -1 \right)C_{\nu\mu}.
	\label{eq:lassojacobian}
\end{align}
The Hessian follows as
\begin{align}
	\derivative[2]{{\ell_\textrm{Lasso}}}{\gamma^{\mu},\gamma^{\nu}} &= 2(X^\text{T}X)_{\mu\nu} + \lambda\sum\limits_\eta 2\delta(C_{\eta\theta}\gamma^\theta)C_{\eta\mu}C_{\eta\nu}, \label{eq:lassohessian}
\end{align}
where $\delta(x)$ is the Dirac-delta distribution.

To compare the Ridge and Lasso approaches, we suppose the data is completely orthonormal and $C_{\mu\nu} = \delta_{\mu\nu}$ is chosen. Then $X^{T}X = I$, and the optimal coefficients can be shown to be \cite{Tibshirani1996}
\begin{align}
	\gamma_{\textrm{Ridge}} ~=&~ \frac{1}{1+\lambda}\gamma_{\textrm{OLS}}, \label{eq:ridgecoeforthogonal} \\
	\gamma_{\textrm{Lasso}} ~=&~ \max \left(0,1+\frac{\lambda}{2\abs{\gamma_{\textrm{OLS}}}}\right) \gamma_{\textrm{OLS}}. \label{eq:lassocoeforthogonal}
\end{align}
Here, we see that Ridge regression scales the OLS coefficients, and Lasso regression either sets the coefficients to be proportional to $\lambda$, or sets them as identically zero.

\subsection{Modified Taylor series} \label{app:implementation_modifiedtaylorseries}
In this work, we often choose to specify the functional representation as a truncated Taylor series expanded around a base point $x_i$. If the partial derivatives were able to be exactly evaluated, this functional representation would be accurate up to higher-order terms, however as discussed, the non-local calculus derivatives are not exact. Their approximation to the differential Taylor series can be improved by introducing the coefficients $\gamma^{\mu_{1}\mu_{2} \dots \mu_{k}}$ for the $\ith[k]$ order derivatives
\begin{align}
	u(x_{j}) \approx u(x_{i}) ~&+~ \sum_{\mu}\gamma^{\mu} \difference[1]{u(x_{i})}{x^{\mu}} \Delta x^{\mu}_{ij} \label{eq:Taylor_approx} \\ 
	&+~\frac{1}{2}\sum_{\mu\nu}\gamma^{\mu\nu} \difference[2]{u(x_{i})}{x^{\mu},x^{\nu}} \Delta x^{\mu}_{ij} \Delta x^{\nu}_{ij} ~+~ O(\Delta x^{\mu}_{ij} \Delta x^{\nu}_{ij} \Delta x^{\theta}_{ij}), \nonumber
\end{align}
where $\Delta x^{\mu}_{ij} = x^{\mu}_{j} - x^{\mu}_{i}$ represents the change in the state variable. 

We aim to ensure the consistency of this model, and demand that the model in the limit 
\begin{align}
\lim_{\Delta x_{ij} \to 0} u(x_{j}) \to u(x_{i}). \label{eq:taylorseries_consistency}
\end{align} 
We therefore fix the base coefficient for $u(x_{i})$ to be identically 1, as per the Taylor series from differential calculus and let the higher order derivative operators be fit using regression.

Care must also be taken when scaling and normalizing the data as per \cref{eq:normalizationconditions} such that the fixed and free coefficients are scaled appropriately such that the base coefficient is unity in the unnormalized basis and the method remains consistent as per \cref{eq:taylorseries_consistency}. Let the operators be partitioned into $X = [\bar{X},\hat{X}]$ with fixed operators $\bar{X}$ and free operators $\hat{X}$. The fixed operators is generally data at fixed vertices $x_{i}$, such as the base vertex for a Taylor series representation, and the free operators are dependent in changes in data between vertices $\Delta x_{ij}$. The linear coefficients and normalizing matrices are also partitioned into free and fixed components such that
\begin{align}
	y = \bar{X}{\bar{\alpha}}{\tilde{\bar{\gamma}}}{\bar{\beta}}^{-1} +&~ \hat{X}\hat{\alpha}{\tilde{\hat{\gamma}}}{\hat{\beta}}^{-1}, \\
	\intertext{and}
	\lim_{\Delta x \to 0} y \to&~ \bar{X}\bar{\gamma}.
\end{align}

If for example, we are given a proposed Taylor series model fit with data $X,y$ with normalizing matrices $\alpha,\beta$ and new data $X^{\prime},~y^{\prime}$ with normalizing matrices $\alpha^{\prime},\beta^{\prime}$, denoted with primed-variables. The coefficients from the original fit $\gamma$ applied to the new data are normalized such that
\begin{align}
	\gamma^{\prime} ~=&~ \gamma \\
	\tilde{\gamma^{\prime}} ~=&~ \alpha^{\prime -1}\alpha\tilde{\gamma}\beta^{-1}\beta^{\prime}	
\end{align}
to ensure 
\begin{align}
	\lim_{\Delta x^{\prime} \to 0} y^{\prime} = \bar{X^{\prime}}\bar{\alpha^{\prime}}\tilde{\bar{\gamma^{\prime}}}\bar{\beta^{\prime}}^{-1} + \hat{X^{\prime}}\hat{\alpha^{\prime}}\tilde{\hat{\gamma^{\prime}}}{\beta^{\prime -1}} \to \bar{X}\bar{\alpha^{\prime}}\tilde{\bar{\gamma^{\prime}}}\bar{\beta^{\prime}}^{-1} = \bar{X}\bar{\gamma} = y,
\end{align}
given we are fitting with the normalized data.

The linear coefficients coefficients for the higher order derivatives $\gamma^{\mu}$, in the case of infinite knowledge of the system, should be identically $\gamma = 1$ if the Taylor series holds exactly. However, as discussed in \cref{sec:graphtheory_consistency}, the non-local partial derivatives may behave fundamentally differently from their differential counterparts, and the resulting Taylor series may require $\gamma^{\mu} \neq 1$ to accurately model the system of interest. A verification of how well the modified Taylor series holds can be conducted by observing which coefficients are close to $1$.

\subsection{Graph theory library} \label{app:implementation_graphtheorylibrary}
A graph theory library has been implemented in the \textit{python} language that represents graphs using the \textit{pandas} dataFrame structure \cite{mckinney-proc-scipy-2010}, allowing for efficient labelling, indexing, and appending of the graph quantities. All settings are contained in \textit{JSON} files, and data is imported into the dataFrame structure, either through user-defined or default importers, that can parse \textit{txt, csv, pickle, json, hdf5} and \textit{mat} files. To calculate operators on the graph, non-local calculus definitions from Gilboa \etal \cite{Gilboa2008} are implemented through vectorized and broadcasting operations using the \textit{numpy} data structure library \cite{NumPy-Array2020}, ensuring efficient sum calculations over the whole graph. 

The library is currently implemented to perform fitting using linear stepwise regression \cite{Kochunas2020}, which means the models are of the form $y = X\gamma$. The linear coefficients $\gamma$, and which basis members are included in the final model are found using an estimator class that inherits and generalizes fit, predict, loss, and cross validation modules from the \textit{scikit-learn} \cite{Pedregosa11a2011} libraries. Several direct linear solvers are available to solve the linear equations from the \textit{numpy} library that are based on LAPACK, including a direct pseudo-inverse, a least squares solver, and an SVD based solver. The library allows the solver to be specified by the user, and due to the different solvers yielding different results for very singular systems, the choice of solver is problem dependent. Additional solvers have been written that use these standard solvers to solve variants of the linear normal equations. For example, when performing Ridge regression, an SVD solver is most appropriate, due to there being a closed form expression for how to modify the singular value decomposition of $X$ with the ridge parameter $\lambda$:
\begin{align}
	X =&~ U \Sigma V^T \to U \frac{\Sigma}{\Sigma^2 + \lambda}V^T.
\end{align}
Here, $U$ and $V$ are the orthogonal matrices of right and left singular vectors, and $\Sigma$ are singular values of $X$. Models that are currently implemented include functional representations of left-hand-sides by either a monomial, polynomial, Chebyshev, Lagrange, or Hermite bases for the right-hand-side, of the form $y = X\gamma$ or by a Taylor Series \cite{Kochunas2020}, up to $\ith[k]$ order. These bases can also be used to form models for differential equations, such as first order dynamics $\derivative[1]{u}{t} = f(u,x,t)$. 

\newpage
\section{Gauss circle problem for discrete high dimensional manifolds} \label{app:implementation_gauss}
\setcounter{equation}{0}
\setcounter{figure}{0}
Given we know how to compute the weight constraints for a continuous graph with a finite radius, the problem must be translated to computing the weight constraints for a discrete graph. This problem is related to a well studied, but still open problem posed by Gauss \cite{Lowry-Duda2011,Ahmed2018}: that of finding the number of discrete points $N(R,p)$ that lie within a $p$-dimensional sphere of radius $R$. This is a counting problem, for how many sets of integers in $\mathbb{Z}^p$ have sums of their squares less than $R$. There are no closed form expressions for $p>2$, where for $p=2$, as per \cref{fig:gausscircle},
\begin{align}
 N(R,p=2) = 1 + 2^p\sum_{i=0}^{\floor{R}} \floor{\sqrt{R^2 - i^2}}.
\end{align}
However there are estimates for higher dimensions based on errors between the continuous volume, and the true value
\begin{figure}[ht]
	\centering
	\includegraphics[width=0.6\textwidth]{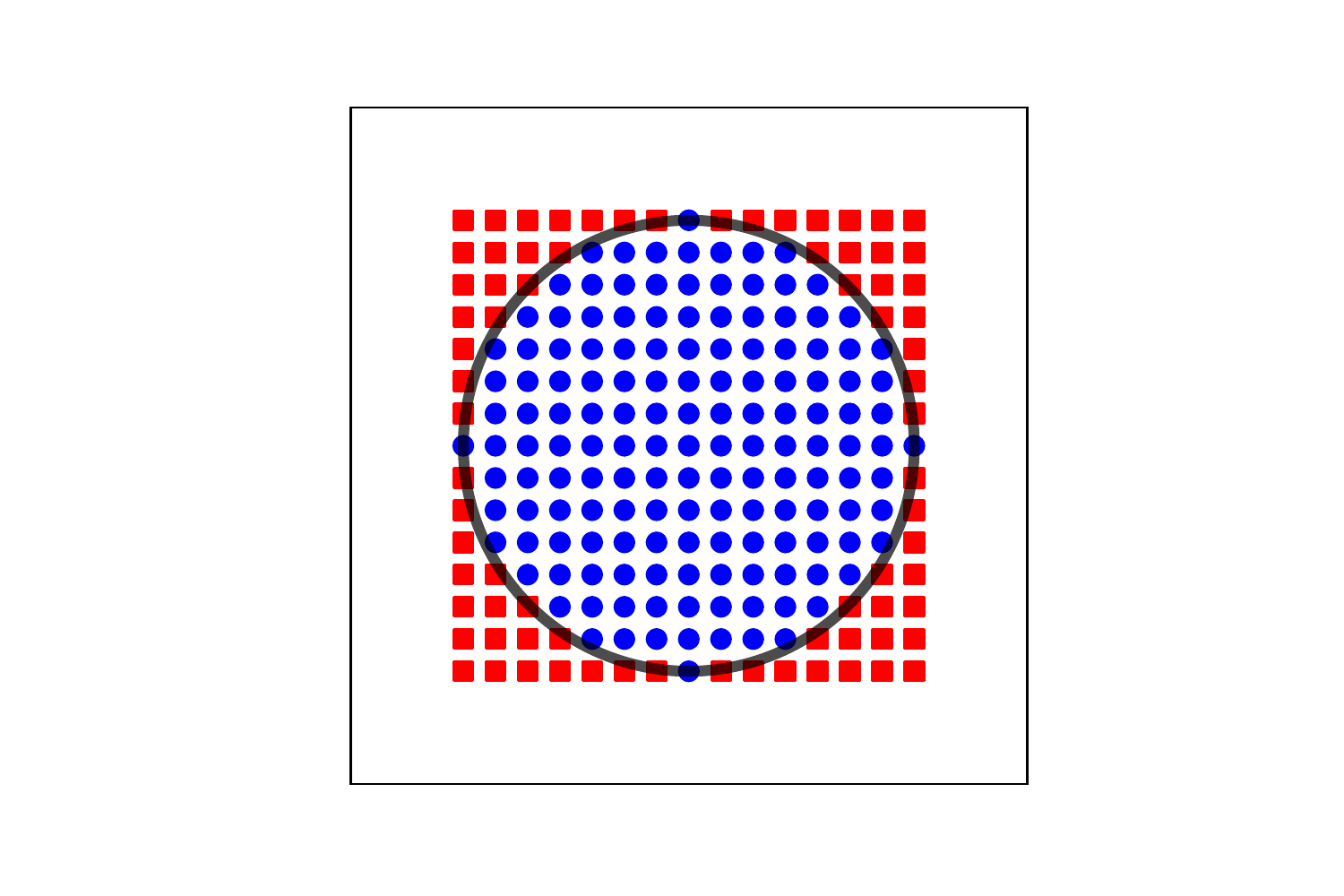}
	\caption{Gauss Circle problem in $d=2$ dimensions, of finding the number of blue points, within the radius of the circle on the grid of red squares.}
	\label{fig:gausscircle}
\end{figure}
\begin{align}
 N(R,p) = V_p + E(R,p).
\end{align}
These errors have been shown to be $E(R,p) = O(R^{\epsilon})$, where $\epsilon = \epsilon_p < 1$. In this work, we use the continuous volume weight constraints given by \cref{eq:constraint_scaled} and \cref{tab:weights}. Improved versions derived from the discrete $p$-dimensional topology will be used in future communications.


\newpage
\section{Error analysis of modified Taylor series for polynomial functions in one dimension}\label{app:error}
\setcounter{equation}{0}
\setcounter{figure}{0}
\subsection{Modified Taylor series model} \label{app:error_model}
We can also study the convergence of the total error $e$ of these methods as the data becomes more finely spaced, with spacing $h$. An error analysis will be conducted for a function $u(x)$ of a $p=1$ dimensional variable $x$, represented by a modified $k$ order Taylor series $u_k(x)$ with non-local derivatives. 
These derivatives are assumed to use polynomial weights, with super-quadratic scaling $\epsilon$.
\noindent For $p=1$, the polynomial weights take the following form for $0 \leq \epsilon < p$:
\begin{align}
	w(x_i,x_j) = C_{\epsilon}\frac{1}{\abs{x_j-x_i}^{2+\epsilon}},
\end{align}
where for the continuous weight constraints, $C_{\epsilon} = (p-\epsilon)R^{\epsilon}$ and $\lim_{\epsilon \to 0} C_{\epsilon} \to p$.

The Taylor series functional representations of the $K$ order function $u(x)$ and $k\leq K$ order model $u_k(x)$, based at a point $\tilde{x}$, are
\begin{align}
	u(x|\tilde{x}) =& u({}\tilde{x}) + \uniderivative[1]{u(\tilde{x})}{x}(x-\tilde{x}) + \frac{1}{2!}\uniderivative[2]{u(\tilde{x})}{x}(x-\tilde{x})^2 + \cdots + \frac{1}{K!}\uninderivative[K]{u(\tilde{x})}{x}(x-\tilde{x})^K, \label{eq:differentialtaylor}\\
	u_k(x|\tilde{x}) =& u(\tilde{x}) + \gamma_{1}(\tilde{x})\unidifference[1]{u(\tilde{x})}{x}(x-\tilde{x}) + \frac{\gamma_{2}(\tilde{x})}{2!}\unidifference[2]{u(\tilde{x})}{x}(x-\tilde{x})^2 + \cdots + \frac{\gamma_{k}(\tilde{x})}{k!}\unindifference[k]{u(\tilde{x})}{x}(x-\tilde{x})^k, \label{eq:modeltaylor}	
\end{align}
where the coefficients $\gamma_{k}(\tilde{x})$ are the linear coefficients fit to the non-local derivative model based around $\tilde{x}$. These linear coefficients depend on the training data used to construct the model. It is assumed that the model is exact at any base point, and so $\gamma_0(\tilde{x}) = 1$ in this unnormalized basis. The coefficients are appropriately scaled $\gamma_l \to \tilde{\gamma}_l$ when the operators are normalized. 

\noindent We must now determine, using the $l$-norm, the total error scaling in powers of $h$ of this $k$ order numerical method for a $K$ order function
\begin{align}
	e_k = \norm{e_k(x)}_l = O(h^{Q[K,k,p]})
\end{align}
based on the mesh used for the training and testing data with length scale $h$, and the pointwise error for the model, defined as
\begin{align}
	e_k(x) =&~ u_k(x) - u(x) = O(h^{q[K,k,p]}).
\end{align}
Based on the functional representations in \cref{eq:differentialtaylor} and \cref{eq:modeltaylor}, the pointwise error will be comprised of the error in the non-local derivatives
\begin{align}
	\unindifference[k]{u(\tilde{x})}{x} =&~ \uninderivative[k]{u(\tilde{x})}{x} + O(h^{r[K,k,p]}), \\
	\intertext{and the scaling of the fit linear coefficients}
	\gamma_{k}(\tilde{x}) =&~ 1 + O(h^{s[K,k,p]}).
\end{align}
The form of the derivatives and coefficients are based on the premise that the non-local calculus modified Taylor series will converge to the differential calculus Taylor series in the limit of infinite data. The partial derivatives will converge to the differential derivatives, and the linear coefficients will converge to unity.
Here the scaling powers $Q,q,r,s$ depend on the order of the function $K$, the order of the model $k$, and the number of dimensions $p$. 

\subsection{Integral summation definitions} \label{app:error_sums}
To conduct the error analysis, we define some quantities related to discrete sums of integers.

\noindent The difference of powers of two integers $i$ and $j$ can be defined using the relation
\begin{align}
	j^l - i^l =&~ (j-i)\sum_{q=1}^{l}i^{q-1}j^{l-q} = (j-i)p_l[i,j], \\
	\intertext{where the polynomial functions are defined as}
	p_l[i,j] \equiv&~\sum_{q=1}^{l}i^{q-1}j^{l-q},
	\intertext{and for example:}
	p_0[i,j] =&~ 0, \\
	p_1[i,j] =&~ 1, \\
	p_2[i,j] =&~ j+i, \\
	p_3[i,j] =&~ j^2 + ji + i^2.
\end{align}

Sums over positive integer powers $l\geq 0$ of integers are defined as the Faulhaber functions
\begin{align}
	\varphi_l[n] \equiv \sum_{j=0}^{n} j^{l} =&~ \frac{1}{l+1}\sum_{j=0}^{l}(-1)^{j}\binom{l+1}{j}B_jn^{l+1-j}, 
\end{align}
where $B_j$ are the Bernoulli numbers. For example,
\begin{align}
	\varphi_0[n] =&~ n, \\
	\varphi_1[n] =&~ \frac{n(n+1)}{2}, \\
	\varphi_2[n] =&~ \frac{n(n+1)(2n+1)}{6}, \\
	\varphi_3[n] =&~ \frac{n^2(n+1)^2}{4}.
\end{align}
We may also consider the normalized sums over powers of integers that do not scale to leading order with $n$:
\begin{align}
	\tilde{\varphi}_l[n] \equiv \frac{1}{n^{l+1}}\varphi_l[n],
\end{align}
where in the limit
\begin{align}
	\lim_{n \to \infty} \tilde{\varphi}_l[n] \to \frac{1}{l+1} + \frac{1}{2}\frac{1}{n} + O(n^{-2}).
\end{align}


\noindent Sums over negative real powers of the integers are defined as the Harmonic numbers
\begin{align}
	H_{\epsilon}[n] =&~ \sum_{j=1}^{n}\frac{1}{j^{\epsilon}},
	\intertext{where in the limits}
	\lim_{\epsilon \to 0} H_{\epsilon}[n] &\to n, \\
	\intertext{and} 
	\lim_{\epsilon \to 1,~n\to \infty} H_{\epsilon}[n] &\to \log{n} + \gamma_{\textrm{Euler-Mascheroni}}.
\end{align}

We will now discuss the appropriate functional representations, the data mesh used for training and testing the model, and the error analysis calculations.

\subsection{Functional representation} \label{app:error_functional}
Assume the function to be represented has a polynomial form with coefficients $\alpha_{l}$ up to $K$ order
\begin{align}
	u(x) =&~ \sum_{l=0}^{K}\alpha_{l}x^{l}.
	\intertext{The derivatives of the polynomial function are}
	\uninderivative[l]{u(x)}{x} =&~ \sum_{q=l}^{K} \frac{q!}{(q-l)!}\alpha_{q}{x}^{q-l} \label{eq:polyderivatives} \\
	\intertext{and so the Taylor series representation in \cref{eq:differentialtaylor} has the form}
	u(x|\tilde{x}) =&~ \sum_{l=0}^{K} \sum_{q=l}^{K} \binom{q}{l}\alpha_{q}{\tilde{x}}^{q-l}(x-\tilde{x})^l. \label{eq:polytaylorseries}
\end{align}

\noindent Sums over products of a polynomial function with a monomial of order $t$ can be written as
\begin{align}
	\sum_{j\in \tilde{V}} u(\tilde{x}_j)\tilde{x}_j^{t} =&~ \sum_{l=0}^{K}\alpha_{l} (2h)^{l+t}\varphi_{l+t}[n]. 
\end{align}

\subsection{Data mesh} \label{app:error_mesh}
\begin{figure}[hpt]
\centering
\includegraphics[width=0.4\textwidth]{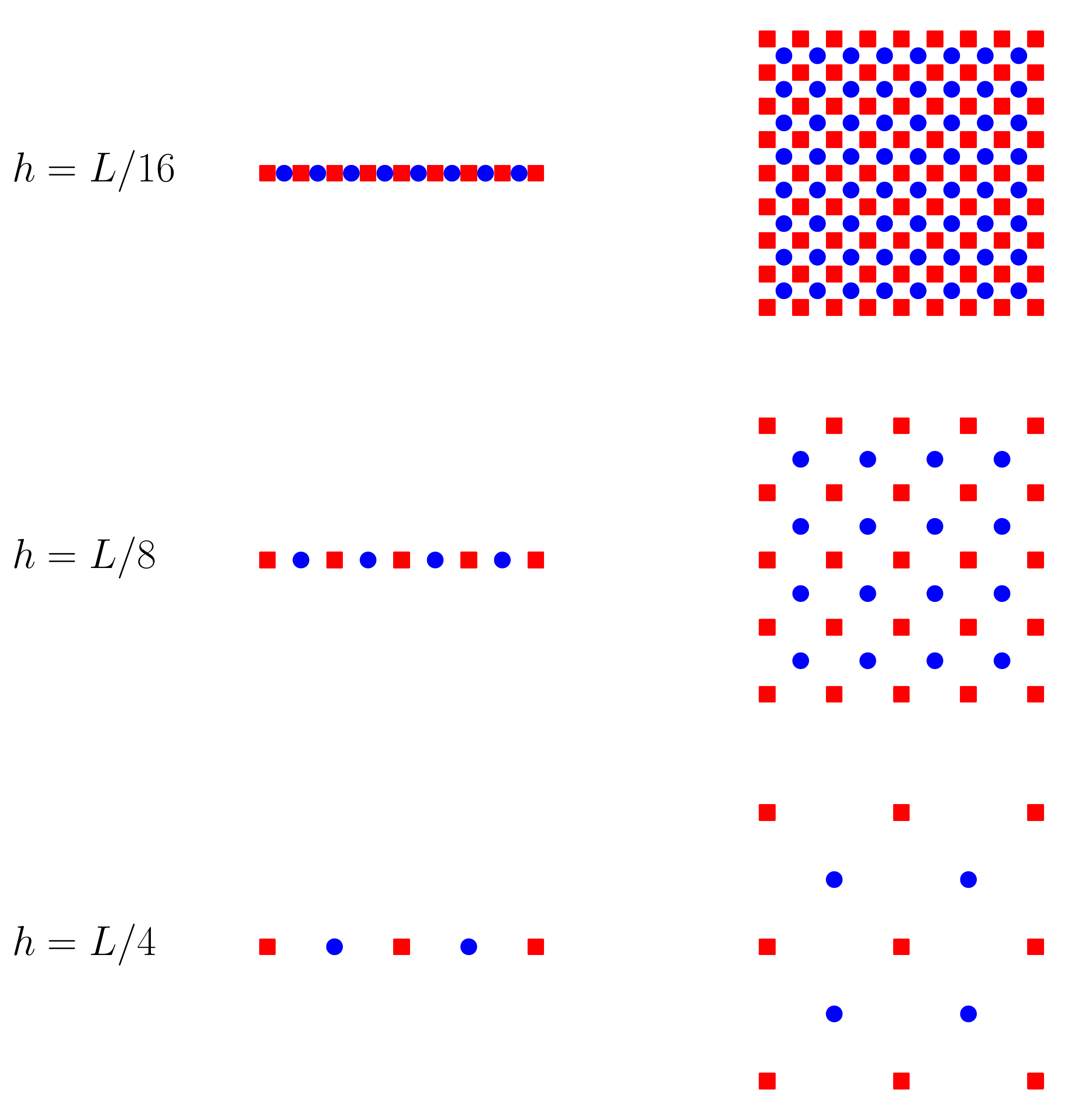}
\caption{Interlacing of red square training points and blue circle testing points for $p=1$ (left) and $p=2$ (right) dimensional datasets. Here $n$ is a power of 2, the domain length $L = 2nh$ is fixed, and $n+1$ training points are placed along each dimension to ensure the boundaries of the domain are always included for all spacings $h$.}
\label{fig:meshes}
\end{figure}

For this calculation, we will assume we have $2n + 1$ data points over an interval $L$ in $p=1$ dimensions, with uniform spacing
\begin{align}
	h =&~ \frac{L}{2n} \label{eq:hspacing}. 
\end{align}
The graph vertices are partitioned into training data $\tilde{V}$ and testing data $V$. For $\tilde{V}$, $n+1$ base points $\tilde{x}$ to fit the model will be spaced by $2h$, including the domain boundaries. For $V$, the training points are interlaced by $n$ testing points $x$ to evaluate the model, also spaced by $2h$ as per \cref{fig:meshes}. The points in $p = 1$ dimensions with uniform spacing $h$ for a given $n$ are at locations:

Training points $\tilde{x} = \{0,~2h,~4h,~6h,\dots,~(2n-2)h,~2nh\}$, 

Testing points $~x = \{h,~3h,~5h,~7h,\dots,~(2n-1)h\}$.

\noindent Note that unless specified, zero-based indexing for the vertices and points is used.

\subsection{Pointwise error of non-local first derivatives} \label{app:error_pointwisederivative}
We will now investigate the error in the non-local derivatives relative to differential derivatives, written as
\begin{align}
	\difference[1]{u(\tilde{x}_i)}{x} =&~ \derivative[1]{u(\tilde{x}_i)}{x} + \varepsilon_1(\tilde{x}_i).
\end{align}
When using the polynomial weight functions, the non-local first derivatives for $n+1$ possible Taylor series base points $\tilde{x} = \{\tilde{x}_{0},\dots,\tilde{x}_{n}\}$, have the form of
\begin{align}
	\difference[1]{u(\tilde{x}_i)}{x} = \frac{C_{\epsilon}}{n}\sum_{j \in \tilde{V}\setminus \{i\}}\frac{(u(\tilde{x}_{j})-u(\tilde{x}_{i})) (\tilde{x}_{j}-\tilde{x}_{i})}{\abs{\tilde{x}_{j}-\tilde{x}_{i}}^{2+\epsilon}}.
\end{align}
Here the polynomial form of the function means the changes in $u$ can be written as
\begin{align}
	u(\tilde{x}_{j})-u(\tilde{x}_{i}) = \sum_{l=0}^{K}\alpha_{l} (\tilde{x}_j^l-\tilde{x}_i^{l}),
\end{align}
and therefore the non-local first derivatives of polynomial functions are
\begin{align}
	\difference[1]{u(\tilde{x}_i)}{x} = \frac{C_{\epsilon}}{n}\sum_{l=0}^{K}\alpha_{l} \sum_{j \in \tilde{V}\setminus \{i\}}\frac{(\tilde{x}_j^l-\tilde{x}_i^{l})(\tilde{x}_{j}-\tilde{x}_{i})^{}}{\abs{\tilde{x}_{j}-\tilde{x}_{i}}^{2+\epsilon}}.
\end{align}
We will assume uniformly spaced data where $\tilde{x}_{j}-\tilde{x}_{i} = 2h(j-i)$. The non-local first derivatives can now be written in terms of the spacing $h$:
\begin{align}
	\difference[1]{u(\tilde{x}_i)}{x} = \frac{C_{\epsilon}}{n}\sum_{l=0}^{K}\alpha_{l} (2h)^{l-1-\epsilon} \sum_{j \in \tilde{V}\setminus \{i\}}\frac{p_l[i,j]}{\abs{j-i}^{\epsilon}}. \label{eq:nonlocalfirstorderpolyepsilon}
\end{align}
In the case of $\epsilon = 0$, the non-local first derivative takes the simplified form
\begin{align}
	\difference[1]{u(\tilde{x}_i)}{x} = \frac{1}{n}\sum_{l=0}^{K}\alpha_{l} (2h)^{l-1} \sum_{j \in \tilde{V}\setminus \{i\}}p_l[i,j]. \label{eq:nonlocalfirstorderpoly}
\end{align}
For polynomial weights and a function that depends on the coefficients $\alpha$, the error is denoted with the appropriate arguments
\begin{align}
	\varepsilon_1(\tilde{x}_i) =&~ \varepsilon_1[h,\epsilon,\alpha](\tilde{x}_i), \\
	\intertext{and given the polynomial form of the function, and the subsequent form of the differential first derivatives in \cref{eq:polyderivatives}}
	\varepsilon_1[h,\epsilon,\alpha](\tilde{x}_i)  =&~\sum_{l=0}^{K} \alpha_l\varepsilon_{1_{l}}[h,\epsilon](\tilde{x}_i) , \label{eq:errorfirstorder}
\end{align}
where the errors associated with each term in the polynomial function are
\begin{align}
	\varepsilon_{1_{l}}[h,\epsilon](\tilde{x}_i) =&~ \left[\left[\frac{C_{\epsilon}}{L}(2h)^{1-\epsilon}\sum_{j \in \tilde{V}\setminus \{i\}}\frac{p_l[i,j]}{\abs{j-i}^{\epsilon}}\right] - \left[\frac{l!}{(l-1)!}i^{l-1}\mathbbm{1}_{l\geq1} \right]\right] (2h)^{l-1}.
\end{align}
In the case of $\epsilon = 0$, the errors take the simplified form
\begin{align}
	\varepsilon_{1_{l}}[h,\epsilon = 0](\tilde{x}_i) =&~ \left[\left[\frac{2h}{L}\sum_{j \in \tilde{V}\setminus \{i\}}p_l[i,j]\right] - \left[\frac{l!}{(l-1)!}i^{l-1}\mathbbm{1}_{l\geq1} \right]\right] (2h)^{l-1}.
\end{align}

\subsection{Quadratic function modeled by linear Taylor series} \label{app:error_example}
As an example, consider a $K=2$ quadratic function modeled by a $k=1$ linear Taylor series with uniformly spaced data. The first derivatives have the form: 
\begin{align}
	\derivative[1]{u(\tilde{x}_i)}{x} =&~ \alpha_{1} + 2\alpha_{2} i2h \\
	\difference[1]{u(\tilde{x}_i)}{x} =&~ \left[\frac{C_{\epsilon}}{n}\sum_{j \in \tilde{V}\setminus \{i\}} \frac{(\tilde{x}_j - \tilde{x}_i)(\tilde{x}_j - \tilde{x}_i)}{ \abs{\tilde{x}_j - \tilde{x}_i}^{2+\epsilon}} \right]  \alpha_{1} +
	 \left[\frac{C_{\epsilon}}{n}\sum_{j \in \tilde{V}\setminus \{i\}} \frac{(\tilde{x}_j^2 - \tilde{x}_i^2)(\tilde{x}_j - \tilde{x}_i)}{\abs{\tilde{x}_j - \tilde{x}_i}^{2+\epsilon}} \right] \alpha_{2} \\
	=&~ (\alpha_1 + 2\alpha_2 i2h) \nonumber \\
	~&+~\left[\frac{C_{\epsilon}}{L}(2h)^{1-\epsilon}\sum_{j \in \tilde{V}\setminus \{i\}} \frac{1}{\abs{j-i}^{\epsilon}} - 1 \right] \alpha_{1} \nonumber \\
	~&+~\left[\frac{C_{\epsilon}}{L}(2h)^{1-\epsilon}\sum_{j \in \tilde{V}\setminus \{i\}} \frac{j+i}{\abs{j-i}^{\epsilon}} - 2i \right]2h \alpha_{2} \nonumber \\ 
	=&~ \derivative[1]{u(\tilde{x}_i)}{x} + \alpha_1\varepsilon_{1_{1}}[h,\epsilon](\tilde{x}_i) + \alpha_2\varepsilon_{1_{2}}[h,\epsilon](\tilde{x}_i).
\end{align}

\noindent The linear $\alpha_1$ error term for the first derivative is
\begin{align}
	\varepsilon_{1_{1}}[h,\epsilon](\tilde{x}_i) =&~ \frac{C_{\epsilon}}{L}(2h)^{1-\epsilon}\sum_{j \in \tilde{V}\setminus \{i\}}  \frac{1}{\abs{j-i}^{\epsilon}} - 1. \label{eq:examplerrorfirstorderalphaone}
\end{align}

\noindent In the limit that $\epsilon = 0$, the $\alpha_1$ error term is zero:
\begin{align}
	\varepsilon_{1_{1}}[h,\epsilon=0](\tilde{x}_i) =&~ \frac{1}{n}\sum_{j \in \tilde{V}\setminus \{i\}}\frac{1}{1} - 1 = 0. \\
	\intertext{In the limit that $\epsilon$ reaches the upper bound of $p=1$ before an infinite amount of data, at $i=0$ the Harmonic number limits can be used:}
	\lim_{\epsilon \to 1,~n\to \infty} \varepsilon_{1_{1}}[h,\epsilon](\tilde{x}_0) =&~ \lim_{\epsilon \to 1,~n\to \infty} (1-\epsilon)\frac{R}{L} H_{\epsilon}[n] - 1 \\    
	\to&~ (1-\epsilon)\frac{R}{L}\log{n} - 1, \nonumber \\
	\intertext{or conversely in terms of $h$}
	\lim_{\epsilon \to 1,~h \to 0} \varepsilon_{1_{1}}[h,\epsilon](\tilde{x}_0) \to&~ (\epsilon-1)\frac{R}{L}\log{h} - 1. \nonumber
\end{align}

\noindent The quadratic $\alpha_2$ error term for the first derivative is
\begin{align}
	\varepsilon_{1_{2}}[h,\epsilon](\tilde{x}_i) =&~ \left[\frac{C_{\epsilon}}{L}(2h)^{1-\epsilon}\sum_{j \in \tilde{V}\setminus \{i\}} \frac{j+i}{\abs{j-i}^{\epsilon}} - 2i \right]2h. \label{eq:examplerrorfirstorderalphatwo}
\end{align}

\noindent In the limit that $\epsilon = 0$, the $\alpha_2$ error term is spatially dependent:
\begin{align}
	\varepsilon_{1_{2}}[h,\epsilon=0](\tilde{x}_i) =&~ \left[\frac{1}{n}\sum_{j \in \tilde{V}\setminus \{i\}} (j+i) - 2i \right]2h \\
	=&~\frac{n+1}{n}\left[\frac{n}{2} - i\right]2h. \nonumber
\end{align}

\noindent The the non-local first derivative point-wise errors are
\begin{align}
	\varepsilon_{1_{0}}[h,\epsilon = 0](\tilde{x}_i) =&~ 0, \\
	\varepsilon_{1_{1}}[h,\epsilon = 0](\tilde{x}_i) =&~ 0, \\
	\varepsilon_{1_{2}}[h,\epsilon = 0](\tilde{x}_i) =&~ \frac{n+1}{n}\left[\frac{L}{2} - \tilde{x}_i\right],\\
	\intertext{and the total error is}
	\varepsilon_{1}[h,\epsilon = 0,\alpha](\tilde{x}_i) =&~ \frac{n+1}{n}\left[\frac{L}{2} - \tilde{x}_i\right]\alpha_2.
\end{align}

\begin{figure}[hpt]
\centering
\includegraphics[width=0.6\textwidth]{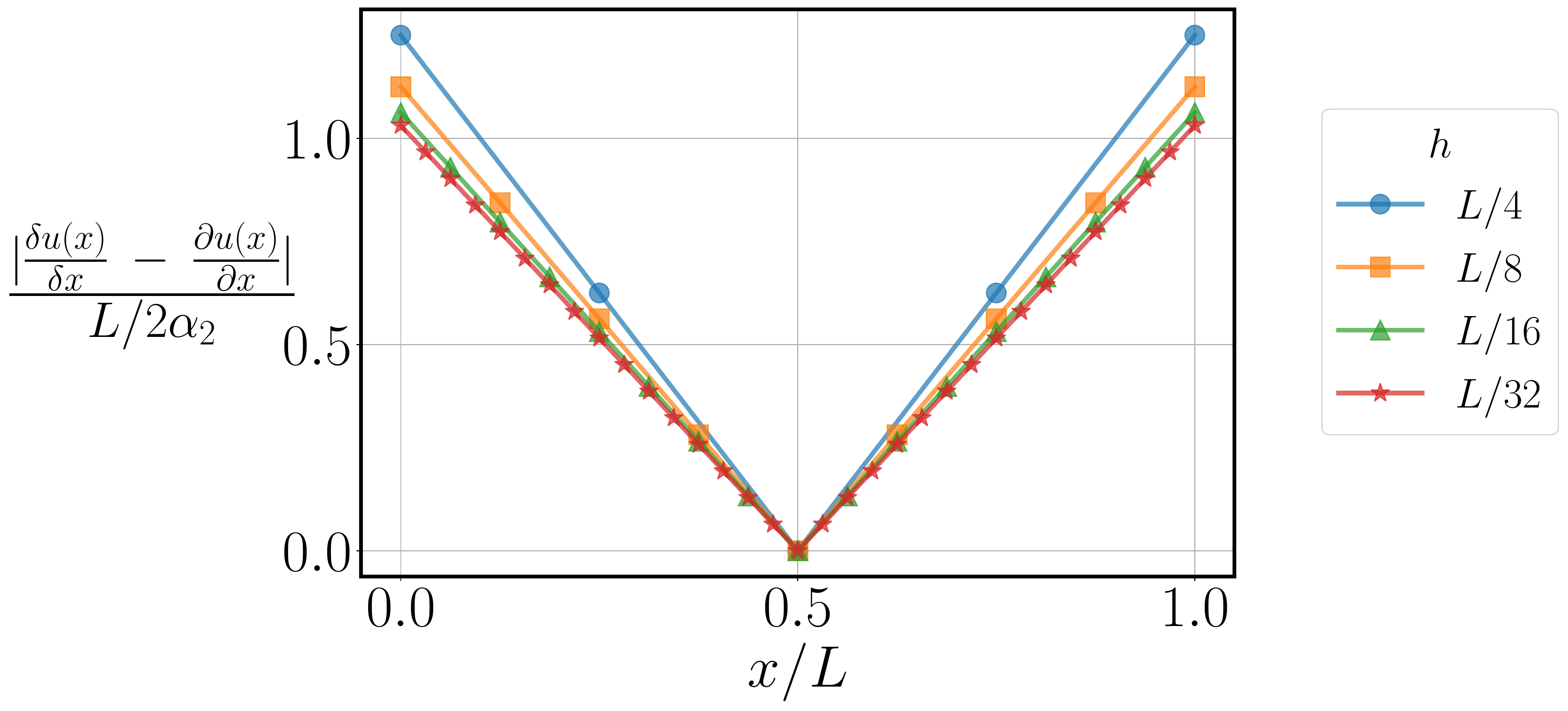}
\caption{Pointwise error of non-local first derivative for $u(x) = \alpha_0 + \alpha_1 x + \alpha_2 x^2$, and $\epsilon = 0$. $n+1$ points along an interval of length $L$ are shown, with different point spacings $h$.}
\label{fig:firstderivativeerror}
\end{figure}

As can be seen in \cref{fig:firstderivativeerror}, there is a persistent, spatially dependent error in the non-local first derivatives. For this quadratic function example, the derivative at the midpoint has minimal zero error and is exact. However as can be seen in the figure as $h$ is decreased for other points along the interval, the error reaches a finite non-zero limit, that is independent of the data spacing, the domain or interest, or the $\alpha$ coefficients. This inherent constant error due to the non-local nature of the calculus limits the pointwise error of the overall method, and future studies should aim to determine whether this affects the validity of this graph theoretic approach.

\subsection{Pointwise error of quadratic function representation} \label{app:error_pointwiseerror}
Given our knowledge of the pointwise error of the non-local derivatives, we may now calculate the pointwise error of this modified Taylor series method. The point-wise error at each point in the testing set $x_j$, for a given associated training set point $\tilde{x}_i$ to base the Taylor series is
\begin{align}
	e(x_j | \tilde{x}_i) =&~ u_k(x_j|\tilde{x}_i) - u(x_j) \label{eq:pointwiseerror}\\
	=&~ \sum_{l=0}^{k}\frac{\gamma_l(\tilde{x}_i) - 1}{l!} \uninderivative[l]{u(\tilde{x}_i)}{x}(x_j - \tilde{x}_i)^l \\
	~+&~ \sum_{l=0}^{k}\frac{\gamma_l(\tilde{x}_i)}{l!} \varepsilon_l(\tilde{x}_i)(x_j - \tilde{x}_i)^l \nonumber \\
	~-&~ \sum_{l=k+1}^{K}\frac{1}{l!} \uninderivative[l]{u(\tilde{x}_i)}{x}(x_j - \tilde{x}_i)^l \nonumber \\ 
	=&~ \sum_{l=0}^{k}\sum_{q=l}^{K}(\gamma_l(\tilde{x}_i) - 1) \binom{q}{l} \alpha_{q} \tilde{x}_i^{q-l}(x_j - \tilde{x}_i)^l \nonumber \\
	~+&~ \sum_{l=0}^{k}\frac{\gamma_l(\tilde{x}_i)}{l!} \varepsilon_l[h,\epsilon,\alpha](\tilde{x}_i)(x_j - \tilde{x}_i)^l \nonumber \\
	~-&~ \sum_{l=k+1}^{K}\sum_{q=l}^{K}\binom{q}{l} \alpha_{q} \tilde{x}_i^{q-l}(x_j - \tilde{x}_i)^l. \nonumber \\	
	=&~ C_0(x_j | \tilde{x}_i) + C_1(x_j | \tilde{x}_i)h + C_2(x_j | \tilde{x}_i)h^2 + \cdots \label{eq:pointwiseerrorconstants}\\
	\intertext{In this example of a $K=2$ quadratic function being modeled by a $k=1$ linear Taylor series, for a given $x_j$, the base point $\tilde{x}_i$ is chosen to be the closest point such that $x_j - \tilde{x}_i = h$, and $\gamma_0(\tilde{x}_i) =1$. The error constants for this example, at the point $x_j$ with base point $\tilde{x}_i$ are}
	C_0(x_j | \tilde{x}_i) =&~ 0 \label{eq:pointwiseerrorconstants_C0}\\
	C_1(x_j | \tilde{x}_i) =&~ \left[\gamma_1(\tilde{x}_i)-1\right]\alpha_1 + \left[\gamma_1(\tilde{x}_i) \frac{n+1}{n}\frac{L}{2}\right] \alpha_2 \label{eq:pointwiseerrorconstants_C1}\\
	C_2(x_j | \tilde{x}_i) =&~  \left[4\left[\gamma_1(\tilde{x}_i)-1\right]i - 2\gamma_1(\tilde{x}_i)\frac{n+1}{n}i - 1\right]\alpha_2 \label{eq:pointwiseerrorconstants_C2}\\
	C_{l>2}(x_j | \tilde{x}_i) =&~ 0 \label{eq:pointwiseerrorconstants_Cl},
\end{align}
and the constants are assumed have upper bounds
\begin{align}
	\abs{C_l(x_j | \tilde{x}_i)} \leq C_l. \label{eq:upperboundserrorconstants}
\end{align}
\noindent We see that the pointwise error in the method is due to the non-local derivatives having an error $\varepsilon_l \neq 0$ from the differential derivatives, due to the fit linear coefficient $\gamma_1(\tilde{x}_i) \neq 1$ being not identically 1 like in a differential calculus Taylor series, and due to the presence of higher order terms in the function when $K>k$. To complete this scaling analysis, we must determine how the fit linear coefficients are dependent on the training data.

\subsection{Pointwise error of linear coefficients} \label{app:error_pointwisgamma}
We must also consider the scaling of the coefficients $\gamma_l(\tilde{x}_i)$ with the spacing of the training data $\tilde{x}$ at a fixed base point $i \in \tilde{V}$. We will consider here again the example of a $k=1$ linear Taylor series with uniformly spaced data, at a fixed base point $\tilde{x}_i$. The linear equation with $p = k+1$ scalar coefficients $\gamma$ therefore takes the form
\begin{align}
	y =&~ X\gamma \\
	\intertext{where the left and right hand sides of the equation can be written as combinations of column vectors}
	\begin{bmatrix} c \end{bmatrix}  =& \begin{bmatrix} a & b \end{bmatrix}\begin{bmatrix} \gamma_0 \\ \gamma_1 \end{bmatrix}. \label{eq:firstordertaylorserieslinearmodel}
\end{align}
Depending on the choice of normalizations $N_{a,b,c}$, the components of the data-dependent quantities are
\begin{align}
	a_j =&~ u(\tilde{x}_i) \to \tilde{a}_j = N_a a_j, \\
	b_j =&~ \unidifference[1]{u(\tilde{x}_i)}{x}(\tilde{x}_j - \tilde{x}_i) \to \tilde{b}_j = N_b b_j, \\
	c_j =&~ u(\tilde{x}_j) \to \tilde{b}_j = N_c c_j,
\end{align}
and the coefficients are scaled accordingly
\begin{align}
	\gamma_0 \to&~ \tilde{\gamma}_0 = \frac{N_c}{N_a}\gamma_0, \\
	\gamma_1 \to&~ \tilde{\gamma}_1 = \frac{N_c}{N_b}\gamma_1.
\end{align}
Due to the vectors being normalized, each vector to leading order has no dependence on the data size or spacing, leading to the linear coefficients being constant to leading order, plus higher order corrections.

Here we will analyze the scaling of the coefficients based on an OLS fitting, as per \cref{eq:olslossmin}. From least squares arguments, if there are no restrictions on the base term $\gamma_0$, the resulting unnormalized coefficients from OLS fitting are
\begin{align}
	\gamma_0 =&~ \frac{1}{a}\frac{\mu_{c}\sigma_{b}^2 - \mu_{b}\sigma_{bc}^2}{\sigma_{b}^2},\\
	\gamma_1 =&~ \frac{\sigma_{bc}^2}{\sigma_{b}^2},\\
	\intertext{where the mean and variance of the variables are}
	\mu_{b} =&~ \frac{1}{n+1}\sum_{j \in \tilde{V}} b_j, \\
	\mu_{c} =&~ \frac{1}{n+1}\sum_{j \in \tilde{V}} c_j, \\
	\sigma_{b} =&~ \frac{1}{n+1}\sum_{j \in \tilde{V}} b_j^2 - \mu_{b}^2, \\
	\sigma_{bc} =&~ \frac{1}{n+1}\sum_{j \in \tilde{V}} b_jc_j - \mu_{b}\mu_{c}.
\end{align}
However for model consistency, we impose that $\gamma_0 = 1$ in the unnormalized basis for all Taylor series models, and so the linear coefficient $\gamma_1$ can be found from the equation
\begin{align}
	c - \gamma_0a =&~ \gamma_1 b \\
	\begin{bmatrix} d \end{bmatrix}  =& \gamma_1\begin{bmatrix} b \end{bmatrix}.
\end{align}
Depending on the choice of normalization $N_{a,b,c}$ in the original equation in \cref{eq:firstordertaylorserieslinearmodel}, the components of the quantities are
\begin{align}
	b_j =&~ \unidifference[1]{u(\tilde{x}_i)}{x}(\tilde{x}_j - \tilde{x}_i) \to \tilde{b}_j = N_b b_j, \\
	d_j =&~ u(\tilde{x}_j) - \gamma_0u(\tilde{x}_i) \to \tilde{d}_j = N_c d_j,
\end{align}
and the coefficients are scaled accordingly
\begin{align}
	\gamma_0 \to&~ \tilde{\gamma}_0 = \frac{N_c}{N_a}\gamma_0, \\
	\gamma_1 \to&~ \tilde{\gamma}_1 = \frac{N_c}{N_b}\gamma_1.
\end{align}
The resulting unnormalized coefficient from OLS fitting is
\begin{align}
	\gamma_1 =&~ \frac{\sum_{j \in \tilde{V}} b_jd_j}{\sum_{j \in \tilde{V}} b_j^2} \\
	=&~  \frac{1}{\unidifference[1]{u(\tilde{x}_i)}{x}} \frac{\sum_{j \in \tilde{V}} (u(\tilde{x}_j) - u(\tilde{x}_i))(\tilde{x}_j - \tilde{x}_i)}{\sum_{j \in \tilde{V}}(\tilde{x}_j - \tilde{x}_i)^2} \label{eq:gamma1_OLS} \\
	=&~  \frac{1}{\uniderivative[1]{u(\tilde{x}_i)}{x} + \varepsilon_1(\tilde{x}_i)} \frac{\sum_{j \in \tilde{V}} (u(\tilde{x}_j) - u(\tilde{x}_i))(\tilde{x}_j - \tilde{x}_i)}{\sum_{j \in \tilde{V}}(\tilde{x}_j - \tilde{x}_i)^2}.
\end{align}

\noindent To calculate this coefficient, we assume the $K$ order polynomial representation for the function, and $\epsilon = 0$ for the polynomial non-local derivative weights. We can then write the numerator and denominator of the coefficient definition when $\tilde{x}_i = 2hi$ as sums of terms that do not explicitly scale with $h$ to leading order:
\begin{align}
	\gamma_1(\tilde{x}_i) =&~ \frac{1}{\sum_{l=0}^{K}\alpha_l(2h)^{l-1}\frac{1}{n}\sum_{j \in \tilde{V}\setminus\{i\}}p_l[i,j]} \sum_{l=0}^{K}\alpha_l(2h)^{l-1}\frac{\sum_{j \in \tilde{V}}(j-i)^2p_l[i,j]}{\sum_{j \in \tilde{V}} (j-i)^2} \\
	=&~ \frac{1+\sum_{l>1}^{K} \frac{\alpha_l}{\alpha_1}A_l[i,n] (2h)^{l-1}}{1 + \sum_{l>1}^{K}\frac{\alpha_l}{\alpha_1}B_l[i,n] (2h)^{l-1}}\\
	\intertext{where the constants $A_l$ and $B_l$ are defined as}
	A_l[i,n] =&~  \frac{\sum_{j \in \tilde{V}}(j-i)^2p_l[i,j]}{\sum_{j \in \tilde{V}} (j-i)^2}, \\
	B_l[i,n] =&~  \frac{1}{n}\sum_{j \in \tilde{V}\setminus \{i\}}p_l[i,j].
\end{align}
Writing $i = \varchi_{i} n$, where $\varchi_{i} \in [0,1]$, we obtain factors that do not scale with $h$ to leading order
\begin{align}
	\tilde{A}_{l}[i,n] =&~  \frac{1}{n^{l-1}} A_l[i,n] \\
	=&~ \frac{\tilde{\varphi}_{l+1} - \varchi_{i}\tilde{\varphi}_{l} - \varchi_{i}^l \tilde{\varphi}_{1} + \varchi_{i}^{l+1}\tilde{\varphi}_{0}}{\tilde{\varphi}_{2} - 2\varchi_{i}\tilde{\varphi}_{1} + \varchi_{i}^2\tilde{\varphi}_{0}}, \\
	\tilde{B}_l[i,n] =&~  \frac{1}{n^{l-1}} B_l[i,n] \\
	=&~ \sum_{q=1}^{l-1}\varchi_{i}^{q-1}\tilde{\varphi}_{l-q} - \left[1 - 2(l-1)\frac{h}{L}\right]\varchi_{i}^{l-1}.
\end{align}
To leading order in $h$
\begin{align}
	\tilde{A}_{l}[i,n] =&~  \frac{\left[\frac{1}{l+2} + \frac{h}{L}\right] - \varchi_{i}\left[\frac{1}{l+1} + \frac{h}{L}\right] - \varchi_{i}^l\left[\frac{1}{2} + \frac{h}{L}\right] + \varchi_{i}^{l+1}\left[1 + 2\frac{h}{L}\right]}{\left[\frac{1}{3} + \frac{h}{L}\right] - 2\varchi_{i}\left[\frac{1}{2} + \frac{h}{L}\right] + \varchi_{i}^2\left[1 + 2\frac{h}{L}\right]} + O\left((h/L)^2\right), \\
	=&~ \tilde{A}_{l_0}[i]\left[1 + \tilde{A}_{l_1}[i]\frac{h}{L} + O\left((h/L)^2\right)\right] \\
	\tilde{B}_l[i,n] =&~  \sum_{q=1}^{l-1}\varchi_{i}^{q-1}\left[\frac{1}{l-q+1} + \frac{h}{L}\right] + \left[1 - 2(l-1)\frac{h}{L}\right]\varchi_{i}^{l-1} + O\left((h/L)^2\right) \\
	=&~ \tilde{B}_{l_0}[i]\left[1 + \tilde{B}_{l_1}[i]\frac{h}{L} + O\left((h/L)^2\right)\right],
	\intertext{and the leading constants are}
	\tilde{A}_{l_0}[i] =&~ \frac{3}{2(l+1)(l+2)}\frac{2(l+1) - 2(l+2)\varchi_{i} - (l+1)(l+2)\varchi_{i}^{l} + 2(l+1)(l+2)\varchi_{i}^{l+1}}{1 - 3\varchi_{i} + 3 \varchi_{i}^2}, \\
	\tilde{A}_{l_1}[i] =&~ \frac{1}{\tilde{A}_{l_0}[i]}\left[1 - \varchi_{i}^{} - \varchi_{i}^{l} + 2\varchi_{i}^{l+1}\right] - \left[1 - 2\varchi_{i}^{} + 2 \varchi_{i}^{2}\right], \\
	\tilde{B}_{l_0}[i] =&~  \sum_{q=1}^{l-1}\varchi_{i}^{q-1}\frac{1}{l-q+1} + \varchi_{i}^{l-1}, \\
	\tilde{B}_{l_1}[i] =&~  \frac{1}{\tilde{B}_{l_0}[i]}\left[\sum_{q=1}^{l-1}\varchi_{i}^{q-1} - 2(l-1)\varchi_{i}^{l-1}\right].
\end{align}
The linear coefficient can finally be written to clearly show the leading order $h$ dependence as
\begin{align}
	\gamma_1(\tilde{x}_i) =&~ \frac{1+\sum_{l>1}^{K} \frac{\alpha_l}{\alpha_1}\tilde{A}_{l_0}[i] L^{l-1}}{1 + \sum_{l>1}^{K}\frac{\alpha_l}{\alpha_1}\tilde{B}_{l_0}[i] L^{l-1}} \label{eq:gamma1_hscaling_full}\\
	=&~ \frac{1 + \sum_{l>1}^{K}\frac{\alpha_l}{\alpha_1}L^{l-1}\tilde{A}_{l_0}[i]\left[1 + \tilde{A}_{l_1}[i]\frac{h}{L} + O((h/L)^2)\right]}{1 + \sum_{l>1}^{K}\frac{\alpha_l}{\alpha_1}L^{l-1}\tilde{B}_{l_0}[i]\left[1 + \tilde{B}_{l_1}[i]\frac{h}{L} + O((h/L)^2)\right]} \label{eq:gamma1_hscaling_h}.
\end{align}
The linear coefficient therefore scales with $h$ as
\begin{align}
	\gamma_1(\tilde{x}_i) =&~ G_0[i] + G_1[i]\left(\frac{h}{L}\right) + G_2[i]\left(\frac{h}{L}\right)^2 + O\left((h/L)^3\right) \label{eq:gamma1_hscaling_explicit},
\end{align}
where the $G_l[i]$ can be found from a Taylor series expansion of the denominator in \cref{eq:gamma1_hscaling_h}, and are constant to leading order in $h$.

\noindent \textbf{Remark}: There are also cases of super-convergence. For example, for the $K = 2$ quadratic function, at the midpoint of the interval, $\varchi_{i} = 1/2$, the leading constants are:
\begin{align}
	\tilde{A}_{l_0}\left[\frac{n}{2}\right] =&~ 1, \\
	\tilde{A}_{l_1}\left[\frac{n}{2}\right] =&~ 0, \\ 
	\tilde{B}_{l_0}\left[\frac{n}{2}\right] =&~ 1, \\ 
	\tilde{B}_{l_1}\left[\frac{n}{2}\right] =&~ 0,
\end{align}
and so to leading order
\begin{equation}
	\gamma_1(L/2) = \left\{ 
		\begin{array}{cc} 
			1 + O\left((h/L)^2\right)  & \frac{\alpha_1}{\alpha_2 L} \neq -1 \\ 
			0 & \frac{\alpha_1}{\alpha_2 L} = -1. 
		\end{array} \right.
\end{equation}
Decomposing $\tilde{x}_{i} = \varchi_{i}L$ becomes most useful when computing the total error and summing over the data points. Sums of the form $n^{-1}\sum_{i=0}^{n}\varchi_{i}^l \sim c_0[l] + c_1[l]h + O(h^2)$ for any powers $l$ and $l$-dependent constants $c[l]$. Therefore averages over $i$ of the $\tilde{A}_{l_q}[i],\tilde{B}_{l_q}[i],~\textrm{and}~\tilde{G}_{q}[i]$ are constant to leading order, and the leading order scaling with respect to $h$ is explicitly shown in the forms for $\gamma_1(\tilde{x}_i)$ in \cref{eq:gamma1_hscaling_explicit}. 

\subsection{Roots of non-local derivatives}
Given that the fit linear coefficients for this Taylor series model are proportional to the reciprocal of the derivatives, we may also study the behavior of $\gamma_l(\tilde{x}_{i^{*}})$ at the roots $\tilde{x}_{i^{*}}$ of the non-local derivatives. These roots must be at integral positions $i^{*} \in \{0,1,2,\cdots,n\}$ along the interval. If the derivatives have roots within the interval for specific function coefficients $\alpha^{*}$, for a given spacing $h$, then the resulting linear problem to be solved will be potentially singular. We have seen above that for the $K=2$ quadratic case, a specific ratio of $\alpha_1/\alpha_2 L = -1$ will impose that $\gamma_1 = 0$ at a specific point. Choosing base points for the Taylor series to be at a root of the derivatives, will reduce the rank of the Taylor series basis for each operator where $\unindifference[l]{u(\tilde{x}_{i^{*}})}{x} = 0$.

\begin{figure}[hpt]
\centering
\begin{subfigure}[t]{0.49\textwidth}
	\centering
	\includegraphics[width=\textwidth]{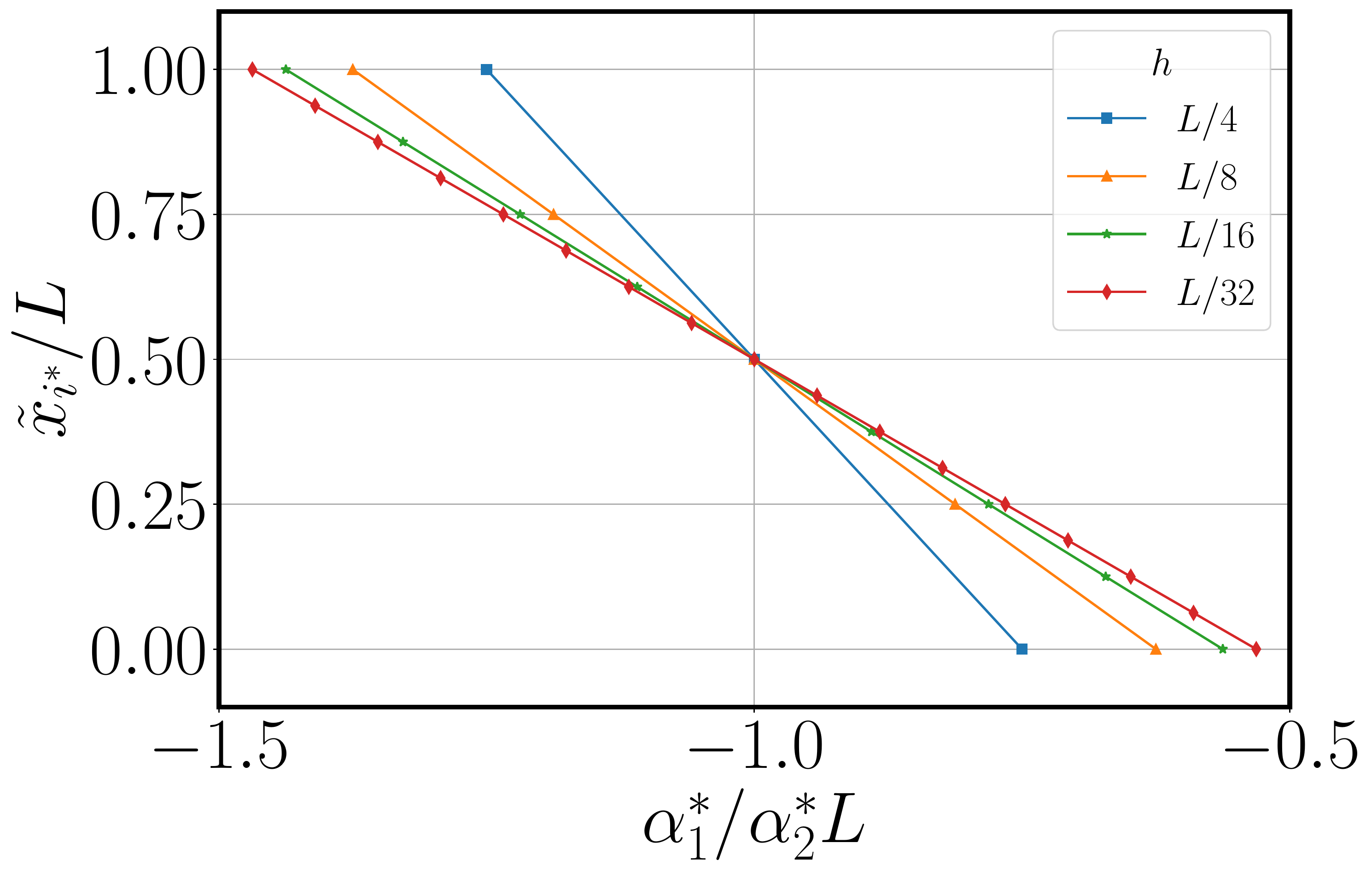}
	\subcaption{Roots of non-local first derivatives.}
	\label{fig:firstderivativebheaviour_roots}	
\end{subfigure}
\hfill
\begin{subfigure}[t]{0.49\textwidth}
	\centering
	\includegraphics[width=1.05\textwidth]{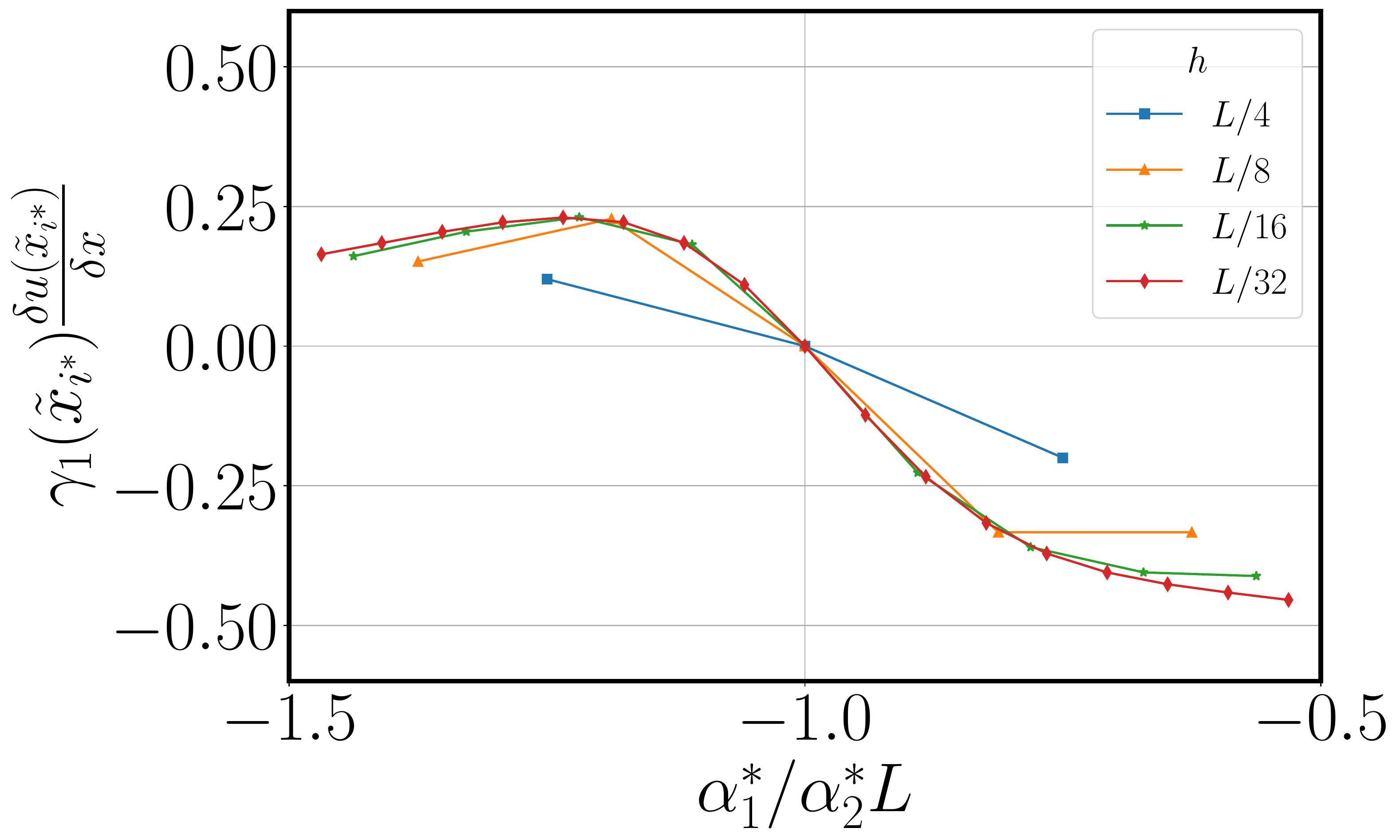}
	\subcaption{Naive product of fit linear coefficient and non-local first derivatives at non-local first derivative roots.}
	\label{fig:firstderivativebehaviour_gamma}
\end{subfigure}
\caption{Behavior of non-local first derivatives at their roots that are integral positions along the interval, for $u(x) = \alpha_0 + \alpha_1 x + \alpha_2 x^2$ and $\epsilon = 0$, and specified $\alpha$ coefficients and spacings $h$.}
\label{fig:firstderivativbehavior}
\end{figure}

For the example of a linear $k=1$ Taylor series modelling a quadratic function $u(x) = \alpha_0 + \alpha_1 x + \alpha_2 x^2$, and $\epsilon = 0$, we can derive that the roots that are integral positions along the interval have the form, for a specific ratio $\alpha_1^{*}/\alpha_2^{*}$,
\begin{align}
	\tilde{x}_{i^{*}} = \frac{n}{n-1}\left(-\frac{\alpha_1^{*}}{\alpha_2^{*}}  - \frac{n+1}{2n}L\right) ~\longleftrightarrow~ \frac{\alpha_1^{*}}{\alpha_2^{*}L} = -\left(\frac{n-1}{n^2}i^{*} + \frac{n+1}{2n}\right).
\end{align}
The behavior of the non-local first derivatives roots is shown in \cref{fig:firstderivativbehavior}.
As can be seen in \cref{fig:firstderivativebheaviour_roots}, the specific $\alpha$ that yield roots at integral positions along the interval are within the region $\alpha_1^{*}/\alpha_2^{*}L \in [-3/2,-1/2]$, which is a compressed region of coefficients compared to for differential first derivatives that have $\alpha_1^{*}/\alpha_2^{*}L \in [-2, 0]$. As shown in the plot, the coefficients associated with roots for a given $h$ are a further subinterval of this $[-3/2,-1/2]$ region, and this subinterval expands to include the limits of these bounds as $h$ decreases. Therefore for finer data, there is a greater probability that the function of interest has coefficients that correspond with a root at an integral position along the interval. This will correspond with the linear Taylor series model at this base point being unable to make predictions away from the base point. The previous point-wise error analysis for \cref{eq:pointwiseerrorconstants,eq:pointwiseerrorconstants_Cl} however still holds for all coefficient values, zero, or non-zero, and the leading order of the error scaling with respect to $h$ will also be unchanged due to the factors of $(1-\gamma_l(\tilde{x})) \sim O(1)$.

The naive product of the non-local first derivative with the fit linear coefficient $\gamma_1$ at the roots is shown in \cref{fig:firstderivativebehaviour_gamma}. This product is essentially the numerator in \cref{eq:gamma1_OLS} if you allow the division by $0$ to be canceled by the product with the first derivative. It must be specified that the chosen solution to this singular problem is $\gamma_1(\tilde{x}_{i^{*}}) = 0$, and the behavior shown in this plot is for demonstrative purposes. This plot therefore indicates that only for the specific coefficients where $\alpha_1^{*}/\alpha_2^{*}L = -1$, does the numerator of the linear coefficient also equal zero, and explicitly yield the correct $\gamma_1(\tilde{x}_{i^{*}}) = 0$. Away from this ratio, the numerator takes what is hypothesized to be a sigmoid-like form.

\begin{figure}[hpt]
\centering
\begin{subfigure}[t]{0.32\textwidth}
	\centering
	\includegraphics[width=0.89\textwidth]{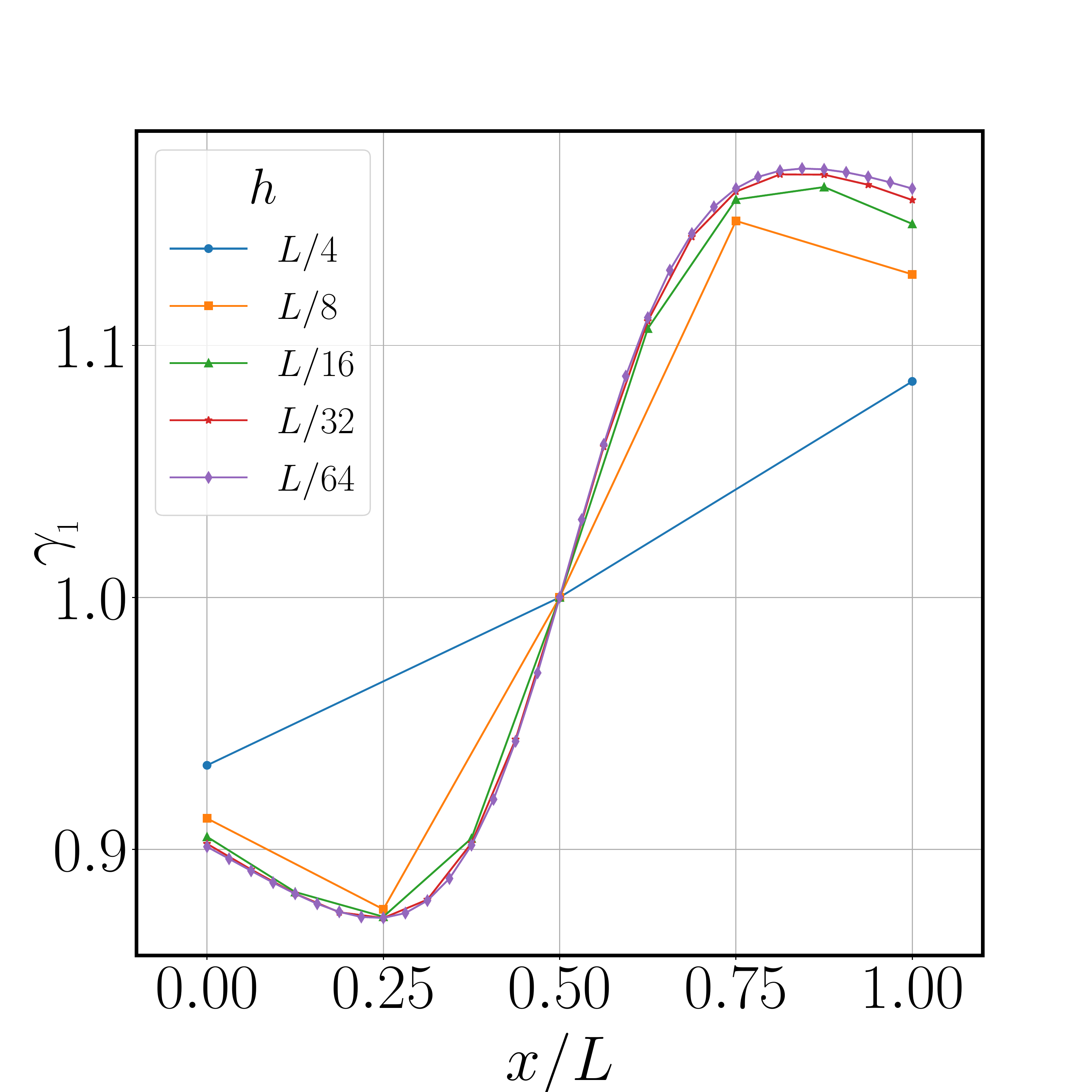}
	\subcaption{$\frac{\alpha_1}{\alpha_2 L} = -3$}
	\label{fig:gamma1_neg3}	
\end{subfigure}
\hfill
\begin{subfigure}[t]{0.32\textwidth}
	\centering
	\includegraphics[width=0.89\textwidth]{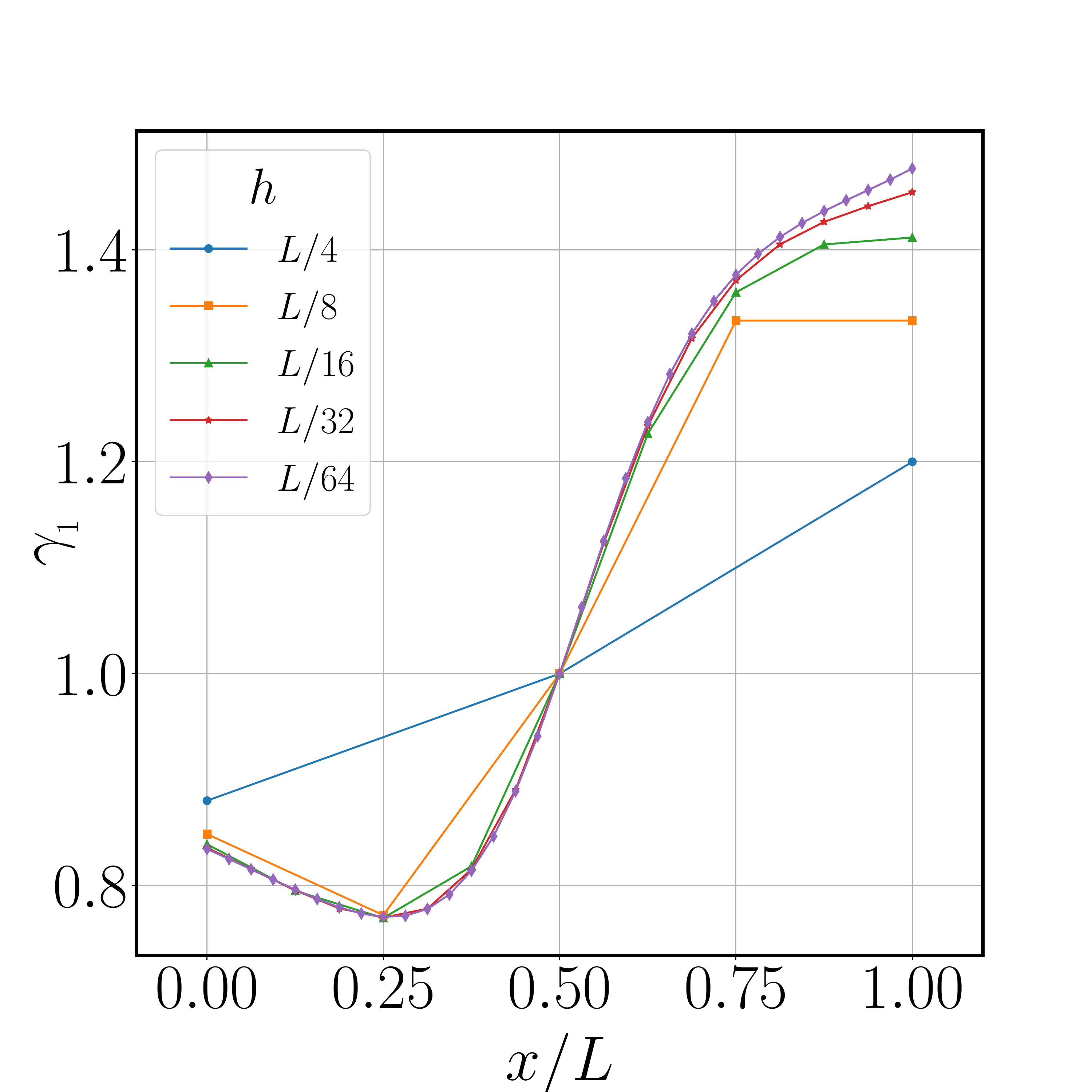}
	\subcaption{$\frac{\alpha_1}{\alpha_2 L} = -2$}
	\label{fig:gamma1_neg2}
\end{subfigure}
\hfill
\begin{subfigure}[t]{0.32\textwidth}
	\centering
	\includegraphics[width=0.89\textwidth]{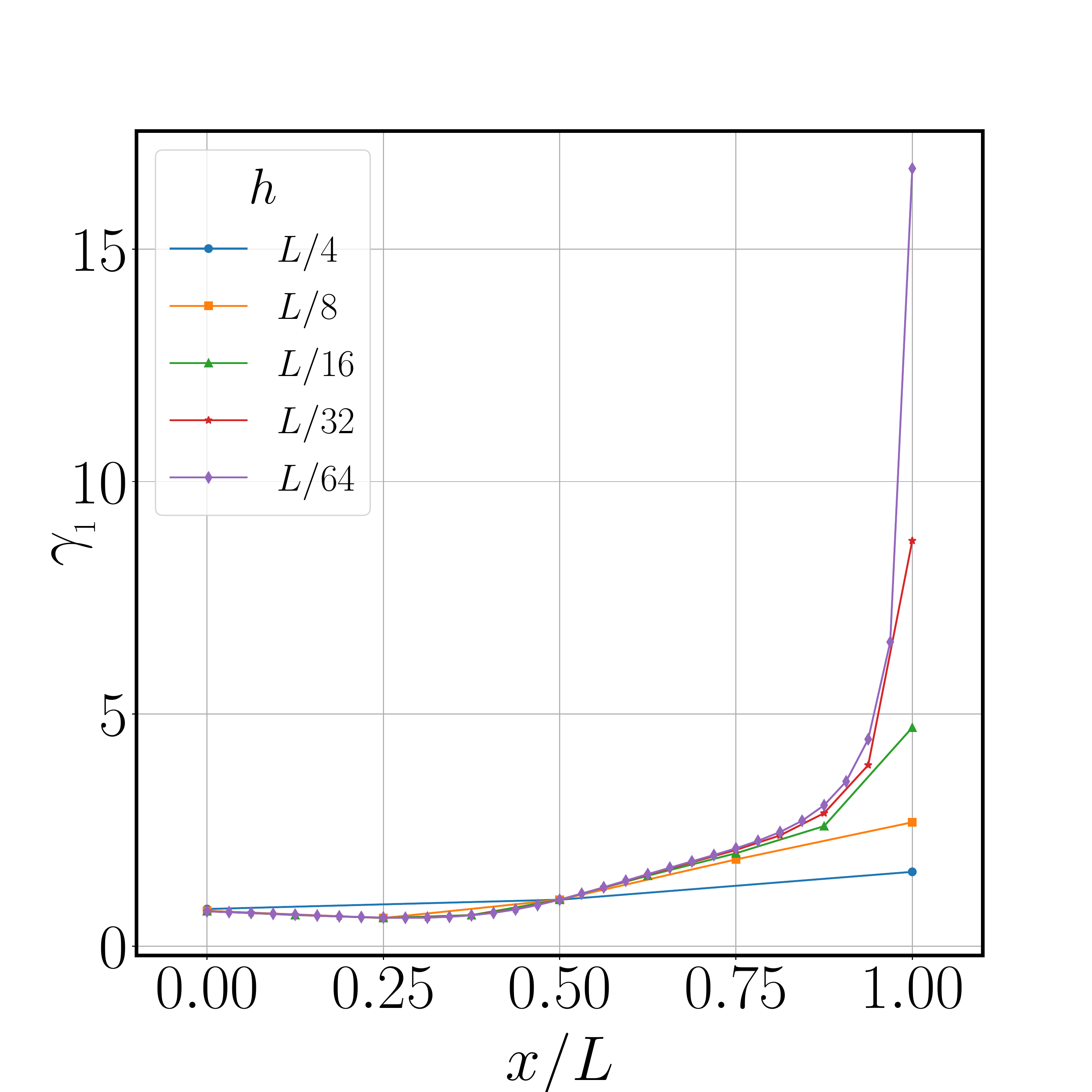}
	\subcaption{$\frac{\alpha_1}{\alpha_2 L} = -3/2$}
	\label{fig:gamma1_neg15}	
\end{subfigure}
\vfill
\begin{subfigure}[t]{0.32\textwidth}
	\centering
	\includegraphics[width=0.89\textwidth]{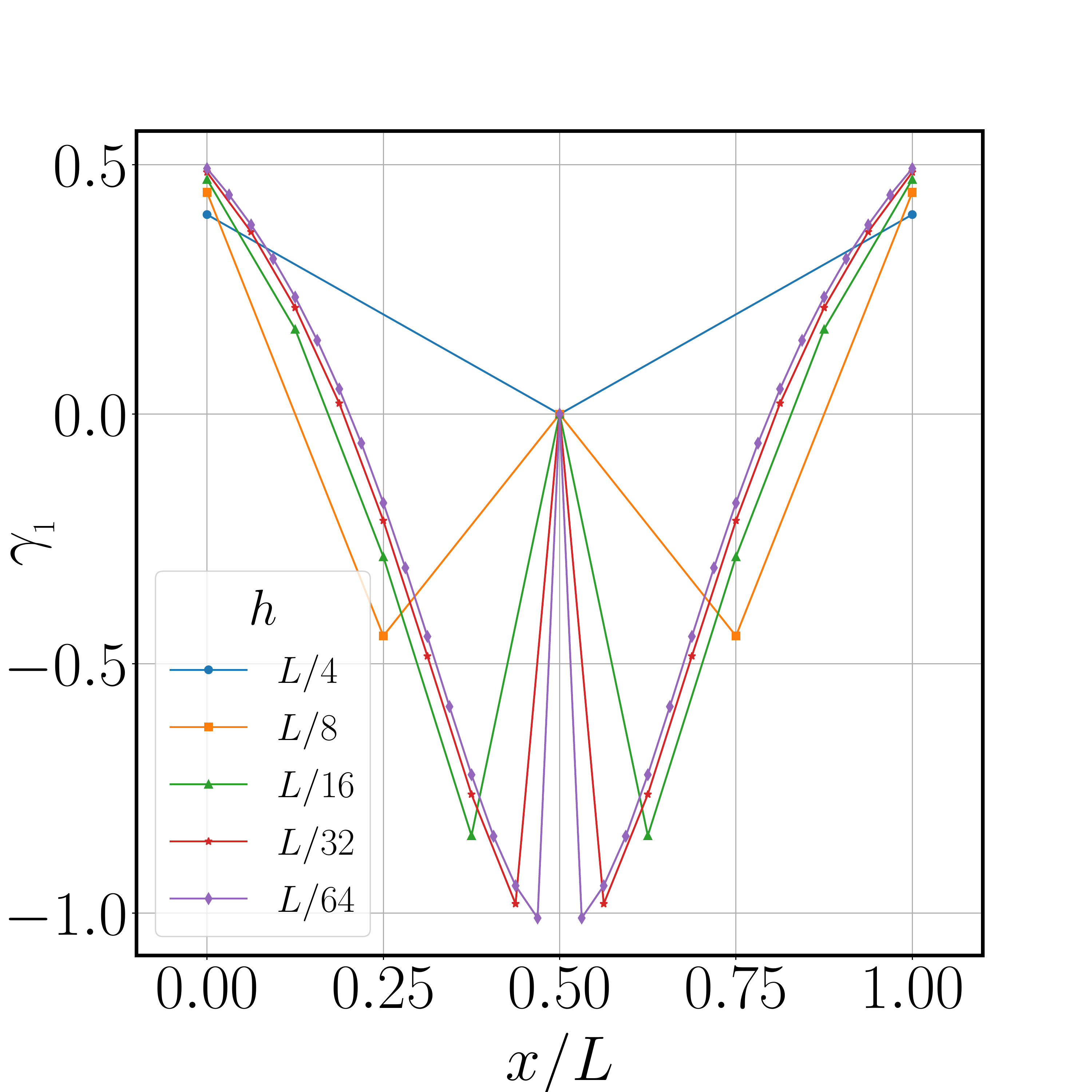}
	\subcaption{$\frac{\alpha_1}{\alpha_2 L} = -1$}
	\label{fig:gamma1_neg1}
\end{subfigure}
\hfill
\begin{subfigure}[t]{0.32\textwidth}
	\centering
	\includegraphics[width=0.89\textwidth]{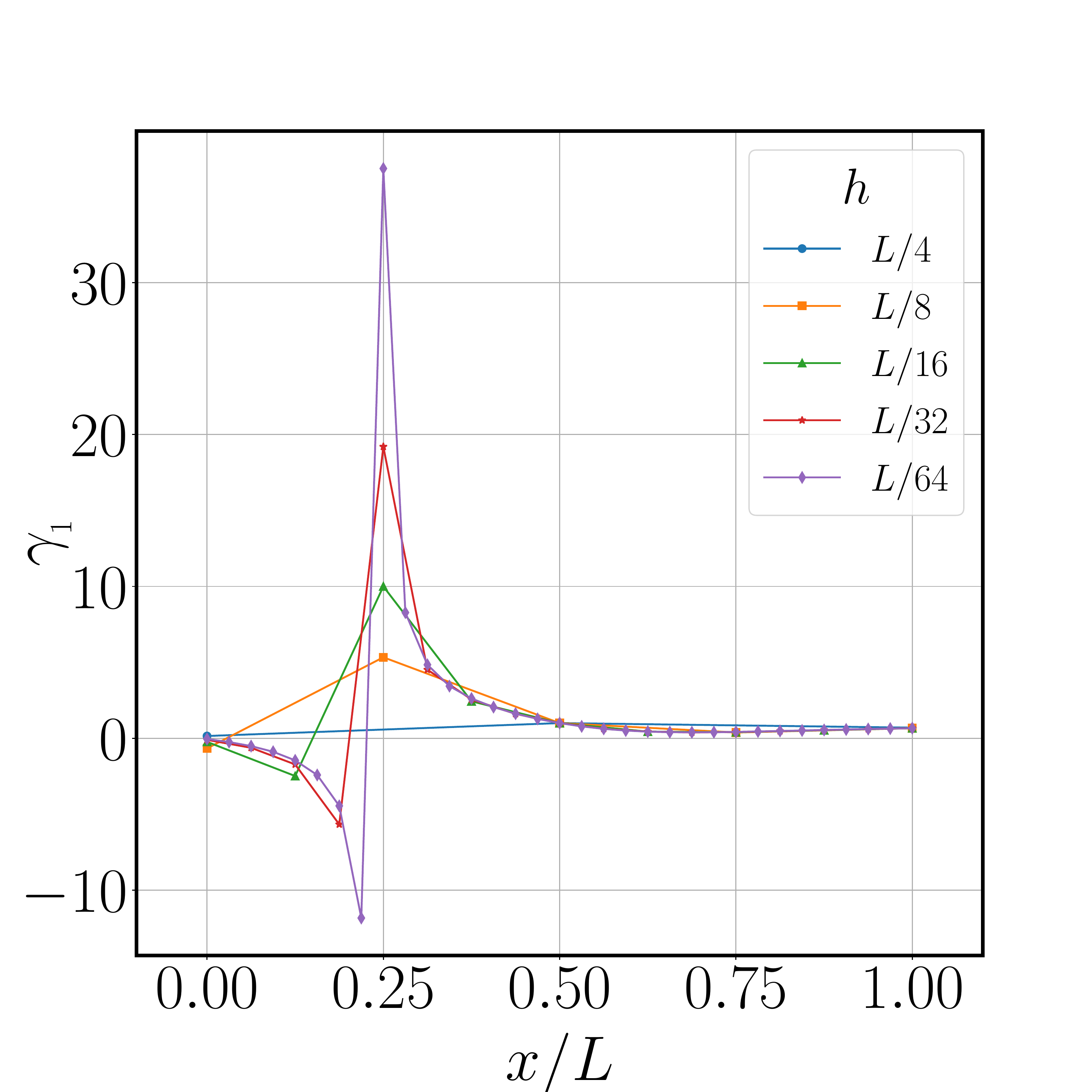}
	\subcaption{$\frac{\alpha_1}{\alpha_2 L} = -3/4$}
	\label{fig:gamma1_neg34}
\end{subfigure}
\hfill
\begin{subfigure}[t]{0.32\textwidth}
	\centering
	\includegraphics[width=0.89\textwidth]{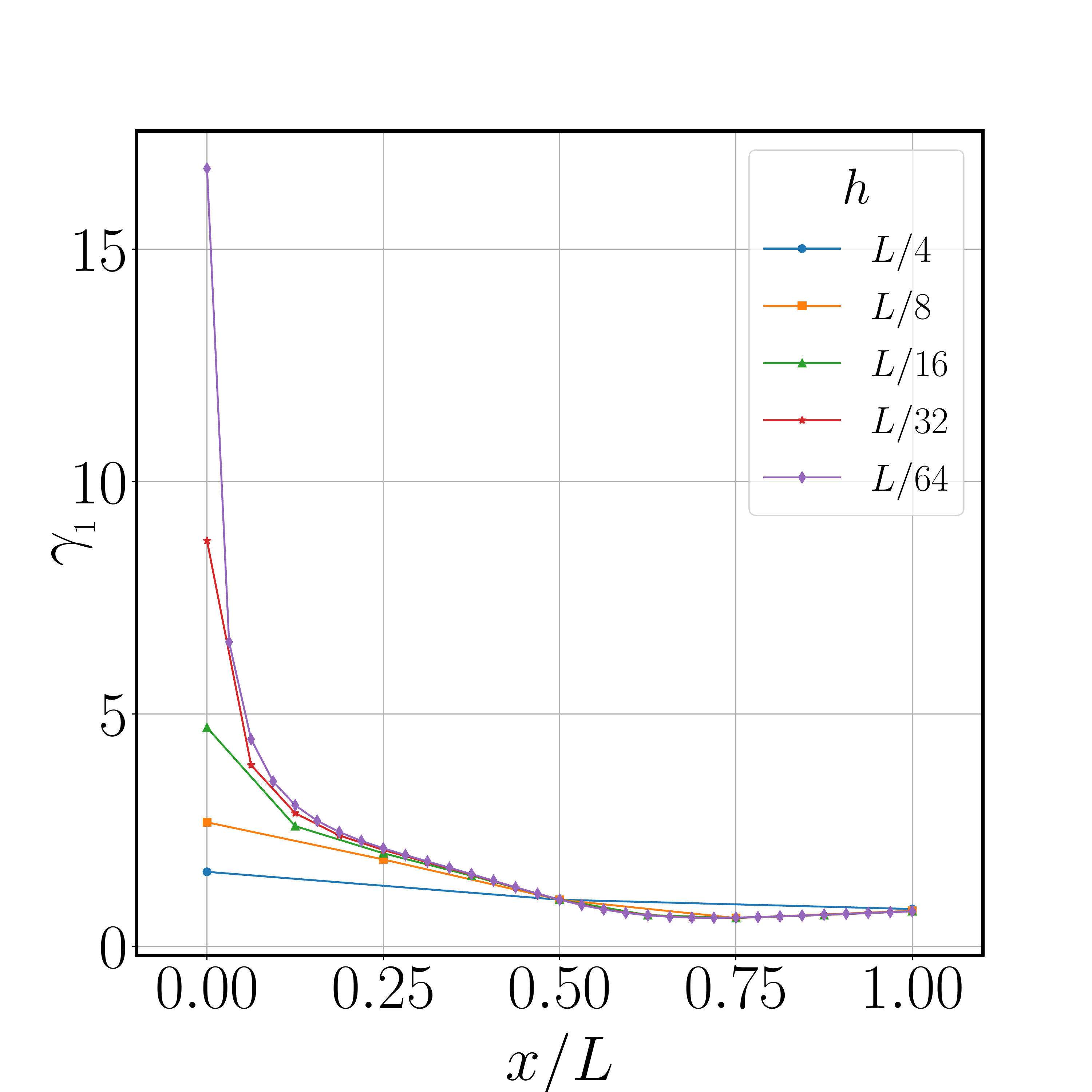}
	\subcaption{$\frac{\alpha_1}{\alpha_2 L} = -1/2$}
	\label{fig:gamma1_neg12}
\end{subfigure}
\vfill
\begin{subfigure}[t]{0.32\textwidth}
	\centering
	\includegraphics[width=0.89\textwidth]{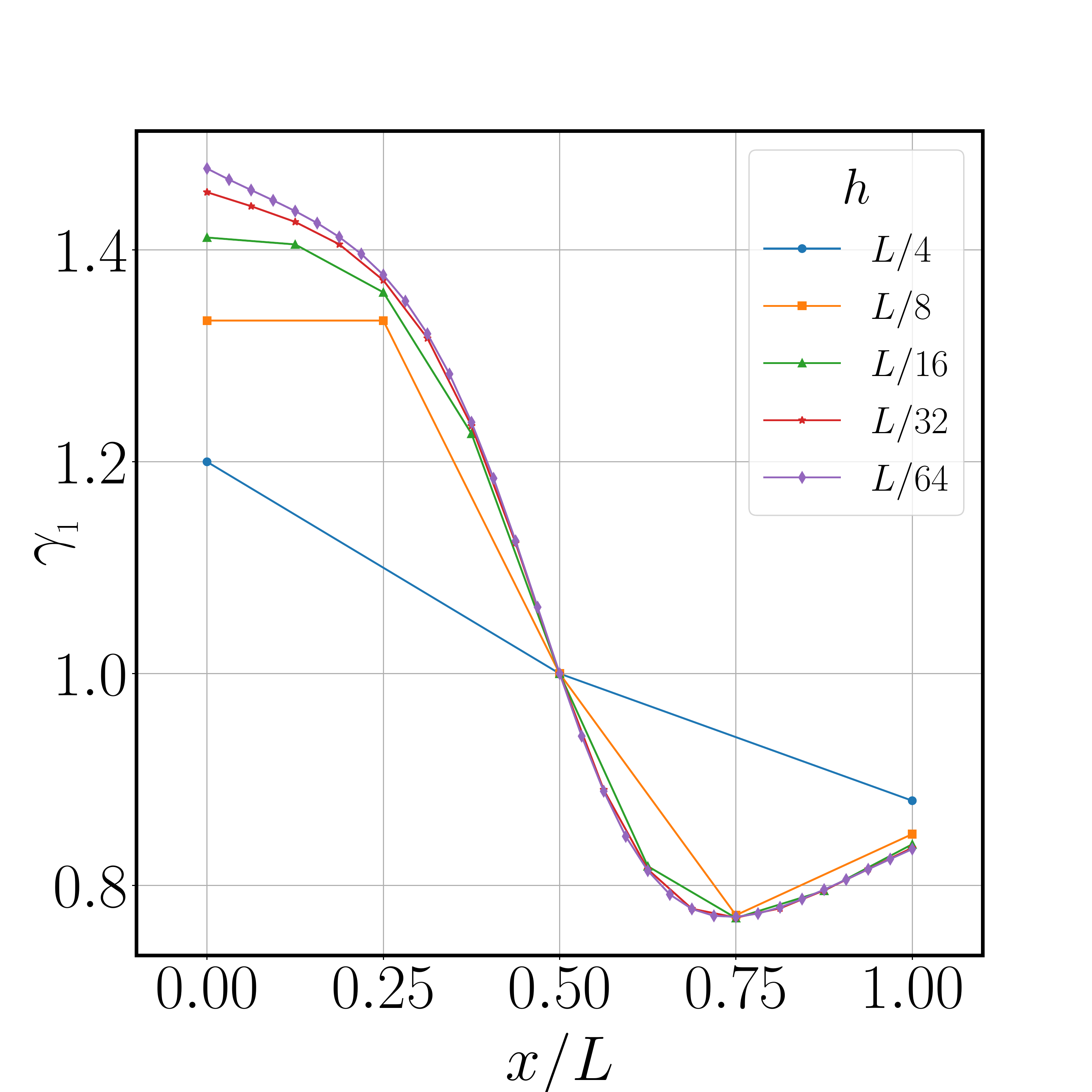}
	\subcaption{$\frac{\alpha_1}{\alpha_2 L} = 0$}
	\label{fig:gamma1_pos0}
\end{subfigure}
\hfill
\begin{subfigure}[t]{0.32\textwidth}
	\centering
	\includegraphics[width=0.89\textwidth]{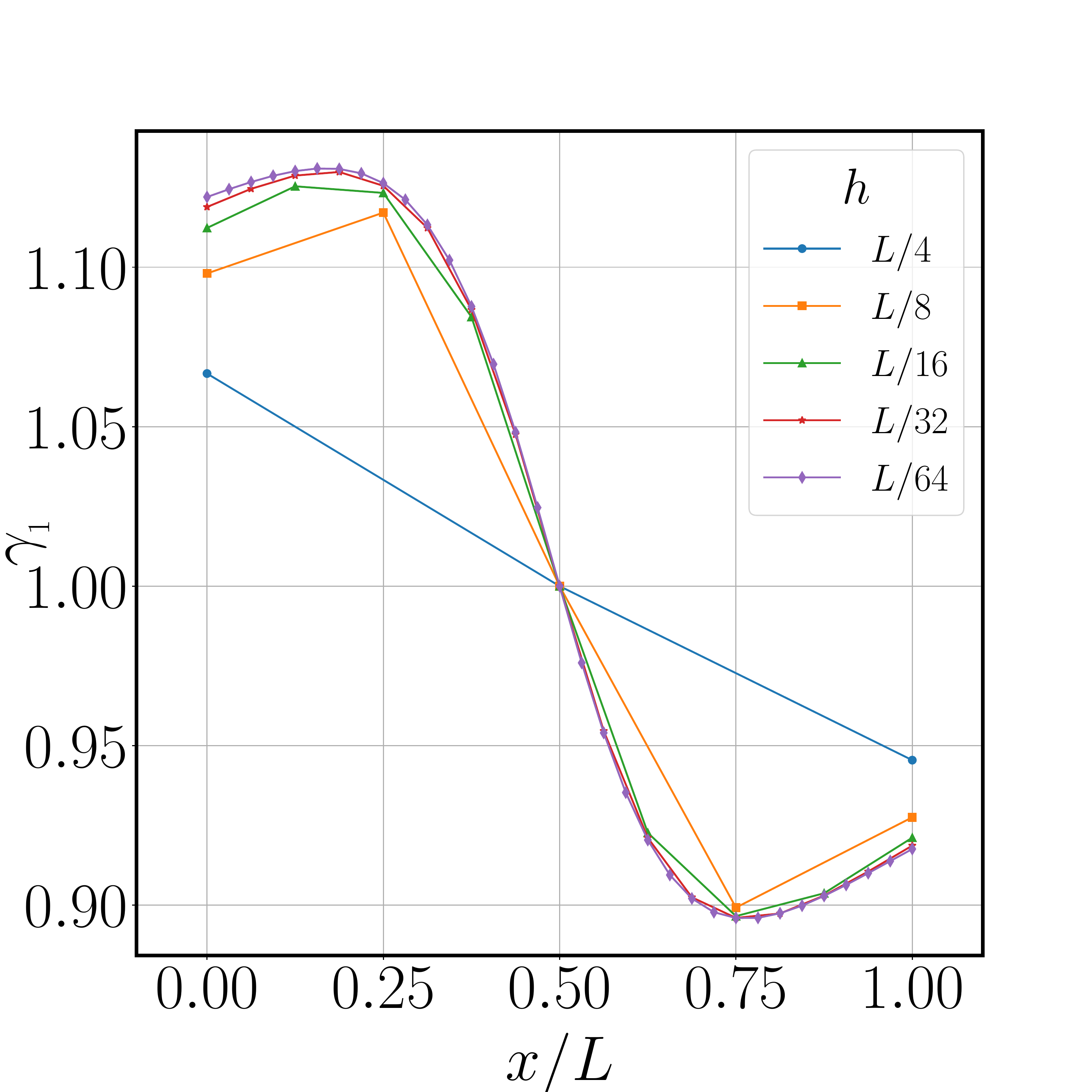}
	\subcaption{$\frac{\alpha_1}{\alpha_2 L} = 3/2$}
	\label{fig:gamma1_pos15}	
\end{subfigure}
\hfill
\begin{subfigure}[t]{0.32\textwidth}
	\centering
	\includegraphics[width=0.89\textwidth]{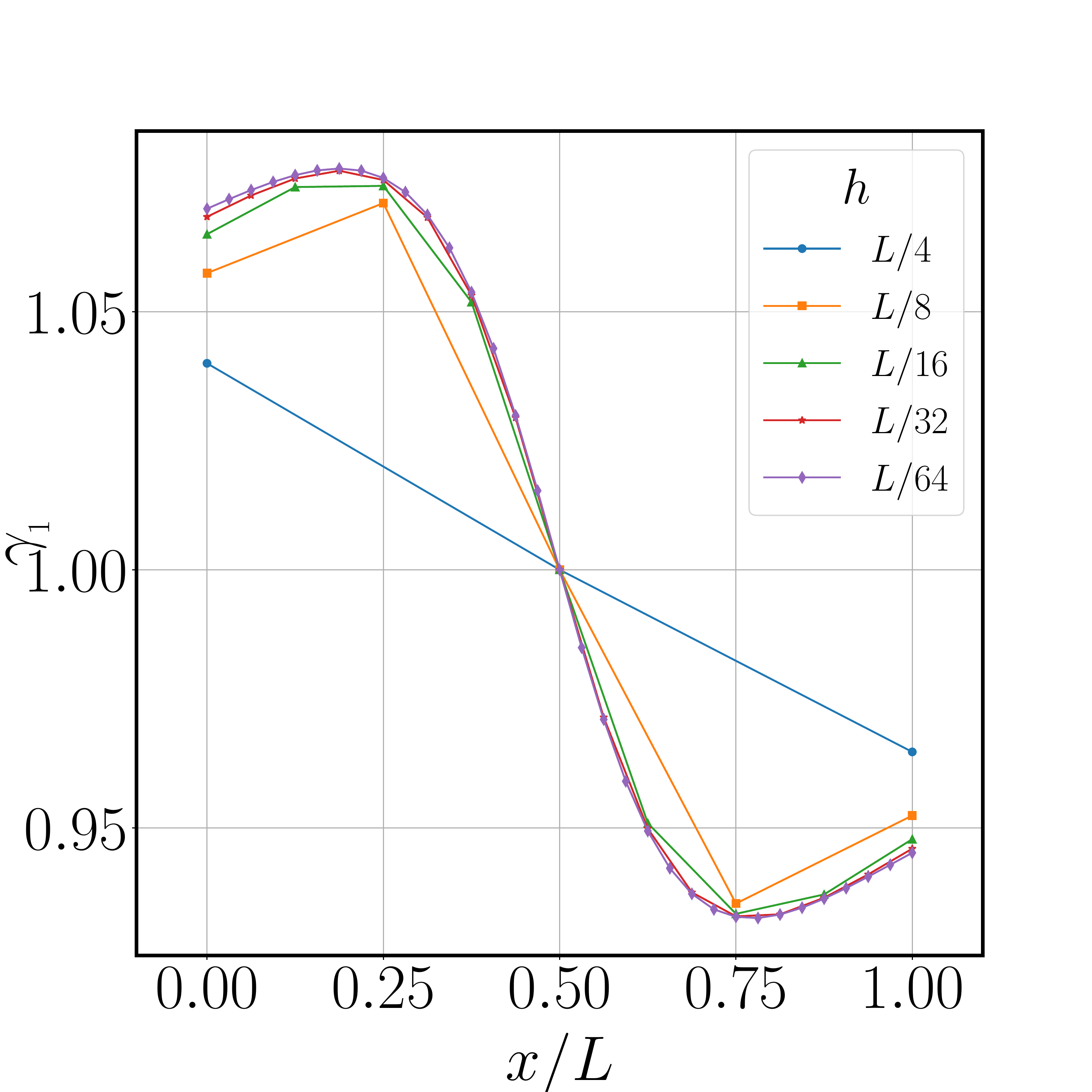}
	\subcaption{$\frac{\alpha_1}{\alpha_2 L} = 3$}
	\label{fig:gamma1_pos3}	
\end{subfigure}
\caption{Spatial behavior of linear coefficient $\gamma_1(\tilde{x})$ at the possible base points $\tilde{x}$, for $u(x) = \alpha_0 + \alpha_1 x + \alpha_2 x^2$, for various $\alpha$ coefficients and spacings $h$.}
\label{fig:gamma1}
\end{figure}

Given this understanding of the roots of the non-local first derivatives, we can plot the linear coefficient $\gamma_1(\tilde{x})$ at the possible base points $\tilde{x}$ for various $\alpha$ and spacing $h$ in \cref{fig:gamma1}. We observe that the $\gamma_1(L/2) = 1$ for all $\alpha$ outside of the region where there are possible roots of the non-local derivatives. In the special case of $\alpha_1/\alpha_2 L = -1$, $\gamma_1(L/2) = 0$, and the curve of $\gamma_1$ values is very different from other curves. For other $\alpha_1/\alpha_2 L \in [-3/2, -1/2]$, we also observe an abrupt divergence in the $\gamma_1$ values at points near the associated roots.

\newpage
\subsection{Total error of linear model for quadratic function} \label{app:error_totalerror}
The total error when considering the finite set of testing points $x_j$, given the training points $\tilde{x}_i$ can be written as being proportional to the average in the $l$-norm of pointwise errors:
\begin{align}
	\norm{e}_l^l \sim&~ \frac{1}{n}\sum_{j \in V}\abs{e(x_j)}^l, \label{eq:totalloss}
\end{align}
and given the form of the point-wise error
\begin{align}
	\norm{e}_l^l \sim&~ \frac{2h}{L}\sum_{j \in V}\abs{C_0(x_j|\tilde{x}_i) + C_1(x_j|\tilde{x}_i)h + C_2(x_j|\tilde{x}_i) h^2 + \cdots}^l \\
	\leq&~ \frac{2h}{L}n\abs{C_0 + C_1h + C_2h^2 + \cdots}^l \\
	\sim&~ \abs{C_0 + C_1h + C_2h^2 + \cdots}^{l},
\end{align}
and therefore the total error scales as 
\begin{align}
	\norm{e}_l \sim C_0 + C_1h + C_2h^2
\end{align}

\textbf{Remark}: To derive the form of the average error in \cref{eq:totalloss}, we note two common approaches to estimate the total $l$-norm error $\norm{e}_l$ across the domain of interest. We assume there is a closed form for the point-wise error $e(x|\tilde{x}) \sim C_q(\tilde{x})(x-\tilde{x})^q \sim C_q(x|\tilde{x})h^q$ that is $q$ order, and the data mesh has spacing $2h$ that splits the domain into $n$ intervals.

The first method involves the total interpolation error of the approximating function along the entire domain:
\begin{align}
	\norm{e}_l^l =&~ \frac{1}{L} \int_{0}^L dx \abs{e(x)}^l \\
	=&~ \frac{1}{L} \sum_{i=0}^{n-1}\int_{\tilde{x}_i}^{\tilde{x}_{i+1}} dx~ \abs{e(x|\tilde{x}_i)}^l \nonumber\\
	=&~ \frac{1}{L} \sum_{i=0}^{n-1} \int_{\tilde{x}_i}^{\tilde{x}_{i+1}} dx~ \abs{C_q(\tilde{x}_i)}^l (x-\tilde{x}_i)^{lq} \nonumber\\
	=&~ \frac{1}{L} \sum_{i=0}^{n-1} \frac{\abs{C_q(\tilde{x}_i)}^l}{lq + 1} (2h)^{lq + 1} \nonumber\\
	\leq&~ \frac{1}{L} \sum_{i=0}^{n-1} \frac{C_q^l}{lq + 1} (2h)^{lq + 1} \nonumber\\
	\leq&~ \frac{n}{L} \frac{C_q^l}{lq + 1} (2h)^{lq + 1} \nonumber \\
	\leq&~ 2^{lq}\frac{C_q^l}{lq + 1} h^{lq} \nonumber \\
	\intertext{ and therefore the total interpolation error is}
	\norm{e}_l \leq&~ \left[2^{q}\frac{C_q}{(lq + 1)^{1/l}}\right] h^{q} \sim O(h^{q}).
\end{align}

We now may determine the scaling for differential and non-local calculus based Taylor series models, where a linear model is approximating a quadratic function, and the $l = 2$ norm for the error is used.

\noindent In the case of differential calculus, a linear Taylor series model has $q[K=2,k=1,p=1]=2$ leading order error, and leading order error constant 
\begin{align}
	C_{2}(\tilde{x}) =&~ \frac{1}{2!}\uniderivative[2]{u(\tilde{x})}{x} \leq C_2,
	\intertext{and for a quadratic function}
	C_2 =&~ \alpha_2.
\end{align} 
\noindent For the modified Taylor series model, where the pointwise error depends on the error in the non-local derivatives and the fit linear coefficients, the error constants take the form of \cref{eq:pointwiseerrorconstants}, there is $q[K=2,k=1,p=1]=1$ leading order error, and the leading order error constant is
\begin{align}
	C_1(\tilde{x}) =&~ \left[\gamma_1(\tilde{x}_i)-1\right]\alpha_1 + \left[\gamma_1(\tilde{x}_i) \frac{n+1}{n}\frac{L}{2}\right] \alpha_2 \\
	\intertext{and given that $\gamma_1(\tilde{x}_i) = G_0[i] + G_1[i]\frac{h}{L} + O\left((h/L)^2\right)$}
	\leq&~ \left[G_0[i]-1 \right]\alpha_1 + \left[G_0[i]\frac{n+1}{n}\frac{L}{2}\right] \alpha_2 \\
	\leq &~ C_1 \\
	=&~ \left[G_{0_{\textrm{max}}} - 1\right]\alpha_1 + \left[G_{0_{\textrm{max}}}\frac{n+1}{n}\frac{L}{2}\right] \alpha_2.
\end{align}

\noindent The total error using this first method for the differential Taylor series therefore scales as
\begin{align}
	\norm{e}_2 \leq~ \left[4\frac{C_{2}}{5^{1/2}}\right] h^{2} \sim O(h^{2}) ~\longrightarrow~ Q[K=2,k=1,p=1]= 2. \label{eq:totalerrordifferential1}
\end{align}

\noindent The total error using this first method for the non-local Taylor series scales as
\begin{align}
	\norm{e}_2 \leq~ \left[2\frac{C_{1}}{3^{1/2}}\right] h \sim O(h) ~\longrightarrow~ Q[K=2,k=1,p=1]= 1. \label{eq:totalerrornonlocal1}
\end{align}
We therefore can see that the non-local Taylor series model scales an order less than the differential Taylor series model
\begin{align}
	Q[K=2,k=1,p=1]_{\textrm{non-local}} < Q[K=2,k=1,p=1]_{\textrm{differential}}.
\end{align}

\newpage
The second method involves approximating the total interpolation error of the model along the intervals of the entire domain by the error at the point $x_j \in [\tilde{x}_i,\tilde{x}_{i+1}]$ in the interior of the interval. The resulting form is what is used in \cref{eq:totalloss}:
\begin{align}
	\norm{e}_l^l =&~ \frac{1}{L} \int_{0}^L dx \abs{e_{j}(x)}^l \\
	=&~ \frac{1}{L} \sum_{i=0}^{n-1}\int_{\tilde{x}_{i}}^{\tilde{x}_{i+1}} dx~ \abs{e(x|\tilde{x}_i)}^l \nonumber\\
	\approx&~ \frac{1}{L} \sum_{i=0}^{n-1} 2h \abs{e(x_j |\tilde{x}_i)}^l \nonumber\\
	=&~ \frac{1}{n} \sum_{i=0}^{n-1} \abs{e(x_j |\tilde{x}_i)}^l \nonumber \\	
	=&~ \frac{1}{n} \sum_{i=0}^{n-1} \abs{C_q(\tilde{x}_i)}^l h^{lq} \nonumber \\
	\leq&~ \frac{1}{n} \sum_{i=0}^{n-1} C_{q}^l h^{lq} \nonumber \\
	\leq&~ C_{q}^l h^{lq} \nonumber \\
	\intertext{ and therefore the approximated total interpolation error is}
	\norm{e}_l \leq&~ \left[C_{q}\right] h^{q} \sim O(h^{q}).
\end{align}
\noindent The total error using this second method for the $q[K=2,k=1,p=1]= 2$ differential Taylor series therefore scales as
\begin{align}
	\norm{e}_2 \leq~ C_2 h^{2} \sim O(h^{2}) ~\longrightarrow~ Q[K=2,k=1,p=1]= 2. \label{eq:totalerrordifferential2}
\end{align}
\noindent The total error using this second method for the $q[K=2,k=1,p=1]= 1$ non-local Taylor series scales as
\begin{align}
	\norm{e}_2 \leq~C_1 h \sim O(h) ~\longrightarrow~ Q[K=2,k=1,p=1]= 1. \label{eq:totalerrornonlocal2}
\end{align}
\noindent Therefore for a linear Taylor series model of a quadratic function, the non-local Taylor series is zero order accurate, and is an order less accurate with spacing $h$ than the differential Taylor series.

The total error and pointwise errors at $\tilde{x} = \{0,~L/4,~L/2,~3L/4,~L\}$ are shown in \cref{fig:totalerror}, for $\alpha_1/\alpha_2 L = 1$, and $\epsilon = 0$. Here, the $O(h)$ scaling for the total error, and majority of pointwise errors is shown, and there is also the region of super-convergence at the midpoint $L/2$ of the interval. This point is super-convergent with the error scaling as $O(h^2)$ due to having zero error between the non-local and differential first derivatives, and $\gamma_1(L/2) = 1$, therefore making the model the exact Taylor series at this point.
\afterpage{
\begin{figure}[hpt]
\centering
\includegraphics[width=0.8\textwidth]{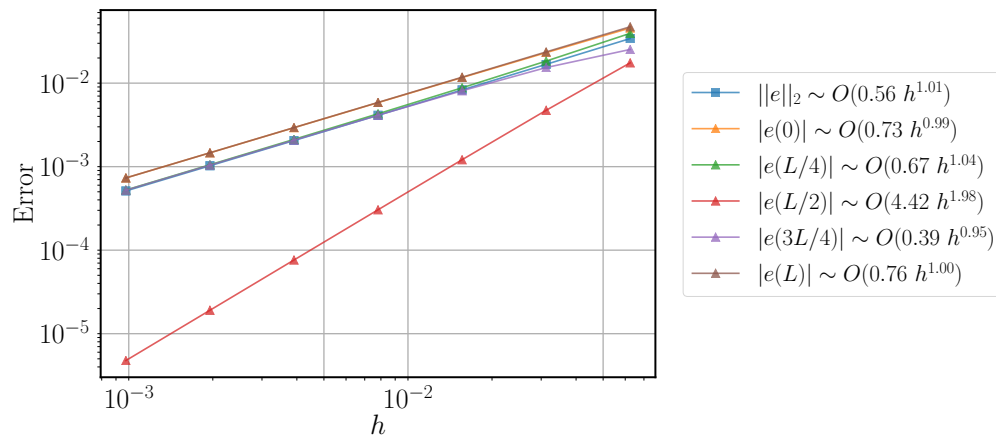}
\caption{Global and pointwise error scaling for linear Taylor series model for $u(x) = \alpha_0 + \alpha_1 x + \alpha_2 x^2$, $\alpha_1/\alpha_2 L = 1$, and $\epsilon = 0$. Leading order scaling with $h$ fits are shown in the legend.}
\label{fig:totalerror}
\end{figure}
}


\newpage
\clearpage
\section{Error analysis of modified Taylor series for polynomial functions in higher dimensions}\label{app:error_p}
\setcounter{equation}{0}
\setcounter{figure}{0}
\subsection{Modified Taylor series model} \label{app:error_model_p}
An error analysis will be conducted for a function $u(x)$ of a $p$ dimensional variable $x = \{x^{\mu}\}$, represented by a modified $k$ order Taylor series $u_k(x)$ with non-local derivatives. 
These derivatives are assumed to use polynomial weights, with super-quadratic scaling $\epsilon$.
\noindent For $p$ dimensions, the polynomial weights take the following form for $0 \leq \epsilon < p$:
\begin{align}
	w(x_i,x_j) = C_{\epsilon}\frac{1}{\left[\sum_{\mu}\abs{x^{\mu}_j-x^{\mu}_i}^2\right]^{\frac{2+\epsilon}{2}}},
\end{align}
where for the continuous weight constraints, $C_{\epsilon} = (p-\epsilon)R^{\epsilon}$ and $\lim_{\epsilon \to 0} C_{\epsilon} \to p$.

The Taylor series functional representations of the $K$ order function $u(x)$ and $k\leq K$ order model $u_k(x)$, based at a point $\tilde{x}$, are
\begin{align}
	u(x|\tilde{x}) =&~ u({}\tilde{x}) + \sum_{\mu}\derivative[1]{u(\tilde{x})}{x^{\mu}}(x^{\mu}-\tilde{x}^{\mu}) + \sum_{\mu,\nu}\frac{1}{2!}\derivative[2]{u(\tilde{x})}{x^{\mu},x^{\nu}}(x^{\mu}-\tilde{x}^{\mu})(x^{\nu}-\tilde{x}^{\nu}) + \cdots \label{eq:differentialtaylor_p}\\
	&~~~~~~~+~ \sum_{\mu_{0},\cdots,\mu_{K-1}}\frac{1}{K!}\nderivative[K]{u(\tilde{x})}{x^{\mu_{0}}}{x^{\mu_{K-1}}}(x^{\mu_{0}}-\tilde{x}^{\mu_{0}})\cdots(x^{\mu_{K-1}}-\tilde{x}^{\mu_{K-1}}), \nonumber\\
	u_k(x|\tilde{x}) =& u(\tilde{x}) + \sum_{\mu}\gamma_{1}^{\mu}(\tilde{x})\difference[1]{u(\tilde{x})}{x^{\mu}}(x^{\mu}-\tilde{x}^{\mu}) + \sum_{\mu,\nu}\frac{\gamma_{2}^{\mu\nu}(\tilde{x})}{2!}\difference[2]{u(\tilde{x})}{x^{\mu},x^{\nu}}(x^{\mu}-\tilde{x}^{\mu})(x^{\nu}-\tilde{x}^{\nu}) + \cdots \label{eq:modeltaylor_p} \\
	&~~~~~~+~ \sum_{\mu_{0},\cdots,\mu_{k-1}}\frac{\gamma_{k}^{\mu_{0}\cdots\mu_{k-1}}(\tilde{x})}{k!}\ndifference[k]{u(\tilde{x})}{x^{\mu_{0}}}{x^{\mu_{k-1}}}(x^{\mu_{0}}-\tilde{x}^{\mu_{0}})\cdots(x^{\mu_{k-1}}-\tilde{x}^{\mu_{k-1}}) \nonumber, 	
\end{align}
where the coefficients $\gamma_{k}(\tilde{x})$ are the linear coefficients fit to the non-local derivative model based around $\tilde{x}$. 

\subsection{Functional representation} \label{app:error_functional_p}
Assume the function to be represented has a polynomial form with coefficients $\alpha_{l}$ up to $K$ order, with independent behavior along each dimension:
\begin{align}
	u(x) =&~ \sum_{l=0}^{K}\sum_{\mu}\alpha_{l\mu} {x^{\mu}}^{l}.
\end{align}
Due to the assumed independence of each dimension, the higher derivatives of this function are only non-zero for derivatives along the same dimension:
\begin{align}
	\nderivative[{l+1}]{u(x)}{{x^{\mu}\partial x^{\mu_{0}}}}{x^{\mu_{l-1}}} =&~ \sum_{q=l+1}^{K} \frac{q!}{(q-l-1)!}\alpha_{q\mu}{x^{\mu}}^{q-l-1}\delta^{\mu\mu_0}\cdots\delta^{\mu\mu_{l-1}} \label{eq:polyderivatives_p} \\
	\intertext{and so the Taylor series representation in \cref{eq:differentialtaylor_p} has the form}
	u(x|\tilde{x}) =&~ \sum_{l=0}^{K}\sum_{\mu} \sum_{q=l}^{K} \binom{q}{l}\alpha_{q\mu}{\tilde{x}^{\mu^{q-l}}}({x^{\mu}-\tilde{x}^{\mu}})^l. \label{eq:polytaylorseries_p}
\end{align}

\subsection{Data mesh} \label{app:error_mesh_p}
For this calculation, we will assume we have $2n + 1$ data points over an interval $L$ along each of the $p$ dimensions, with uniform spacing
\begin{align}
	h =&~ \frac{L}{2n} \label{eq:hspacing_p}, 
\end{align}
giving $N = (n+1)^p + n^p$ total points. Each point will be assigned the index vector $j$ such that
\begin{align}
	x_j = 2h(j^0,\cdots,j^{p-1}).
\end{align}
The graph vertices are partitioned into a training data $\tilde{V}$ with $(n+1)^p$ uniformly spaced points, interlaced by testing data $V$ with $n^p$ uniformly spaced points.

\subsection{Pointwise error of non-local first derivatives} \label{app:error_pointwisederivative_p}
We will now investigate the error in the non-local derivatives relative to differential derivatives, written as
\begin{align}
	\difference[1]{u(\tilde{x}_i)}{x^{\mu}} =&~ \derivative[1]{u(\tilde{x}_i)}{x^{\mu}} + \varepsilon_1^{\mu}(\tilde{x}_i).
\end{align}
When using the polynomial weight functions, the non-local first derivatives for $(n+1)^p$ possible Taylor series base points $\tilde{x}^{\mu}$, have the form of
\begin{align}
	\difference[1]{u(\tilde{x}_i)}{x^{\mu}} = \frac{C_{\epsilon}}{(n+1)^p-1}\sum_{j \in \tilde{V}\setminus \{i\}}\frac{(u(\tilde{x}_{j})-u(\tilde{x}_{i})) (\tilde{x}_{j}^{\mu}-\tilde{x}_{i}^{\mu})}{\left[\sum_{\mu}\abs{x^{\mu}_j-x^{\mu}_i}^2\right]^{\frac{2+\epsilon}{2}}}.
\end{align}
Here the polynomial form of the function, and the uniform spacing of the data $\tilde{x}_{j}-\tilde{x}_{i} = 2(j-i)h$ means the changes in $u$ can be written as
\begin{align}
	u(\tilde{x}_{j})-u(\tilde{x}_{i}) =&~ \sum_{l=0}^{K}\sum_{\mu}\alpha_{l\mu} (\tilde{x}_j^{\mu^l}-\tilde{x}_i^{\mu^l}), \\
	=&~ \sum_{l=0}^{K}\sum_{\mu}\alpha_{l\mu} (2h)^{l} (j^{\mu^l}-i^{\mu^l}).
\end{align}
and therefore the non-local first derivatives of polynomial functions are
\begin{align}
	\difference[1]{u(\tilde{x}_i)}{x^{\mu}} = \frac{C_{\epsilon}}{(n+1)^p-1}\sum_{l=0}^{K}\sum_{\nu}\alpha_{l\nu}(2h)^{l-1-\epsilon} \sum_{j \in \tilde{V}\setminus \{i\}}\frac{(j^{\nu^l}-i^{\nu^l})(j^{\mu}-i^{\mu})}{\left[\sum_{\theta}\abs{j^{\theta}-i^{\theta}}^2\right]^{\frac{2+\epsilon}{2}}},\\
	\intertext{and in the case $\epsilon=0$,}
	\difference[1]{u(\tilde{x}_i)}{x^{\mu}} = \frac{p}{(n+1)^p-1}\sum_{l=0}^{K}\sum_{\nu}\alpha_{l\nu}(2h)^{l-1} \sum_{j \in \tilde{V}\setminus \{i\}}\frac{(j^{\nu^l}-i^{\nu^l})(j^{\mu}-i^{\mu})}{\sum_{\theta}\abs{j^{\theta}-i^{\theta}}^2}.
\end{align}
We now can define the error for this multi dimensional polynomial as
\begin{align}
	\varepsilon_1^{\mu}(\tilde{x}_i) =&~ \varepsilon_1^{\mu}[h,\epsilon,\alpha](\tilde{x}_i), \\
	\intertext{and given the polynomial form of the function, and the subsequent form of the differential first derivatives in \cref{eq:polyderivatives_p}}
	\varepsilon_1^{\mu}[h,\epsilon,\alpha](\tilde{x}_i)  =&~ \sum_{l=0}^{K}\sum_{\nu} \alpha_{l\nu}\varepsilon_{1_{l}}^{\mu\nu}[h,\epsilon](\tilde{x}_i) , \label{eq:errorfirstorder_p}
\end{align}
where the errors associated with each term in the polynomial function are
\begin{align}
	\varepsilon_{1_{l}}^{\mu\nu}[h,\epsilon](\tilde{x}_i) =&~ \left[\left[\frac{C_{\epsilon}}{((n+1)^p-1)(2h)^{\epsilon}}\sum_{j \in \tilde{V}\setminus \{i\}}\frac{(j^{\nu^l}-i^{\nu^l})(j^{\mu}-i^{\mu})}{\left[\sum_{\theta}\abs{j^{\theta}-i^{\theta}}^2\right]^{\frac{2+\epsilon}{2}}}\right] \right.  \\
	&~\quad\quad -\left.\left[\frac{l!}{(l-1)!}i^{\mu^{l-1}}\mathbbm{1}_{l\geq1}\delta^{\mu\nu} \right]\right] (2h)^{l-1}, \nonumber \\
	\intertext{and for $\epsilon = 0$, the error takes the simplifed form of } 
	\varepsilon_{1_{l}}^{\mu\nu}[h,\epsilon=0](\tilde{x}_i) =&~ \left[\left[\frac{p}{((n+1)^p-1)h^{\epsilon}}\sum_{j \in \tilde{V}\setminus \{i\}}\frac{(j^{\nu^l}-i^{\nu^l})(j^{\mu}-i^{\mu})}{\sum_{\theta}\abs{j^{\theta}-i^{\theta}}^2}\right] \right. \\
	&~\quad\quad -\left. \left[\frac{l!}{(l-1)!}i^{\mu^{l-1}}\mathbbm{1}_{l\geq1}\delta^{\mu\nu} \right]\right] (2h)^{l-1}.	\nonumber
\end{align}
We can see that although the polynomial function has no coupling between different dimensions, the non-local derivative does impose that the dimensions are coupled with non-zero $\varepsilon_{1_{l}}^{\mu\nu}(\tilde{x}_i)$ for $\mu \neq \nu$. This potentially induces additional error from the differential derivatives which have no coupling across dimensions.

\subsection{Quadratic function modeled by linear Taylor series} \label{app:error_example_p}
\begin{figure}[hpt]
\centering
\includegraphics[width=\textwidth]{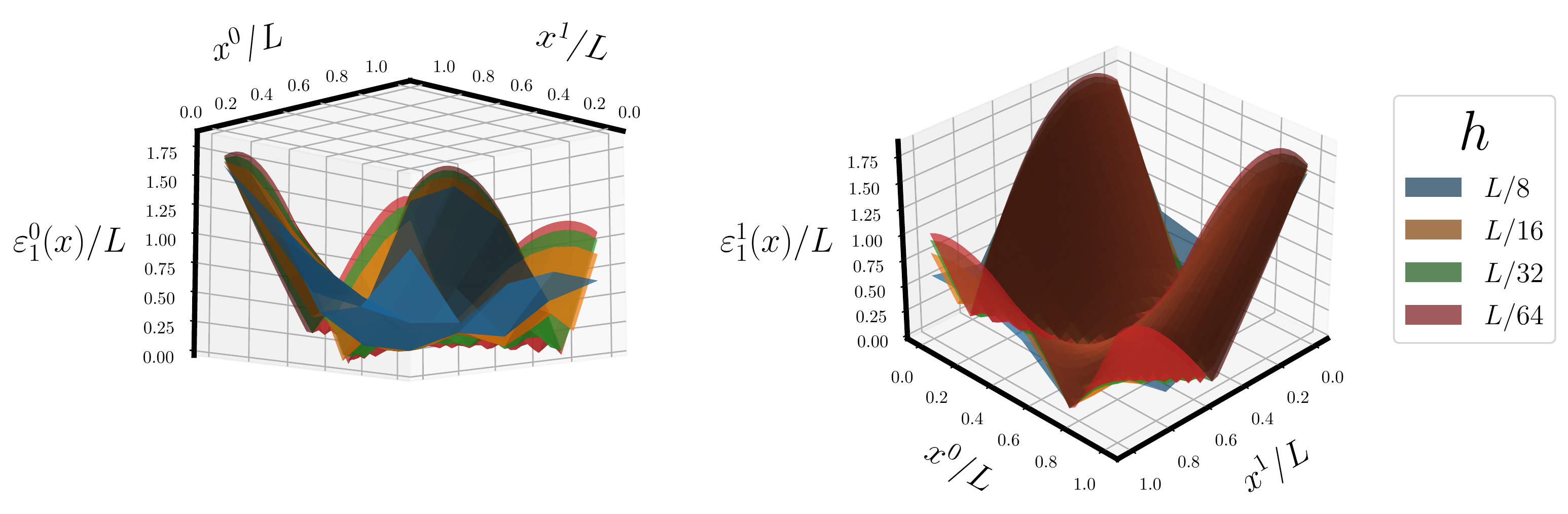}
\caption{Pointwise error of non-local first derivatives for $u(x) = \alpha_0 + \alpha_{10} x^{0} + \alpha_{20} x^{0^2} + \alpha_{11} x^{1} + \alpha_{21} x^{1^2}$, $\alpha_{l\mu} = 1$ and $\epsilon = 0$. $n+1$ points along an interval of length $L$ in each dimension are shown, with different spacings $h$.}
\label{fig:firstderivativeerror_p}

\end{figure}

\begin{figure}[hpt]
\centering
\includegraphics[width=0.7\textwidth]{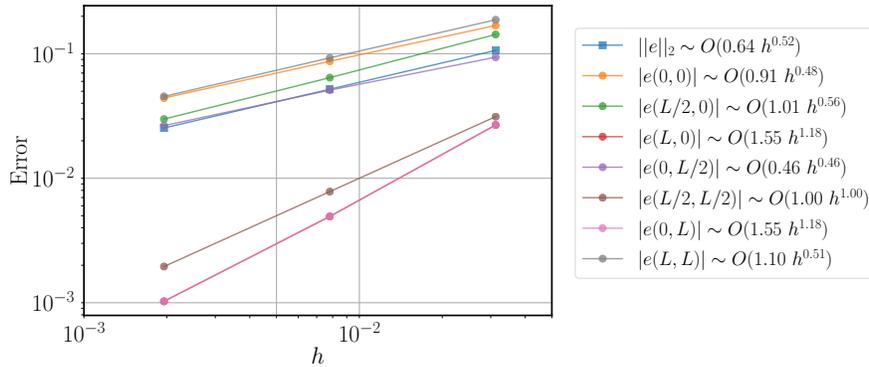}
\caption{Global and pointwise error scaling for linear Taylor series model $u(x) = \alpha_0 + \alpha_{10} x^{0} + \alpha_{20} x^{0^2} + \alpha_{11} x^{1} + \alpha_{21} x^{1^2}$, $\alpha_{l\mu} = 1$ and $\epsilon = 0$. Leading order scaling with $h$ fits are shown in the legend.}
\label{fig:totalerror_p}
\end{figure}

As an example, consider a $K=2$ quadratic function modeled by a $k=1$ linear Taylor series with uniformly spaced data. The error in the first derivative is shown in \cref{fig:firstderivativeerror_p} for $p=2$ dimensions. Similar to the one dimensional case, there is a persistent, spatially dependent error in the non-local first derivatives, that is peaked at the edges of the domain, and minimum at the center of the domain. The errors are not orders of magnitude larger than the length scale $L$, however are consistent and non-zero only for certain conditions on $\alpha$.
The total error and pointwise errors at $\tilde{x} = \{(0,0),~(L/2,0),~(L,0),~(0,L/2),~(L/2,L/2),~(0,L),~(L,L)\}$ for $p=2$ dimensions are shown in \cref{fig:totalerror_p}, for $\alpha_{l\mu} = 1$, and $\epsilon = 0$. Here, the total total error and pointwise error is shown to scale as 
\begin{align}
\norm{e}_2 \sim O(h^{\frac{1}{2}}) ~\longrightarrow~ Q[K=2,k=1,p=2]= \frac{1}{2} ,
\end{align}
with some regions at the edge of the domain scaling as $O(h)$. These regions of super-convergence are likely related to the error in the derivatives being minimal near those regions, particularly at the midpoint $(L/2,L/2)$, making the modified Taylor series more similar to the differential Taylor series. The error scaling has decreased for this higher dimensional function to be sub-linear, and is now more than an order less convergent than a differential Taylor series. This is likely due to the added complexity in the weight function in higher dimensions, and its support over the entire graph. 

\noindent It is conjectured that the error will at least be 
\begin{align}
	\norm{e}_2 \lessapprox&~ O(h^{Q[K,k,p]}) \\
	\intertext{such that} 
	Q[K,k,p] \leq&~ \frac{1}{2}.
\end{align}
for even higher dimensional $p>2$ problems, suggesting alternate weight function definitions should be sought for improved convergence of this method. Higher dimensional non-local derivatives have likely more regions of non-zero error, particularly for functions with coupled dependencies on the spatial dimensions. Weight functions that decay sufficiently and have finite extent across nearest neighbors in the graph, as opposed to having support over the entire graph, will contribute to the non-local derivatives behaving more similarly to a finite difference stencil approach. The effect of these weight functions on the resulting Taylor series model, and whether this contributes to $\gamma_l(\tilde{x}) \approx 1$, and the error scaling more similarly to other interpolation methods will be studied in future communications. \\

\end{appendices}

\end{document}

%% file: figures/model_vol_OLS_Ridge17.tex
\begin{align}
\frac{{\delta}^{} {\bar{\varphi}}_{\textrm{OLS}}}{{\delta} {{{t}}_{}}^{}} =&~ {\gamma}^{{{\bar{E}}_{11}}^{}{{\bar{E}}_{12}}^{}{{\bar{l}}_{1}}^{}}{{\bar{E}}_{11}}^{} {{\bar{E}}_{12}}^{} {{\bar{l}}_{1}}^{}~+~{\gamma}^{{{\bar{E}}_{12}}^{}{{\bar{\varphi}}_{}}^{}{{\bar{l}}_{1}}^{}}{{\bar{E}}_{12}}^{} {{\bar{\varphi}}_{}}^{} {{\bar{l}}_{1}}^{} \label{eq:vol_model_OLS} \\
+&~ {\gamma}^{{{\bar{E}}_{11}}^{}{{\bar{E}}_{22}}^{}{{\bar{l}}_{2}}^{}}{{\bar{E}}_{11}}^{} {{\bar{E}}_{22}}^{} {{\bar{l}}_{2}}^{}~+~{\gamma}^{{{\bar{E}}_{11}}^{}{{\bar{E}}_{22}}^{}}{{\bar{E}}_{11}}^{} {{\bar{E}}_{22}}^{} \nonumber \\ 
+&~  {\gamma}^{{{\bar{E}}_{22}}^{2}}{{\bar{E}}_{22}}^{2}~+~{\gamma}^{{{\bar{E}}_{22}}^{2}{{\bar{l}}_{2}}^{}}{{\bar{E}}_{22}}^{2} {{\bar{l}}_{2}}^{} \nonumber \\ 
+&~  {\gamma}^{{{\bar{\varphi}}_{}}^{}{{\bar{l}}_{1}}^{}}{{\bar{\varphi}}_{}}^{} {{\bar{l}}_{1}}^{}~+~{\gamma}^{{{\bar{\varphi}}_{}}^{2}{{\bar{l}}_{1}}^{}}{{\bar{\varphi}}_{}}^{2} {{\bar{l}}_{1}}^{} \nonumber \\ 
+&~  {\gamma}^{{{\bar{E}}_{11}}^{2}{{\bar{l}}_{1}}^{}}{{\bar{E}}_{11}}^{2} {{\bar{l}}_{1}}^{}~+~{\gamma}^{{{\bar{E}}_{11}}^{}{{\bar{l}}_{1}}^{}}{{\bar{E}}_{11}}^{} {{\bar{l}}_{1}}^{}, \nonumber \\ \nonumber \\
\frac{{\delta}^{} {\bar{\varphi}}_{\textrm{Ridge}}}{{\delta} {{{t}}_{}}^{}} =&~ {\gamma}^{{{\bar{E}}_{11}}^{}{{\bar{E}}_{12}}^{}}{{\bar{E}}_{11}}^{} {{\bar{E}}_{12}}^{}~+~{\gamma}^{{{\bar{E}}_{12}}^{}{{\bar{l}}_{1}}^{}}{{\bar{E}}_{12}}^{} {{\bar{l}}_{1}}^{} \label{eq:vol_model_Ridge17} \\
+&~ {\gamma}^{{{\bar{E}}_{11}}^{2}}{{\bar{E}}_{11}}^{2}~+~{\gamma}^{{{\bar{E}}_{11}}^{2}{{\bar{l}}_{2}}^{}}{{\bar{E}}_{11}}^{2} {{\bar{l}}_{2}}^{} \nonumber \\ 
+&~  {\gamma}^{{{\bar{E}}_{22}}^{}{{\bar{l}}_{2}}^{}}{{\bar{E}}_{22}}^{} {{\bar{l}}_{2}}^{}~+~{\gamma}^{{{\bar{E}}_{22}}^{2}}{{\bar{E}}_{22}}^{2} \nonumber \\ 
+&~  {\gamma}^{{{\bar{E}}_{11}}^{}}{{\bar{E}}_{11}}^{}~+~{\gamma}^{{{\bar{E}}_{12}}^{}{{\bar{E}}_{22}}^{}{{\bar{l}}_{1}}^{}}{{\bar{E}}_{12}}^{} {{\bar{E}}_{22}}^{} {{\bar{l}}_{1}}^{} \nonumber \\ 
+&~ {\gamma}^{{{\bar{\varphi}}_{}}^{2}{{\bar{l}}_{1}}^{}}{{\bar{\varphi}}_{}}^{2} {{\bar{l}}_{1}}^{}~+~ {\gamma}^{{{\bar{E}}_{11}}^{}{{\bar{E}}_{22}}^{}{{\bar{l}}_{1}}^{}}{{\bar{E}}_{11}}^{} {{\bar{E}}_{22}}^{} {{\bar{l}}_{1}}^{}. \nonumber
\end{align}

%% file: figures/model_psi_OLS_Ridge17.tex
\begin{align}
{\Psi}_{\textrm{OLS}} =&~ {{\Psi}_{}}_{}({{\mathbf{\bar{E}}}}_{0},{{{\bar{\varphi}}}}_{0}) + {\gamma}^{{\bar{\varphi}}_{}{\bar{\varphi}}_{}}\frac{1}{2!}\frac{{\delta}^{2} {\Psi}_{}({{\mathbf{\bar{E}}}}_{0},{{{\bar{\varphi}}}}_{0})}{{\delta {\bar{\varphi}}_{}}^{2}}{\Delta {\bar{\varphi}}_{}}^{2} \label{eq:psi_model_OLS}\\
+&~ {\gamma}^{{\bar{E}}_{12}{\bar{\varphi}}_{}}\frac{1}{2!}\frac{{\delta}^{2} {\Psi}_{}({{\mathbf{\bar{E}}}}_{0},{{{\bar{\varphi}}}}_{0})}{\delta {\bar{E}}_{12} \delta {\bar{\varphi}}_{}}{\Delta {\bar{E}}_{12}}^{} {\Delta {\bar{\varphi}}_{}}^{} + {\gamma}^{{\bar{E}}_{12}{\bar{\varphi}}_{}{\bar{\varphi}}_{}}\frac{1}{3!}\frac{{\delta}^{3} {\Psi}_{}({{\mathbf{\bar{E}}}}_{0},{{{\bar{\varphi}}}}_{0})}{\delta {\bar{E}}_{12} \delta {\bar{\varphi}}_{} \delta {\bar{\varphi}}_{}}{\Delta {\bar{E}}_{12}}^{} {\Delta {\bar{\varphi}}_{}}^{2} \nonumber \\
+&~ {\gamma}^{{\bar{E}}_{11}{\bar{E}}_{22}}\frac{1}{2!}\frac{{\delta}^{2} {\Psi}_{}({{\mathbf{\bar{E}}}}_{0},{{{\bar{\varphi}}}}_{0})}{\delta {\bar{E}}_{11} \delta {\bar{E}}_{22}}{\Delta {\bar{E}}_{11}}^{} {\Delta {\bar{E}}_{22}}^{} + {\gamma}^{{\bar{\varphi}}_{}{\bar{E}}_{22}{\bar{E}}_{11}}\frac{1}{3!}\frac{{\delta}^{3} {\Psi}_{}({{\mathbf{\bar{E}}}}_{0},{{{\bar{\varphi}}}}_{0})}{\delta {\bar{\varphi}}_{} \delta {\bar{E}}_{22} \delta {\bar{E}}_{11}}{\Delta {\bar{\varphi}}_{}}^{} {\Delta {\bar{E}}_{22}}^{} {\Delta {\bar{E}}_{11}}^{} \nonumber \\
+&~ {\gamma}^{{\bar{E}}_{11}{\bar{E}}_{22}{\bar{\varphi}}_{}{\bar{\varphi}}_{}}\frac{1}{4!}\frac{{\delta}^{4} {\Psi}_{}({{\mathbf{\bar{E}}}}_{0},{{{\bar{\varphi}}}}_{0})}{\delta {\bar{E}}_{11} \delta {\bar{E}}_{22} \delta {\bar{\varphi}}_{} \delta {\bar{\varphi}}_{}}{\Delta {\bar{E}}_{11}}^{} {\Delta {\bar{E}}_{22}}^{} {\Delta {\bar{\varphi}}_{}}^{2} + {\gamma}^{{\bar{E}}_{22}{\bar{\varphi}}_{}}\frac{1}{2!}\frac{{\delta}^{2} {\Psi}_{}({{\mathbf{\bar{E}}}}_{0},{{{\bar{\varphi}}}}_{0})}{\delta {\bar{E}}_{22} \delta {\bar{\varphi}}_{}}{\Delta {\bar{E}}_{22}}^{} {\Delta {\bar{\varphi}}_{}}^{} \nonumber \\
+&~ {\gamma}^{{\bar{E}}_{11}{\bar{\varphi}}_{}}\frac{1}{2!}\frac{{\delta}^{2} {\Psi}_{}({{\mathbf{\bar{E}}}}_{0},{{{\bar{\varphi}}}}_{0})}{\delta {\bar{E}}_{11} \delta {\bar{\varphi}}_{}}{\Delta {\bar{E}}_{11}}^{} {\Delta {\bar{\varphi}}_{}}^{} + {\gamma}^{{\bar{\varphi}}_{}{\bar{\varphi}}_{}{\bar{E}}_{22}}\frac{1}{3!}\frac{{\delta}^{3} {\Psi}_{}({{\mathbf{\bar{E}}}}_{0},{{{\bar{\varphi}}}}_{0})}{\delta {\bar{\varphi}}_{} \delta {\bar{\varphi}}_{} \delta {\bar{E}}_{22}}{\Delta {\bar{\varphi}}_{}}^{2} {\Delta {\bar{E}}_{22}}^{}, \nonumber \\ \nonumber \\
{\Psi}_{\textrm{Ridge}} =&~ {{\Psi}_{}}_{}({{\mathbf{\bar{E}}}}_{0},{{{\bar{\varphi}}}}_{0}) + {\gamma}^{{\bar{\varphi}}_{}{\bar{\varphi}}_{}}\frac{1}{2!}\frac{{\delta}^{2} {\Psi}_{}({{\mathbf{\bar{E}}}}_{0},{{{\bar{\varphi}}}}_{0})}{{\delta {\bar{\varphi}}_{}}^{2}}{\Delta {\bar{\varphi}}_{}}^{2} \label{eq:psi_model_Ridge17} \\
+&~ {\gamma}^{{\bar{\varphi}}_{}{\bar{\varphi}}_{}{\bar{\varphi}}_{}{\bar{\varphi}}_{}}\frac{1}{4!}\frac{{\delta}^{4} {\Psi}_{}({{\mathbf{\bar{E}}}}_{0},{{{\bar{\varphi}}}}_{0})}{{\delta {\bar{\varphi}}_{}}^{4}}{\Delta {\bar{\varphi}}_{}}^{4} + {\gamma}^{{\bar{\varphi}}_{}{\bar{E}}_{22}{\bar{E}}_{11}}\frac{1}{3!}\frac{{\delta}^{3} {\Psi}_{}({{\mathbf{\bar{E}}}}_{0},{{{\bar{\varphi}}}}_{0})}{\delta {\bar{\varphi}}_{} \delta {\bar{E}}_{22} \delta {\bar{E}}_{11}}{\Delta {\bar{\varphi}}_{}}^{} {\Delta {\bar{E}}_{22}}^{} {\Delta {\bar{E}}_{11}}^{} \nonumber \\
+&~ {\gamma}^{{\bar{\varphi}}_{}{\bar{E}}_{11}{\bar{E}}_{11}{\bar{E}}_{11}}\frac{1}{4!}\frac{{\delta}^{4} {\Psi}_{}({{\mathbf{\bar{E}}}}_{0},{{{\bar{\varphi}}}}_{0})}{\delta {\bar{\varphi}}_{} \delta {\bar{E}}_{11} \delta {\bar{E}}_{11} \delta {\bar{E}}_{11}}{\Delta {\bar{\varphi}}_{}}^{} {\Delta {\bar{E}}_{11}}^{3} + {\gamma}^{{\bar{E}}_{11}{\bar{E}}_{22}}\frac{1}{2!}\frac{{\delta}^{2} {\Psi}_{}({{\mathbf{\bar{E}}}}_{0},{{{\bar{\varphi}}}}_{0})}{\delta {\bar{E}}_{11} \delta {\bar{E}}_{22}}{\Delta {\bar{E}}_{11}}^{} {\Delta {\bar{E}}_{22}}^{} \nonumber \\
+&~ {\gamma}^{{\bar{E}}_{12}}\frac{{\delta}^{} {\Psi}_{}({{\mathbf{\bar{E}}}}_{0},{{{\bar{\varphi}}}}_{0})}{{\delta {\bar{E}}_{12}}^{}}{\Delta {\bar{E}}_{12}}^{} + {\gamma}^{{\bar{E}}_{11}{\bar{\varphi}}_{}}\frac{1}{2!}\frac{{\delta}^{2} {\Psi}_{}({{\mathbf{\bar{E}}}}_{0},{{{\bar{\varphi}}}}_{0})}{\delta {\bar{E}}_{11} \delta {\bar{\varphi}}_{}}{\Delta {\bar{E}}_{11}}^{} {\Delta {\bar{\varphi}}_{}}^{} \nonumber \\
+&~ {\gamma}^{{\bar{E}}_{22}{\bar{\varphi}}_{}}\frac{1}{2!}\frac{{\delta}^{2} {\Psi}_{}({{\mathbf{\bar{E}}}}_{0},{{{\bar{\varphi}}}}_{0})}{\delta {\bar{E}}_{22} \delta {\bar{\varphi}}_{}}{\Delta {\bar{E}}_{22}}^{} {\Delta {\bar{\varphi}}_{}}^{} + {\gamma}^{{\bar{\varphi}}_{}{\bar{\varphi}}_{}{\bar{E}}_{22}}\frac{1}{3!}\frac{{\delta}^{3} {\Psi}_{}({{\mathbf{\bar{E}}}}_{0},{{{\bar{\varphi}}}}_{0})}{\delta {\bar{\varphi}}_{} \delta {\bar{\varphi}}_{} \delta {\bar{E}}_{22}}{\Delta {\bar{\varphi}}_{}}^{2} {\Delta {\bar{E}}_{22}}^{}. \nonumber 
\end{align}